%% file: main.tex
\documentclass[A4paper]{article}

\newcommand{\footremember}[2]{%
    \footnote{#2}
    \newcounter{#1}
    \setcounter{#1}{\value{footnote}}%
}
\newcommand{\footrecall}[1]{%
    \footnotemark[\value{#1}]%
} 

\usepackage{verbatim}

\usepackage{graphicx}

\usepackage{amsmath}
\usepackage{amssymb}

\usepackage{bm}
\usepackage{mathtools}

\usepackage{bm}
\newcommand{\LL}{\mathbf{L}}
\newcommand{\WW}{\mathbf{W}}
\newcommand{\N}{\mathbb{N}}
\newcommand{\R}{\mathbb{R}}
\newcommand{\K}{\bm{K}}
\newcommand{\D}{\bm{D}}
\newcommand{\PPi}{\bm{\Pi}}
\newcommand{\V}{\bm{V}}
\newcommand{\x}{\bm{x}}
\newcommand{\uu}{\bm{u}}
\newcommand{\g}{\bm{g}}
\newcommand{\bbeta}{\bm{\beta}}
\newcommand{\vv}{\bm{v}}
\newcommand{\nn}{\bm{n}}
\newcommand{\pp}{\bm{p}}
\newcommand{\mm}{\bm{m}}
\newcommand{\qq}{\bm{q}}
\newcommand{\bb}{\bm{b}}
\newcommand{\vvarphi}{\bm{\varphi}}
\newcommand{\Pk}{\mathbb{P}_k\left(E\right)}
\newcommand{\Pkk}{\left(\mathbb{P}_k\left(E\right)\right)^2}
\newcommand{\ggVk}{{^{\g}\mathbf{V}}^k}
\newcommand{\ppVk}{{^{\pp}\mathbf{V}}^k}

\newcommand{\diver}{\operatorname{div}}
\newcommand{\rot}{\operatorname{rot}}
\newcommand{\dof}{\operatorname{dof}}

\newcommand{\ppHk}{{^{\pp}\mathbf{H}}^k}
\newcommand{\mmHk}{{^{\mm}\mathbf{H}}^k}
\newcommand{\qqHk}{{^{\qq}\mathbf{H}}^k}

\newcommand{\oneG}{{^{\bm{1}}\mathbf{G}}^k}
\newcommand{\twoG}{{^{\bm{2}}\mathbf{G}}^k}
\newcommand{\fourG}{{^{\bm{3}}\mathbf{G}}^k}

\newcommand{\twoGset}{{^{\bm{2}}\mathcal{G}}^k}

\newcommand{\oneTnablak}{{^{\bm{1}}\mathbf{T}}^{\nabla,k}}
\newcommand{\twoTnablak}{{^{\bm{2}}\mathbf{T}}^{\nabla,k}}
\newcommand{\fourTnablak}{{^{\bm{3}}\mathbf{T}}^{\nabla,k}}

\newcommand{\oneTbotk}{{^{\bm{1}}\mathbf{T}}^{\bot,k}}
\newcommand{\twoTbotk}{{^{\bm{2}}\mathbf{T}}^{\bot,k}}
\newcommand{\fourTbotk}{{^{\bm{3}}\mathbf{T}}^{\bot,k}}

\newcommand{\twognablak}{{^{\bm{2}}\g}^{\nabla,k}}
\newcommand{\fourgnablak}{{^{\bm{3}}\g}^{\nabla,k}}

\newcommand{\twogbotk}{{^{\bm{2}}\g}^{\bot,k}}
\newcommand{\fourgnablakm}{{^{\bm{3}}\g}^{\nabla,k-1}}

\usepackage{lineno}

\usepackage{subfigure}
\usepackage{graphicx}

\usepackage{booktabs}
\usepackage{multirow}
\usepackage{longtable}
\usepackage{pdflscape}  
\usepackage{lscape}
\usepackage{adjustbox}

\usepackage{comment}
\usepackage{hyperref}

\usepackage{amsthm}
\newtheorem{remark}{Remark}[section]
\newtheorem{theo}{Theorem}[section]

\begin{document}


\title{Orthogonal polynomial bases in the Mixed Virtual Element Method}

\author{S. Berrone\footremember{trailer}{Department of Applied Mathematics, Politecnico di Torino, Italy
  (stefano.berrone@polito.it, stefano.scialo@polito.it, gioana.teora@polito.it).}, S. Scial\`o\footrecall{trailer}{}, G. Teora\footrecall{trailer}{}}
\maketitle

\begin{abstract}
The use of orthonormal polynomial bases has been found to be efficient in preventing ill-conditioning of the system matrix in the primal formulation of Virtual Element Methods (VEM) for high values of polynomial degree and in presence of badly-shaped polygons. However, we show that using the natural extension of a orthogonal polynomial basis built for the primal formulation is not sufficient to cure ill-conditioning in the mixed case. Thus, in the present work, we introduce an orthogonal vector-polynomial basis which is built ad hoc for being used in the mixed formulation of VEM and which leads to very high-quality solution in each tested case. Furthermore, a numerical experiment related to simulations in Discrete Fracture Networks (DFN), which are often characterised by very badly-shaped elements, is proposed to validate our procedures.
\end{abstract}


\input{sections}


\bibliographystyle{IEEEtran}
\bibliography{biblio.bib}

\end{document}

%% file: sections.tex
\section{Introduction}

The Mixed Virtual Element Methods were introduced originally in \cite{basicMixed} for the Poisson problem in the two-dimensional case and then were extended to more general elliptic equations in \cite{Mixed}.
In the Mixed Virtual Element Space, two discrete spaces are introduced for approximating the pressure variable and the velocity field, respectively. The first space is a scalar-polynomial space, while a vector-polynomial basis is required to build the local projection matrices and for defining the internal degrees of freedom needed to obtain an approximation of the velocity field. It was observed that using the classical choice of scaled monomials in the definition of internal degrees of freedom in the primal Virtual Element construction \cite{basicMixed,basicvem,LBe14,MixedImplem}, the system matrix could become ill-conditioned in presence of \textit{badly-shaped} polygons \cite{Sbe17} (e.g. collapsing edges and bulks) and for high values of local polynomial degree \cite{mascotto}. 

The present work aims at defining new polynomial bases for the mixed VEM construction yielding well-conditioned local projection matrices also in presence of badly-shaped elements. More precisely, we propose two different approaches for building a vector-polynomial basis, which we briefly called ``Partial'' and ``Ortho'', respectively. The first one is the natural extension to the mixed case of the approach presented in \cite{mascotto} for the primal VEM, which allows us to build a vector-polynomial basis that is only partially $L^2$-orthonormalized. We show that the use of such basis is not sufficient to cure the ill-conditioning of the system matrix related to the mixed formulation of VEM in all circumstances \cite{DFNbook}, even if, in the primal setting, it reveals to be the best alternative.
Thus, we introduce the Ortho approach which aims to orthogonalize the gradients of a proper scalar-polynomial basis in order to obtain a full orthonormal vector-polynomial basis. We show that this approach leads to the best local and global performances, throughout different numerical experiments characterised by challenging geometries.

The outline of the present paper is the following. We define the model problem in Section \ref{sec:modelProblem} and its mixed VEM approximation in Section \ref{sec:mixedForm}. In Section \ref{sec:polynomialBasis}, we describe how to build the new polynomial bases, while in Section \ref{sec:implementation} we show an efficient implementation of the method, totally matrix-based. Finally, in Section \ref{sec:numericalExp}, we perform some numerical experiments that show the advantages of using the new polynomial bases.

Let us introduce some notations used throughout the paper. 
Given $k\in\N$, we use $(\cdot,\cdot)_{k,\omega}$ and $\|\cdot \|_{k,\omega}$ to indicate the inner product and the norm in the Sobolev space $H^k(\omega)$ on some open subset $\omega \subset \R^2$, respectively.
Furthermore, if $\vv = \begin{bmatrix} v_1, v_2
\end{bmatrix}^T$ and $\ \uu = \begin{bmatrix} u_1, u_2
\end{bmatrix}^T$ are vectors in $L^2(\omega) \times L^2(\omega)$, we define
\begin{equation}
    \left(\vv,\uu\right)_{0,\omega} = \int_{\omega} (v_1 u_1 + v_2 u_2), \quad \|\vv \|_{0,\omega} = \sqrt{\left(\vv,\vv\right)_{0,\omega}}.
\end{equation}

Let $\Omega \subset \mathbb{R}^2$ be a bounded convex polygonal domain with boundary $\Gamma$ and let $\nn_{\Gamma}$ be the outward unit normal vector to the boundary, then we define the functional spaces
\begin{equation}
    H(\diver; \Omega) = \big\{ \vv \in L^2(\Omega) \times L^2(\Omega) : \nabla \cdot \vv \in L^2(\Omega)\big\},
\end{equation}
\begin{equation}
    H_{0,\Gamma_N}(\diver; \Omega) = \big\{ \vv \in H(\diver, \Omega) : \vv \cdot \nn_{\Gamma} = 0 \text{ on } \Gamma_N \subseteq \Gamma \big\},
\end{equation}
\begin{equation}
    H(\rot; \Omega) = \{\vv\in L^2(\Omega) \times L^2(\Omega) : \rot \vv \in L^2(\Omega)\}.
\end{equation}
Furthermore, let $H^{-\frac{1}{2}}(\Gamma)$ be the dual space of the Sobolev space $H^{\frac{1}{2}}(\Gamma)$, the symbol $\langle \cdot, \cdot \rangle_{\pm \frac{1}{2},\Gamma}$ denotes the duality pairing between $H^{-\frac{1}{2}}(\Gamma)$ and $H^{\frac{1}{2}}(\Gamma)$.

\section{The continuous problem and the mixed variational formulation}\label{sec:modelProblem}

Let $\D(\x)$ be a symmetric uniformly positive definite tensor over $\Omega$, $\gamma$  a sufficiently smooth function $\Omega \to \R$ and $\bb$ a smooth vector valued function $\Omega \to \R^2$, then the following problem is considered:
\begin{equation}
    \begin{cases}
        \nabla \cdot \left(-\D \nabla p + \bb p \right) + \gamma p = f & \text{ in } \Omega,\\
        p = g_D & \text{ on } \Gamma_D,\\
        \left(-\D \nabla p + \bb p \right) \cdot \nn_{\Gamma_N} = g_N & \text{ on } \Gamma_N,
    \end{cases}
    \label{eq:continuousProblem}
\end{equation}
where $\Gamma_D$ and $\Gamma_N$ are the Dirichlet and the Neumann boundary, respectively, such that $ \Gamma_D \cup \Gamma_N = \Gamma $ and $\vert \Gamma_D \cap \Gamma_N \vert= 0$.

In order to introduce the mixed variational formulation, we define

\begin{equation}
    \K = \D ^{-1}, \quad \bbeta = \K \bb,
\end{equation}
and we re-write problem \eqref{eq:continuousProblem} as
\begin{equation}
    \begin{cases}
        \K \uu =  -\nabla p + \bbeta p & \text{ in } \Omega,\\
        \nabla \cdot \uu + \gamma p = f & \text{ in } \Omega,\\
        p = g_D & \text{ on } \Gamma_D,\\
        \uu \cdot \nn_{\Gamma_N} = g_N & \text{ on } \Gamma_N.
    \end{cases}
    \label{eq:mixedContinuousProblem}
\end{equation}

Thus, the mixed variational formulation of \eqref{eq:continuousProblem} reads:

\textit{Find } $\uu = \uu_0 + \uu_N,$ \textit{ with } $\uu_0 \in \V = H_{0,\Gamma_N}(\diver;\Omega), \textit{ and } p \in Q = L^2(\Omega)$ \textit{ such that }
\begin{equation}
    \begin{cases}
        \left(\K \uu, \vv \right)_{0,\Omega} -  \left(p, \nabla \cdot \vv \right)_{0,\Omega}  - \left(\bbeta  p, \vv \right)_{0,\Omega}  = - \langle g_D, \vv \cdot \nn_{\Gamma_D}\rangle_{\frac{1}{2},\Gamma_D} & \forall \vv \in \V \\
        \left(\nabla \cdot \uu, q \right)_{0,\Omega}  +\left( \gamma p,q \right)_{0,\Omega} = \left(f, q \right)_{0,\Omega} & \forall q \in  Q
    \end{cases}
    \label{eq:mixedVariationalFomrulation}
\end{equation}
where $\uu_N \in H(\diver,\Omega)$ is a chosen function that satisfies $\uu_N \cdot \nn_{\Gamma_N} = g_N$.

\section{The Mixed Virtual Element Method}\label{sec:mixedForm}

Let $\mathcal{T}_h$ be a decomposition of $\Omega$ into star-shaped polygons $E$.
We will denote by $\x_C$, $h_E$ and $\mathcal{E}_{E,h}$ the centroid, the diameter and the set of edges of $E$, respectively. We further set $N^{E,e} = \#\mathcal{E}_{E,h}$, and, as usual, we fix $h = \max_{E \in \mathcal{T}_h} h_E$.

Moreover, $\Pk$ is the set of all polynomials defined on $E$ of degree less or equal to $k\geq 0$ and $n_k = \dim \Pk = \frac{(k+1)(k+2)}{2}$. For simplicity, we fix $\mathbb{P}_{-1} = \{0\}$ and $n_{-1} = 0$.

A classical choice of a basis for the $\Pk$ space, that can be found in virtual element literature (see \cite{basicvem, LBe14, MixedImplem}), is the set of the scaled monomials defined as
\begin{equation}
    \mathcal{M}_k\left(E\right) = \left\{ \left(\frac{\x-\x_C}{h_E}\right)^{\bm{\alpha}}, \forall \bm{\alpha} = (\alpha_x,\alpha_y)\in \mathbb{N}^2 \text{ s.t. } 0 \leq \vert \bm{\alpha} \vert \leq k \right\}.
    \label{eq:scaledMonomial}
\end{equation}
For simplicity, we will make extensive use of the following function $\mathbb{N}^2 \to \mathbb{N}$

\begin{equation*}
    (0,0) \mapsto 1,\quad (1,0) \to 2, \quad (0,1) \mapsto 3,\quad (2,0) \mapsto 4, \dots
\end{equation*}

As in \cite{Mixed}, we introduce the (vector) polynomial space
\begin{equation}
    \mathcal{G}_k^{\nabla}(E) = \nabla \mathbb{P}_{k+1}(E) = \big\{\g^{\nabla,k}_{\alpha}\big\}_{\alpha = 1}^{n_k^{\nabla}}
    \label{eq:GNablak}
\end{equation}
and its complement $\mathcal{G}_k^{\bot}(E) = \big\{\g^{\bot,k}_{\alpha}\big\}_{\alpha = 1}^{n_k^{\bot}}$ in $\left(\mathbb{P}_{k}(E)\right)^2$, which satisfies
\begin{equation}
    \left(\mathbb{P}_{k}(E)\right)^2 = \mathcal{G}_k^{\nabla}(E) \bigoplus \mathcal{G}_k^{\bot}(E),
    \label{eq:complemCondition}
\end{equation}
where $\bigoplus$ is the direct sum operator, and
\begin{equation}
    \dim \left(\mathbb{P}_{k}(E)\right)^2 =2 n_k,
\end{equation}
\begin{equation}
   n^{\nabla}_k = \dim \mathcal{G}_k^{\nabla} = n_k + (k+1),
\end{equation}
\begin{equation}
   n_{k}^{\bot} = \dim \mathcal{G}_k^{\bot}(E) = n_k - (k+1).
\end{equation}

Now, following \cite{Mixed}, for any integer $k\geq 0$, we define the local mixed virtual element space for the velocity variable $\uu$ as

\begin{multline}
    \V_k(E) = \big\{ \vv_h \in H(\diver; E) \cap H(\rot; E) \text{ s.t. } \vv_h \cdot \nn_{e} \in \mathbb{P}_k(e)\ \forall e \in \mathcal{E}_{h,E}, \\
    \nabla \cdot \vv_h \in \mathbb{P}_k(E),\ \rot \vv_h \in \mathbb{P}_{k-1}(E) \big\}.
    \label{eq:mixedVemSpace}
\end{multline}
It is easy to see that $\left(\mathbb{P}_{k}(E)\right)^2 \subset  \V_k(E)$.

The following set of local degrees of freedom is unisolvent for $\V_k(E)$ (see \cite{MixedImplem,HdivHcurl}): given $\vv_h \in \V_k(E)$, 
\begin{itemize}
    \item \textbf{Edge dofs}: chosen $k+1$ Gauss quadrature points $\x^{e,Q}_i$ internal on each edge $e \in \mathcal{E}_{h,E}$:
    \begin{equation}
        \dof_i^e(\vv_h) =\left( \vv_h\cdot \nn_e \right) (\x^{e,Q}_i) \quad \forall i = 1,\dots, k+1.
        \label{eq:DOFedges}
    \end{equation}
    We note that this choice automatically ensures the continuity of the flux $\vv_h\cdot \nn$ across two adjacent elements.
    \item \textbf{Internal $\nabla$ dofs}:
    \begin{equation}
        \dof^{\nabla}_{\alpha}(\vv_h) =\frac{1}{\vert E \vert} \int_E \vv_h \cdot \g^{\nabla,k-1}_{\alpha}\quad \forall \alpha = 1,\dots, n_{k-1}^{\nabla}.
        \label{eq:DOFnabla}
    \end{equation}
    \item \textbf{Internal $\bot$ dofs}:
        \begin{equation}
        \dof^{\bot}_{\alpha}(\vv_h) =\frac{1}{\vert E \vert} \int_E \vv_h \cdot \g^{\bot,k}_{\alpha}\quad \forall \alpha = 1,\dots, n_{k}^{\bot}.
        \label{eq:DOFbigoplus}
    \end{equation}
\end{itemize}

Let it be $N^{\operatorname{dof}}_E= \dim\left(\V_k\left(E\right)\right) = N^{E,e} (k+1) +  n_{k-1}^{\nabla} + n_{k}^{\bot}$, we denote henceforth  the local Lagrangian mixed VE basis corresponding to the defined degrees of freedom:
\begin{equation}
 \left\{\vvarphi_i\right\}_{i=1}^{N^{\operatorname{dof}}_E} = \left\{\left\{\vvarphi_i^{e}\right\}_{i=1}^{ k+1}\right\}_{e \in \mathcal{E}_{h,E}} \cup \left\{\vvarphi_i^{\nabla}\right\}_{i=1}^{n_{k-1}^{\nabla}} \cup \left\{\vvarphi_i^{\bot}\right\}_{i=1}^{n_{k}^{\bot}},
\end{equation}
where the dofs numbering first counts the edge dofs, then the internal $\nabla$ dofs and lastly the internal $\bot$ dofs.

\bigskip

As in \cite{Mixed}, we define the local mixed virtual element space $Q_k(E)$ for the pressure variable $p$ as the space of polynomials $\Pk$, i.e. $Q_k(E) = \Pk$. In the next section, we will give more details about the choice of the local basis functions for the local space for the pressure variable.

\bigskip

Finally, we define the global mixed virtual element spaces for both velocity and pressure variables as
\begin{equation}
    \V_h= \{\vv_h \in H_{0,\Gamma_N}(\diver; \Omega) \text{ s.t. } \vv_{h|E} \in \V_k(E)\ \forall E \in \mathcal{T}_h\},
\end{equation}
\begin{equation}
    Q_h = \{q_h \in L^2(\Omega) \text{ s.t. } q_{h|E} \in Q_k(E)\ \forall E \in \mathcal{T}_h\}.
\end{equation}

\subsection{The discrete mixed variational formulation }

The $L^2(E)$-projection operator $\PPi^0_k : \V_h \to \Pkk$ is defined locally as
\begin{equation}
    \left(\PPi^0_k \vv_h, \pp_k\right)_{0,E} = \left(\vv_h, \pp_k\right)_{0,E} \ \forall \pp_k \in \Pkk \text{ and } \forall \vv_h \in \V_k(E).
    \label{eq:Pi0k}
\end{equation}
and, as shown in \cite{HdivHcurl}, the projection $\PPi^0_k \vv_h$ can be explicitly computed from the knowledge of the degrees of freedom (\ref{eq:DOFedges})-(\ref{eq:DOFbigoplus}) of $\vv_h \in \V_h$.

Now, the local discrete counterpart of 
\begin{equation}
    a(\uu,\vv) = \left(\K \uu, \vv \right)_{0,\Omega},\quad \forall \uu \in \V,\ \vv \in \V
\end{equation}
reads

\begin{align}
    a_h(\uu_h,\vv_h) &= \sum_{E \in \mathcal{T}_h}a_h^E\left(\uu_h, \vv_h \right) \\
    &= \sum_{E \in \mathcal{T}_h}\left((\K \PPi^0_k \uu_h, \PPi^0_k \vv_h)_{0,E} + S^E(\uu_h,\vv_h)\right),
    \label{eq:ahE}
\end{align}
where the stabilization term $S^E(\cdot,\cdot)$ is any symmetric and positive definite bilinear form that satisfies, $\forall \vv_h \in \V_h$

\begin{equation*}
    \alpha_{\ast} a_{| E}(\vv_h, \vv_h) \leq S^E(\vv_h,\vv_h) \leq \alpha^{\ast} a_{| E}(\vv_h, \vv_h)
\end{equation*}
for some constants $\alpha_{\ast},\alpha^{\ast}> 0$ that are depending on $\K$ but independent of $h$.
As in \cite{Mixed,MixedImplem}, we will choose

\begin{equation*}
    S^E(\uu_h,\vv_h) = \overline{\K} \vert E \vert \sum_{r=1}^{N^{\dof}_E} \dof_r\left((I-\PPi^0_k)\uu_h\right)\dof_r\left((I-\PPi^0_k)\vv_h\right),
\end{equation*}
where $\overline{\K}$ is the largest singular value of $\K$ on $E$.

Finally, the mixed VEM approximation of \eqref{eq:mixedVariationalFomrulation} is given by:

\textit{Find } $\uu_h = \uu_{0,h} + \uu_{N,h},$ \textit{ with } $\uu_{0,h} \in \V_h, \textit{ and } p_h \in Q_h$ \textit{ such that } $\forall \vv_h \in \V_h$ \textit{and} $\forall q_h \in  Q_h$:
\begin{equation}
    \begin{cases}
        a_h\left( \uu_h, \vv_h \right) -  \left(p_h, \nabla \cdot \vv_h \right)_{0,\Omega}  - \left(\bbeta  p_h, \PPi^0_k \vv_h \right)_{0,\Omega}  = - \langle g_D, \vv_h \cdot \nn_{\Gamma_D}\rangle_{\pm \frac{1}{2},\Gamma_D} \\
        \left(\nabla \cdot \uu_h, q_h \right)_{0,\Omega}  +\left( \gamma p_h,q_h \right)_{0,\Omega} = \left(f, q_h \right)_{0,\Omega}
    \end{cases}
    \label{eq:discreteMixedVariationalFomrulation}
\end{equation}
where $\uu_{N,h} \in \{\vv \in H(\diver,\Omega):  \vv \in \V_k(E), \forall E \in \mathcal{T}_h \}$ is a proper function that satisfies $\operatorname{dof}^e_i(\uu_{N,h}) = \operatorname{dof}^e_i(\uu_N)$, for each boundary edge $e \in \mathcal{E}_h$ and $\forall i = 1,\dots, k+1$.

The problem (\ref{eq:discreteMixedVariationalFomrulation}) has unique solution $(\uu_h,p_h) \in \V_h \times Q_h$ and, for $h$ sufficiently small, the following a priori error estimates holds true

\begin{equation}
    \| p -p_h\|_{0} = O( h^{k+1}),\quad 
    \| \uu - \uu_h\|_0 = O( h^{k+1}).
\end{equation}
Furthermore, the following superconvergence result holds true.

\begin{theo}[Superconvergence result]
Let $ p_h$ the solution to (\ref{eq:discreteMixedVariationalFomrulation}) and let $ p_I \in Q_h$ be the interpolant of $p$. Then, for $h$ sufficiently small, 
\begin{equation}
    \| p_I - p_h\|_0 = O(h^{k+2}).
\end{equation}

\end{theo}

\section{Polynomial basis}\label{sec:polynomialBasis}

We will now show some procedures for building polynomial bases for both $\Pk$ and $\Pkk$.

As mentioned in the previous section, the standard choice for the polynomial basis is the set of scaled monomials, defined in \eqref{eq:scaledMonomial} (see \cite{basicvem,LBe14,MixedImplem}). Orthogonal polynomial bases have already been introduced in the virtual element literature, since it has been proven that modifying the definition of internal moments by choosing an $L^2(E)$-orthonormal basis for $\Pk$ can largely improve the reliability of the virtual element method for higher order approximations and in presence of badly-shaped polygons (\cite{Sbe17,mascotto}).

An efficient procedure for building an $L^2(E)$-orthonormal polynomial basis for $\Pk$ is presented in \cite{mgs_basis}.
This procedure consists in the application of the  modified Gram-Schmidt (MGS) orthogonalization process to the monomial Vandermonde matrix $^m\mathbf{V}^{E,k} \in \R^{N^{E,Q} \times n_k}$ associated to a quadrature formula $\left(\{\mathbf{x}_i^{E,Q}\}, \{w_i^{E,Q}\} \right)_{i=1}^{N^{E,Q}}$ on a given element $E$ with $N^{E,Q}$ nodes. As suggested in \cite{mgs_basis,mgs_reortho}, the process must be applied twice in order to make the orthonomalization error $\| \mathbf{I} - {^{q}\mathbf{H}^{E,k}} \|$ independent of the condition number of $^m\mathbf{V}^{E,k}$ matrix, where ${^{q}\mathbf{H}^{E,k}} \in \R^{n_k \times n_k}$ is the mass matrix of the resulting $L^2(E)$-orthonormal polynomial basis $\{q^k_{\alpha}\}^{n_k}_{\alpha=1}$ for $\Pk$.

In particular, the whole procedure defines a matrix $\LL^{E,k} \in \R^{n_k \times n_k}$ on each element $E$ such that 
\begin{equation}
    ^{q}\mathbf{V}^{E,k} = {^{m}\mathbf{V}^{E,k}} \left(\LL^{E,k}\right)^T,
    \label{eq:MGS_basis}
\end{equation}
where $^{q}\mathbf{V}^{E,k}$ is the Vandermonde matrix associated to the new polynomial basis. More precisely, we first apply MGS process to the $^m\mathbf{V}^{E,k}$ matrix, i.e. we define an upper triangular matrix $\mathbf{R}_1^{E,k} \in \R^{n_k \times n_k}$ and an orthonormal matrix $\mathbf{Q}_1^{E,k}  \in \R^{N^{E,Q} \times n_k}$ such that

\begin{equation*}
    ^m\mathbf{V}^{E,k} = \mathbf{Q}_1^{E,k}  \mathbf{R}_1^{E,k} , 
\end{equation*}
and then we apply MGS to the $\mathbf{Q}_1^{E,k} $ matrix properly rescaled by quadrature weights in order to obtain an $L^2(E)$-orthonormal basis:

\begin{equation*}
    \left(\mathbf{W}^{E,Q}\right)^{1/2}\mathbf{Q}_1^{E,k}  = \mathbf{Q}_2^{E,k} \mathbf{R}_2^{E,k} , 
\end{equation*}
where $\mathbf{W}^{E,Q}\in \R^{N^{E,Q} \times N^{E,Q}}$ is the diagonal matrix whose diagonal entries are the quadrature weights.
Being $\LL^{E,k} = \left( \mathbf{R}_2^{E,k}  \mathbf{R}_1^{E,k}  \right)^{-T}$, we note that

\begin{align*}
    \mathbf{Q}_2^{E,k}  &= \left(\mathbf{W}^{E,Q}\right)^{1/2}\mathbf{Q}_1^{E,k}  \left(\mathbf{R}_2^{E,k}\right)^{-1}\\
    &= \left(\mathbf{W}^{E,Q}\right)^{1/2}\ {^{m}\mathbf{V}^{E,k}} \left(\LL^{E,k}\right)^T \\
    &= \left(\mathbf{W}^{E,Q}\right)^{1/2}\ {^{q}\mathbf{V}^{E,k}}
\end{align*}
which means that $\mathbf{Q}_2^{E,k} $ is the Vandermonde matrix ${^{q}\mathbf{V}^{E,k}}$ rescaled by the square root of the quadrature weights, 
thus,

\begin{equation}
        {^q\mathbf{H}^{E,k}} = \left({^{q}\mathbf{V}^{E,k}}\right)^T \mathbf{W}^{E,Q}\ {^{q}\mathbf{V}^{E,k}} = \left(\mathbf{Q}_2^{E,k} \right)^T \mathbf{Q}_2^{E,k} = \mathbf{I}.
\end{equation}

\begin{remark}
    Note that $^{q}\mathbf{V}^{E,k} = {^{m}\mathbf{V}^{E,k}} \left(\LL^{E,k}\right)^T$, and the columns of $\left(\LL^{E,k}\right)^T$ are the coefficients that provide each orthonormal polynomial as a linear combination of the monomials.
\end{remark}

\subsection{Polynomial bases for \texorpdfstring{$\Pkk$}{the polynomial space}}

For simplicity, in the following, we will drop the superscript $E$, when no ambiguity occurs. Furthermore, the left superscript will denote the underlying used polynomial basis for the space $\Pk$. In particular, we will use the symbols $m$ and $q$ to indicate the scaled monomial and the MGS basis, respectively. Instead, we will use the symbol $p$ to represent a generic polynomial basis.

Now, we introduce some auxiliary matrices that we will use in the following. Let ${^{p}}\mathbf{D}^{k+1,x}$, ${^{p}}\mathbf{D}^{k+1,y} \in \mathbb{R}^{ n_{k+1} \times n_k}$ be the matrices which collect the coefficients of the partial derivatives of $p^{k+1}_{\alpha}$, i.e.

\begin{equation}
    \frac{\partial p^{k+1}_{\alpha}}{\partial x} = \sum_{\beta = 1}^{n_k} {^{p}}\mathbf{D}^{k+1,x}_{\alpha \beta} p^{k}_{\beta}, \quad     \frac{\partial p^{k+1}_{\alpha}}{\partial y} = \sum_{\beta = 1}^{n_k} {^{p}}\mathbf{D}^{k+1,y}_{\alpha \beta} p^{k}_{\beta}, \quad \forall \alpha = 1,\dots,n_{k+1}.
    \label{eq:gradMatrix}
\end{equation}

Note that, if the MGS basis is used, matrices ${^q\mathbf{D}^{k+1,\ast}}$ can be derived from the (easily computable) monomial ones ${^m\mathbf{D}^{k+1,\ast}}$ as:

\begin{equation*}
     {^q\mathbf{D}^{k+1,\ast}} = \LL^{k+1} {^m\mathbf{D}^{k+1,\ast}}  (\LL^{k})^{-1} \quad \forall \ast \in \mathcal{D} =\{x,y\},
\end{equation*}
since, $\forall \alpha = 1,\dots,n_{k+1}$,

\begin{align*}
    \frac{\partial q_{\alpha}^{k+1}}{\partial \ast} &= \sum_{\beta = 1}^{n_{k+1}} \LL^{k+1}_{\alpha \beta} \frac{\partial m_{\beta}^{k+1}}{\partial \ast} = \sum_{\beta = 1}^{n_{k+1}} \sum_{\gamma = 1}^{n_{k}} \LL^{k+1}_{\alpha \beta}  {^m\mathbf{D}^{k+1,\ast}_{\beta \gamma}} m_{\gamma}^k \\
    &= \sum_{\beta = 1}^{n_{k+1}} \sum_{\gamma = 1}^{n_{k}} \sum_{s=1}^{n_k} \LL^{k+1}_{\alpha \beta}  {^m\mathbf{D}^{k+1,\ast}_{\beta \gamma}} (\LL^{k})^{-1}_{\gamma s} q_{s}^k.
\end{align*}

\begin{remark}
As highlighted in \cite{mgs_basis}, the modified Gram-Schmidt algorithm allows to obtain a hierarchical sequence of bases $\{\{q^k_{\alpha}\}_{\alpha=1}^{n_k}\}_{k \geq 0}$, i.e.
\begin{equation}
    \{q^k_{\alpha}\}_{\alpha=1}^{n_k} \subset \{q^{k+1}_{\alpha}\}_{\alpha=1}^{n_{k+1}}.
\end{equation}
As a consequence, we only need to compute $\mathbf{L}^{k+1}$ and then we can set
\begin{equation}
    \mathbf{L}^{k} =\mathbf{L}^{k+1}(1:n_k, 1:n_k),
\end{equation}
being $\mathbf{L}^{k+1}(1:n_k, 1:n_k)$ the matrix obtained from the first $n_k$ rows and columns of $\mathbf{L}^{k+1}$. 
Indeed, we define the Vandermonde matrices of both $k$ and $k+1$ orders with respect to the same quadrature formula.
\end{remark}

Starting from a generic polynomial basis for $\Pk$, as shown in \cite{MixedImplem}, an easily computable basis $\{\pp_I^k\}_{I=1}^{2n_k}$ for $\Pkk$ can be built as

\begin{equation}
    \pp_I^k =\begin{cases}
    \begin{bmatrix} p^k_I \\ 0\end{bmatrix} I=1,\dots,n_k\\ \\
    \begin{bmatrix} 0 \\ p^k_{I-n_k}\end{bmatrix} I=n_k+1,\dots,2n_k
    \end{cases}
\end{equation}

The Vandermonde matrix $\ppVk \in \mathbb{R}^{2N^{Q}\times 2n_k }$ associated to the $\{\pp_I^k\}_{I=1}^{2n_k}$ polynomial basis functions can be written as

\begin{equation}
    \ppVk =   \begin{bmatrix} {^{p}\mathbf{V}^{k}} & \mathbf{O} \in \mathbb{R}^{N^{Q} \times n_k}\\
    \mathbf{O} \in \mathbb{R}^{N^{Q} \times n_k} & {^{p}\mathbf{V}^{k}} \end{bmatrix},
\end{equation}
where
\begin{itemize}
    \item the top left ${^{p}\mathbf{V}^{k}}$ matrix contains the evaluations of the $x$-components of $\pp_I^k$ for all $I=1,\dots,n_k$;
    \item the top right $\mathbf{O}$ matrix contains the evaluations of the $y$-components of $\pp_I^k$ for all $I=1,\dots,n_k$, that are all zeros;
    \item the bottom left $\mathbf{O}$ matrix contains the evaluations of the $x$-components of $\pp_I^k$ for all $I=n_k+1,\dots,2n_k$, that are all zeros;
    \item the bottom right ${^{p}\mathbf{V}^{k}}$ matrix contains the evaluations of the $y$-components of $\pp_I^k$ for all $I=n_k+1,\dots,2n_k$.
\end{itemize}

Note that if $\{p^k_{\alpha}\}_{\alpha=1}^{n_k}$ is an $L^2(E)$-orthonormal polynomial basis for $\mathbb{P}_k(E)$, then 
$\{\pp^k_{I}\}_{I=1}^{2n_k}$ is an $L^2(E)$-orthonormal polynomial basis for $\Pkk$.

\bigskip

Now, functions belonging to $\mathcal{G}^{\nabla}_k$, as defined in \eqref{eq:GNablak}, can be written as
\begin{equation}
    \g^{\nabla,k}_{\alpha} = \nabla p^{k+1}_{\alpha + 1} = \sum_{I = 1}^{2 n_k} \mathbf{T}^{\nabla,k}_{\alpha I} \pp_I^k \quad \forall \alpha = 1,\dots, n_k^{\nabla},
\end{equation}
where $\mathbf{T}^{\nabla,k} \in \R^{n^{\nabla}_k \times 2n_k}$ is the coefficient matrix of gradients of the polynomial functions $ \{p^{k+1}_{\alpha}\}_{\alpha=2}^{n_{k+1}}$ with respect to the polynomial basis $\{\pp^k_{I}\}_{I=1}^{2n_k}$ of $\Pkk$ and the corresponding Vandermonde matrix is
\begin{equation}
    {^{\g}\mathbf{V}^{\nabla,k}} = \ppVk \left(\mathbf{T}^{\nabla,k}\right)^T.
\end{equation}
Based on definitions \eqref{eq:gradMatrix}, the $\mathbf{T}^{\nabla,k}$ matrix reads

\begin{equation*}
    \mathbf{T}^{\nabla,k} = \begin{bmatrix}
        {^{p}}\mathbf{D}^{k+1,x}(2:n_{k+1},:) & {^{p}}\mathbf{D}^{k+1,y}(2:n_{k+1},:)
    \end{bmatrix},
\end{equation*}
with ${^{p}}\mathbf{D}^{k+1,\ast}(2:n_{k+1},:)$ the sub-matrix of ${^{p}}\mathbf{D}^{k+1,\ast}$ obtained extracting rows from $2$ to $n_{k+1}$ and all columns.

Now, we can complete a basis $\mathcal{G}_k$ for $\Pkk$, by adding the set of functions $\mathcal{G}^{\bot}_k = \{\g^{\bot, k}_{\alpha}\}_{\alpha = 1}^{n^{\bot}_k}$ such that \eqref{eq:complemCondition} is satisfied. As suggested in \cite{MixedImplem}, $\g^{\bot,k}_{\alpha}$ function can be defined as
\begin{equation}
    \g^{\bot,k}_{\alpha}  = \sum_{I = 1}^{2 n_k} \mathbf{T}^{\bot,k}_{\alpha I} \pp_I^k \quad \forall \alpha = 1,\dots, n_k^{\bot},
\end{equation}
where $\mathbf{T}^{\bot,k} \in \R^{n^{\bot}_k \times 2n_k}$ is the matrix whose rows define an \textit{euclidean} orthonormal basis for the nullspace of $\mathbf{T}^{\nabla,k}$ matrix. Thus, by considering the Singular Value Decomposition of $\mathbf{T}^{\nabla,k} = \mathbf{U}\mathbf{\Sigma}\left(\mathbf{V}\right)^T$, we can define $\mathbf{T}^{\bot,k}$ as 
\begin{equation}
    \mathbf{T}^{\bot,k} = \left(\mathbf{V}(:,n^{\nabla}_k+1:2n_k)\right)^T,
\end{equation}
where $\mathbf{V}(:,n^{\nabla}_k+1:2n_k)$ is the submatrix of $\mathbf{V}$ made up of all its rows and of the columns running from the $(n_{k}^{\nabla}+1)$-th to the $2n_{k}$-th. As a consequence,
\begin{equation}
    \mathbf{T}^{\nabla,k} \left(\mathbf{T}^{\bot,k}\right)^T = \mathbf{0}.
    \label{eq:propertyTortho}
\end{equation}

The Vandermonde matrix $\ggVk \in \mathbb{R}^{2N^{Q} \times 2n_k}$ associated to the basis functions $\mathcal{G}^k = \{\g^{k}_I\}_{I=1}^{2n_k} = \{\g^{\nabla,k}_{\alpha}\}_{\alpha = 1}^{n_k^{\nabla}} \cup \{\g^{\bot,k}_{\beta}\}_{\beta = 1}^{n_k^{\bot}}$ reads

\begin{equation}
    \ggVk = \begin{bmatrix}
        {^{\g}\mathbf{V}^{\nabla,k}} & {^{\g}\mathbf{V}^{\bot,k}}
    \end{bmatrix} = \ppVk\ \begin{bmatrix}
    \left(\mathbf{T}^{\nabla,k}\right)^T & \left(\mathbf{T}^{\bot,k}\right)^T
    \end{bmatrix}.
\end{equation}

We define $\mathbf{G}^k \in \mathbb{R}^{2n_k \times 2n_k}$ as the mass matrix of the $\mathcal{G}^k$ basis, whose entries are given by
\begin{equation}
    \mathbf{G}_{IJ}^k = \int_E \g_I^k \cdot \g_J^k,\quad \forall I,J=1,\dots,2n_k.
\end{equation}
In matrix form, $\mathbf{G}^k $ reads

\begin{align}
    \mathbf{G}^k 
    &=  \nonumber \left(\ggVk\right)^T \begin{bmatrix}
    \mathbf{W}^{Q} & \mathbf{O}\\
    \mathbf{O} & \mathbf{W}^{Q}
    \end{bmatrix} \ggVk\\
    &= \nonumber \begin{bmatrix}
     \mathbf{T}^{\nabla,k}\ \ppHk\left(\mathbf{T}^{\nabla,k}\right)^T & \mathbf{T}^{\nabla,k}\ \ppHk \left(\mathbf{T}^{\bot,k}\right)^T \\
    \mathbf{T}^{\bot,k}\ \ppHk \left(\mathbf{T}^{\nabla,k}\right)^T & \mathbf{T}^{\bot,k}\ \ppHk \left(\mathbf{T}^{\bot,k}\right)^T
    \end{bmatrix} \\
    &= \begin{bmatrix}
    \mathbf{G}^{\nabla,\nabla} & \mathbf{G}^{\nabla,\bot} \\
    \mathbf{G}^{\bot,\nabla} & \mathbf{G}^{\bot,\bot}
    \end{bmatrix},
\end{align}
where 

\begin{equation*}
    \ppHk = \left(\ppVk\right)^T \begin{bmatrix}
    \mathbf{W}^{Q} & \mathbf{O}\\
    \mathbf{O} & \mathbf{W}^{Q}
    \end{bmatrix} \ppVk 
\end{equation*}
is the mass matrix related to $\{\pp_I\}_{I=1}^{2n_k}$ basis.

We note that
\begin{enumerate}
    \item if we choose the standard set of scaled monomials as the basis for $\Pk$, it is known that the corresponding $\oneG$ matrix, i.e.
    \begin{equation}
        \oneG=  \begin{bmatrix}
        \oneTnablak \ \mmHk \left(\oneTnablak\right)^T & \oneTnablak\ \mmHk \left(\oneTbotk\right)^T \\
        \oneTbotk\ \mmHk \left(\oneTnablak\right)^T & \oneTbotk \ \mmHk \left(\oneTbotk\right)^T
        \end{bmatrix}      
        \label{eq:monomialGMatrix}
    \end{equation}
    will be ill-conditioned for high polynomial degrees. Anyway, we want to highlight that, in the so-called \textit{monomial approach}, the rows of $\oneTbotk$ are orthonormal to each other and are orthogonal to the rows of $\oneTnablak$ with respect to the euclidean scalar product, by construction.
    
    \item On the other hand, if we choose the MGS basis, i.e. the $L^2(E)$-orthonormal polynomial basis for $\Pk$ introduced in \cite{mgs_basis}, the corresponding mass matrix $\qqHk$ will be the identity matrix $\mathbf{I}$. Thus, in infinite precision, $\twoG$ takes the form

    \begin{align}
        \twoG = \begin{bmatrix}
         \twoTnablak \ \mathbf{I} \left(\twoTnablak\right)^T & \twoTnablak\ \mathbf{I} \left(\twoTbotk\right)^T \\
        \twoTbotk\ \mathbf{I} \left(\twoTnablak\right)^T & \twoTbotk\ \mathbf{I} \left(\twoTbotk\right)^T
        \end{bmatrix} 
         = \begin{bmatrix}
         \twoTnablak \left(\twoTnablak \right)^T & \mathbf{O} \\
        \mathbf{O} & \mathbf{I}
        \end{bmatrix},
        \label{eq:partialOrthoGMatrix}
    \end{align}
    where the last equation is a consequence of property \eqref{eq:propertyTortho}.
     Thus, if an $L^2(E)$-orthonormal polynomial basis for $\Pk$ are used, $\twoGset$ will be partially $L^2(E)$-orthonormalized, since the $\twogbotk_{\alpha}$ functions are a set of $L^2(E)$-orthonormal functions and are $L^2(E)$-orthogonal to $\twognablak_{\alpha}$ functions, but the $\twognablak_{\alpha}$ are not naturally $L^2(E)$-orthogonal to each other. We denote this as \textit{partial-orthonormal approach}.
\end{enumerate}

\subsubsection{A full \texorpdfstring{$L^2(E)$}{L2E}-orthonormal approach}

In this section, we will show a procedure that allows to full $L^2(E)$-orthonormalize the $\mathcal{G}^k$ basis. For this purpose, we note from \eqref{eq:partialOrthoGMatrix} that in order to make $\mathcal{G}^{\nabla}_k$ an $L^2(E)$-orthonormal set of functions, it is sufficient to orthonormalize the rows of $\twoTnablak$ matrix with respect to the euclidean scalar product. 
In order to do this, we apply only once the modified Gram-Schmidt algorithm to $\left(\twoTnablak\right)^T$. More precisely, we factorize $\left(\twoTnablak\right)^T$ as

\begin{equation*}
    \left(\twoTnablak\right)^T  = \mathbf{Q}^{\nabla,k}\mathbf{R}^{\nabla,k},
\end{equation*}
and, then, we set
\begin{equation}
    \fourTnablak = (\mathbf{Q}^{\nabla,k})^T = \LL^{\nabla,k}\ \twoTnablak, \quad \fourTbotk = \tilde{\mathbf{V}}(:, n^{\nabla}_k+1:2n_k)^T,
    \label{eq:Tnabla_full}
\end{equation}
where $\fourTnablak = \tilde{\mathbf{U}}\tilde{\mathbf{\Sigma}}\tilde{\mathbf{V}}^T$ and $\LL^{\nabla,k} = \left(\mathbf{R}^{\nabla,k}\right)^{-T} \in \R^{n^{\nabla}_k\times n^{\nabla}_k}$.

By proceeding in this way, the $\fourG$ matrix will become the identity matrix. We will call this procedure as \textit{full-orthonormal approach}.

\begin{remark}\label{rem:gnabla_hierarchical}
It is worth mentioning that in order to orthonormalize the rows of $\twoTnablak$ we could resort to its computed Singular Value Decomposition and set

\begin{align*}
    \fourTnablak = \mathbf{V}(&:, 1:n^{\nabla}_k)^T = \left(\mathbf{\Sigma}(:,1:n^{\nabla}_k)\right)^{-1}\mathbf{U}^T\ \twoTnablak, \\
    \fourTbotk &= \mathbf{V}(:, n^{\nabla}_k+1:2n_k)^T = \twoTbotk .
\end{align*}
In this way, we obtain that $\fourG = \mathbf{I}$, but the SVD process is not hierarchical and, as a consequence, the set of functions $\{\fourgnablakm_{\alpha}\}_{\alpha=1}^{n^{\nabla}_{k-1}}$ used for defining the internal $\nabla$ dofs is not $L^2(E)$-orthonormalized for free. As highlighted in \cite{Sbe17,mascotto}, this type of change can improve only locally the condition number of elemental matrices, but does not ensure to improve the global performances of the method.

Instead, since the MGS is a hierarchical procedure, we obtain also that the functions
\begin{equation}
   \{\fourgnablakm_{\alpha}\}_{\alpha=1}^{n^{\nabla}_{k-1}} \subset \{\fourgnablak_{\alpha}\}_{\alpha=1}^{n^{\nabla}_k}
\end{equation}
used in the definition of internal $\nabla$ dofs are a set of $L^2(E)$-orthonormal functions.
\end{remark}

\section{Some implementation details}\label{sec:implementation}

In this section, we will show how to compute the local matrices needed for assembling the local system matrix related to the discrete problem \eqref{eq:discreteMixedVariationalFomrulation} with the aforementioned approaches, following a procedure similar to the one shown in \cite{MixedImplem}.

\subsection{Term \texorpdfstring{$\left(\nabla \cdot \uu_h, q_h\right)_{0,E}$}{which includes divergence}}\label{sec:Wmatrix}

By choosing $\uu_h = \vvarphi_i$ and $q_h = p^k_{\alpha}$, we define $\WW \in \R^{n_k \times N^{\dof}_E}$ as the matrix whose entries read, $\forall i =1,\dots,N^{\dof}_E, \ \alpha = 1,\dots,n_k$,

\begin{align}
    \WW_{\alpha i} &= \int_E \nabla \cdot \vvarphi_i\ p^k_{\alpha} \\
    &= - \int_E \vvarphi_i \cdot \nabla p^k_{\alpha} + \int_{\partial E} \vvarphi_i \cdot \nn_{\partial E}\ p^k_{\alpha} \label{eq:GaussGreenW}\\ 
    &\coloneqq \left(\WW_1\right)_{\alpha i} + \left(\WW_2\right)_{\alpha i} \nonumber. 
\end{align}

Since

\begin{equation*}
     \nabla p^k_{\alpha}  = \g^{\nabla,k-1}_{\alpha-1} \quad \forall \alpha = 2,\dots,n_{k},
\end{equation*}
in the monomial and partial-orthonormal approaches, by exploiting the internal $\nabla$ degrees of freedom, we have

\begin{equation*}
    \WW_1 = \begin{bmatrix}
    \mathbf{O} \in \R^{1 \times N^{\dof}_E}\\
    \begin{matrix}
    \mathbf{O} \in \R^{n^{\nabla}_{k-1} \times N^{E,e}(k+1)} & - \vert E \vert \mathbf{I} \in \R^{n^{\nabla}_{k-1} \times n^{\nabla}_{k-1} }& \mathbf{O} \in \R^{n^{\nabla}_{k-1} \times n^{\bot}_k}
    \end{matrix}
    \end{bmatrix}.
\end{equation*}
In the full-orthonormal approach, $\forall \alpha = 1,\dots,n_{k-1}^{\nabla}$, we have

\begin{align*}
    \nabla p^k_{\alpha + 1} &= \sum_{I = 1}^{2n_{k-1}} \twoTnablak_{\alpha I} \pp_I^{k-1} 
    = \sum_{I = 1}^{2n_{k-1}} \sum_{\beta =1}^{n^{\nabla}_{k-1}} \left(\mathbf{L}^{\nabla,k-1}\right)^{-1}_{\alpha \beta} \fourTnablak_{\beta I} \pp_I^{k-1} \\
    &= \sum_{\beta =1}^{n^{\nabla}_{k-1}} \left(\mathbf{L}^{\nabla,k-1}\right)^{-1}_{\alpha \beta} \fourgnablakm_{\beta}.
\end{align*}
Thus, thanks to the hierarchical property shown in Remark \ref{rem:gnabla_hierarchical}, the $\WW_1$ matrix becomes

\begin{equation*}
    \WW_1 = \begin{bmatrix}
    \mathbf{O} \in \R^{1 \times N^{\dof}_E}\\
    \begin{matrix}
    \mathbf{O} \in \R^{n^{\nabla}_{k-1} \times N^{E,e}(k+1)} & - \vert E \vert \left(\LL^{\nabla,k}\right)^{-1}(1:n^{\nabla}_{k-1}, 1:n^{\nabla}_{k-1}) & \mathbf{O} \in \R^{n^{\nabla}_{k-1} \times n^{\bot}_k}
    \end{matrix}
    \end{bmatrix}.
\end{equation*}

Let now $^{p}\mathbf{V}^{\partial,k} \in \R^{N^{E,e}(k+1) \times n_k}$ be the Vandermonde matrix of polynomials of order less or equal to $k$ related to the boundary quadrature formula $\left(\{\{\x_i^{e,Q}\},\{w_i^{e,Q}\}_{i=1}^{k+1}\}\right)_{e \in \mathcal{E}_{h,E}}$ defined on $\partial E$, whose quadrature nodes coincide with the edge dofs defined in \eqref{eq:DOFedges}. Note that the integrand function of the second integral of \eqref{eq:GaussGreenW} is a known polynomial of degree $2k$ on each edge of $E$. Thus, it can be computed exactly (up to machine precision) by evaluating the integrand at the ($k+1$)-edge dofs on each edge.
The second term $\mathbf{W}_2$ then reads

\begin{equation*}
    \mathbf{W}_2 = \begin{bmatrix}
    \left(^{p}\mathbf{V}^{\partial,k}\right)^T \mathbf{W}^{\partial, Q} & \mathbf{O} \in \R^{n_k \times (n^{\nabla}_{k-1} + n^{\bot}_k)}
    \end{bmatrix},
\end{equation*}
where $\mathbf{W}^{\partial,Q} \in \R^{N^{E,e}(k+1) \times N^{E,e}(k+1)}$ is the diagonal matrix whose non-zero entries coincide with the boundary quadrature weights properly arranged.

\subsection{Computation of \texorpdfstring{$L^2(E)$}{L2E}-projection of basis functions of \texorpdfstring{$\bm{V}_k(E)$}{VEM space for the velocity variable}}

Now, we compute the projection $\PPi^0_k \vvarphi_i$ of basis functions of $\V_k(E)$ in terms of the $\mathcal{G}^k$ basis of $\Pkk$, i.e.
\begin{equation}
   \PPi^0_k \vvarphi_i = \sum_{I=1}^{2n_k}  \left(\mathbf{\Pi^0_k}\right)_{Ii} \g^k_I.
   \label{eq:PioKFunctionG}
\end{equation}
By replacing (\ref{eq:PioKFunctionG}) in the definition (\ref{eq:Pi0k}) of $\PPi^0_k$, we obtain, $\forall i =1,\dots,N^{\dof}_E$, $\forall J =1,\dots,2n_k$,

\begin{equation*}
    \sum_{I=1}^{2n_k} \left(\mathbf{\Pi^0_k}\right)_{Ii} \left(\g^k_I,\g^k_J\right)_{0,E} = \left(\vvarphi_i,\g^k_J\right)_{0,E}.
\end{equation*}
or, analogously, in matrix form

\begin{equation*}
    \mathbf{G}^k \mathbf{\Pi^0_k} = \mathbf{B},
\end{equation*}
where $\mathbf{B} \in \R^{2n_k \times N^{\dof}_E}$ is the matrix whose entries are defined as

\begin{equation*}
    \mathbf{B}_{Ji} =  \left(\vvarphi_i,\g^k_J\right)_{0,E}, \quad \forall i =1,\dots,N^{\dof}_E,\ \forall J =1,\dots,2n_k.
\end{equation*}
Now, as in \cite{MixedImplem}, we split $\mathbf{B}$ as

\begin{equation*}
    \mathbf{B} = \begin{bmatrix}
    \mathbf{B}^{\nabla}\\
    \mathbf{B}^{\bot}
    \end{bmatrix}.
\end{equation*}
The term $\mathbf{B}^{\bot} \in \R^{n^{\bot}_k \times N^{\dof}_E}$, whose entries are

\begin{equation*}
       \mathbf{B}_{\alpha i}^{\bot} =  \left(\vvarphi_i,\g^{\bot,k}_{\alpha}\right)_{0,E}, \quad \forall \alpha =1,\dots,n^{\bot}_k, \ \forall i =1,\dots,N^{\dof}_E,
\end{equation*}
can be readily computed as

\begin{equation*}
    \mathbf{B}^{\bot} = \begin{bmatrix}
    \mathbf{O} \in \R^{n^{\bot}_k \times (N^{E,e}(k+1)+n^{\nabla}_{k-1})} & \vert E \vert \mathbf{I} \in R^{n^{\bot}_k \times n^{\bot}_k}
    \end{bmatrix}.
\end{equation*}
Concerning the first term $\mathbf{B}^{\nabla} \in  \R^{n^{\nabla}_k \times N^{\dof}_E}$, we note that, $\forall \alpha =1,\dots,n^{\nabla}_k$, $\forall i =1,\dots,N^{\dof}_E$
\begin{itemize}
    \item in the monomial and in the partial-orthonormal approaches,
    
    \begin{align}
       \nonumber \mathbf{B}_{\alpha i}^{\nabla} &=  \int_E \vvarphi_i \cdot \g^{\nabla,k}_{\alpha} = \nonumber \int_E \vvarphi_i \cdot \nabla p^{k+1}_{\alpha + 1}\\
       &= \label{eq:BNablaGaussGreen}-\int_E \nabla \cdot \vvarphi_i \ p^{k+1}_{\alpha + 1} + \int_{\partial E} \vvarphi_i \cdot \nn_{\partial E}\ p^{k+1}_{\alpha +1 } \\
       &\nonumber \coloneqq \left(\mathbf{B}^{\nabla}_1\right)_{\alpha i} + \left(\mathbf{B}^{\nabla}_2\right)_{\alpha i},
    \end{align}
    \item in the full-orthonormal approach
    
    \begin{align*}
       \mathbf{B}_{\alpha i}^{\nabla} &=  \int_E \vvarphi_i \cdot \g^{\nabla,k}_{\alpha}
       = \sum_{\beta} \mathbf{L}^{\nabla,k}_{\alpha \beta} \int_E \vvarphi_i \cdot \nabla p^{k+1}_{\beta + 1}
       =  \left(\mathbf{L}^{\nabla,k}\mathbf{B}^{\nabla}_1\right)_{\alpha i} + \left(\mathbf{L}^{\nabla,k}\mathbf{B}^{\nabla}_2\right)_{\alpha i}.
    \end{align*}
\end{itemize}
Since the integrand of the second term of \eqref{eq:BNablaGaussGreen} is a known polynomial of degree $2k+1$ on each edge of $E$, it can be integrated exactly by exploiting the edge dofs. Thus,

\begin{equation*}
    \mathbf{B}^{\nabla}_2 = \begin{bmatrix}
    \left(^{p}\mathbf{V}^{\partial,k+1}(:,2:n_{k+1})\right)^T \mathbf{W}^{\partial, Q} & \mathbf{O} \in \R^{n^{\nabla}_k \times (n^{\nabla}_{k-1} + n^{\bot}_k)}
    \end{bmatrix}.
\end{equation*}
In order to compute the term $\mathbf{B}^{\nabla}_1$, we first note that $\nabla \cdot \vvarphi_i$ is a known polynomial of degree $k$ on $E$. Then, it can be written as

\begin{equation*}
    \nabla \cdot \vvarphi_i = \sum_{\alpha = 1}^{n_k} \mathbf{\Lambda}_{\alpha i } p^k_{\alpha}, \quad \forall i = 1,\dots,N^{\dof}_E. 
\end{equation*}
whose coefficient matrix $\mathbf{\Lambda} \in \R^{n_k \times N^{\dof}_E}$ can be retrieved from 

\begin{equation}
    \sum_{\alpha = 1}^{n_k} \mathbf{\Lambda}_{\alpha i } \int_E p^k_{\alpha} p^k_{\beta} = \int_E \nabla \cdot \vvarphi_i\ p^k_{\beta}, \quad \forall \beta = 1,\dots,n_k,\ \forall i = 1,\dots, N^{\dof}_E.
    \label{eq:coeffDiver}
\end{equation}
Equation \eqref{eq:coeffDiver} can be written in matrix form as

\begin{equation*}
    {^{p}\mathbf{H}^k} \mathbf{\Lambda} = \WW,
\end{equation*}
since the term at right hand side of \eqref{eq:coeffDiver} coincides with the definition of the $\mathbf{W}$ matrix given in Section \ref{sec:Wmatrix}.
We highlight that, in the partial and in the full orthonormal approaches, in infinite precision, it holds

\begin{equation*}
    \mathbf{\Lambda} = \WW.
\end{equation*}
Finally, matrix $\mathbf{B}^{\nabla}_1$ is given by
\begin{align*}
    \left(\mathbf{B}^{\nabla}_1\right)_{\alpha i} &= - \int_E \nabla \cdot \vvarphi_i p^{k+1}_{\alpha +1}
    = - \sum_{\beta = 1}^{n_k} \mathbf{\Lambda}_{\beta i} \int_E p^k_{\beta} p^{k+1}_{\alpha + 1},
\end{align*}
or, equivalently,
\begin{align*}
    \mathbf{B}^{\nabla}_1 &= - \ {^{p}\mathbf{H}^{k+1}}(2:n_{k+1},1:n_k) \mathbf{\Lambda}
    = -\ {^{p}\mathbf{H}^{k+1}}(2:n_{k+1},1:n_k) ({^{p}\mathbf{H}^k})^{-1} \WW.
\end{align*}

\subsection{The diffusion term \texorpdfstring{$a^E_h(\uu_h,\vv_h)$}{}}

This term is defined by equation \eqref{eq:ahE} as the sum of a consistency and a stability term, thus, as in \cite{MixedImplem}, we define a consistency matrix 

\begin{equation*}
    \left(\mathbf{K}^a_C\right)_{ij} = (\K \PPi^0_k \vvarphi_i, \PPi^0_k \vvarphi_j)_{0,E}
\end{equation*}
and a stability matrix

\begin{equation*}
    \left(\mathbf{K}^a_S\right)_{ij} =\overline{\K} \vert E \vert \sum_{r=1}^{N^{\dof}_E} \dof_r\left((I-\PPi^0_k) \vvarphi_i \right)\dof_r\left((I-\PPi^0_k)\vvarphi_j\right).
\end{equation*}
The consistency matrix can be computed as

\begin{equation*}
    \mathbf{K}^a_C = \left(\mathbf{\Pi^0_k}\right)^T \mathbf{G}^{\bm{K}} \mathbf{\Pi^0_k},
\end{equation*}
where

\begin{equation*}
    \mathbf{G}^{\bm{K}} = \left(\ggVk\right)^T \begin{bmatrix}
    \mathbf{W}^{Q,K_{xx}} & \mathbf{W}^{Q,K_{xy}} \\
    \mathbf{W}^{Q,K_{xy}} & \mathbf{W}^{Q,K_{yy}} 
    \end{bmatrix} \ggVk,
\end{equation*}
\begin{equation*}
    \bm{K} = \begin{bmatrix}
K_{xx} & K_{xy}\\
K_{xy} & K_{yy}
\end{bmatrix},
\end{equation*}
and $\mathbf{W}^{Q,K_{\ast}}\in \R^{N^Q \times N^Q}$ is the diagonal matrix whose non zero entries are

\begin{equation*}
    \mathbf{W}^{Q,K_{\ast}}_{ii} = K_{\ast}(\x_i^Q)w_i^Q,\quad \forall i = 1,\dots,N^Q\ \forall \ast \in \{xx, xy, yy\}.
\end{equation*}

Now, we define matrix $\mathbf{D} \in \R^{N^{\dof}_E \times 2n_k}$, whose entries are

\begin{equation*}
    \mathbf{D}_{iI} = \dof_i(\g^k_I),\quad \forall i=1,\dots,N^{\dof}_E,\ \forall I=1,\dots,2n_k.
\end{equation*}
Matrix $\mathbf{D}$ can be split as

\begin{equation*}
    \mathbf{D} = \begin{bmatrix}
    \mathbf{D}^{\partial}\\
    \mathbf{D}^{\nabla}\\
    \mathbf{D}^{\bot}
    \end{bmatrix},
\end{equation*}
being matrices $\mathbf{D}^{\nabla} \in \R^{n^{\nabla}_{k-1} \times 2n_k}$ and $\mathbf{D}^{\bot} \in \R^{n^{\bot}_{k} \times 2n_k}$ given by

\begin{equation*}
    \mathbf{D}^{\nabla} = \frac{1}{\vert E \vert }\mathbf{G}^k(1:n^{\nabla}_{k-1},:), \quad 
    \mathbf{D}^{\bot} = \frac{1}{\vert E \vert }\mathbf{G}^k(n^{\nabla}_k +1:2n_k,:).
\end{equation*}
Matrix $\mathbf{D}^{\partial} \in \R^{N^{E,e}(k+1) \times 2n_k}$, instead, $\forall i =1,\dots,k+1$, $\forall e \in \mathcal{E}_{h,E}$, $\forall I =1,\dots,2n_k,$, reads

\begin{equation*}
    \mathbf{D}^{\partial}_{iI} = \dof_i^e(\g_I^k) = \left(\g_I^k \cdot \nn_e\right)(\x^{e,Q}_i),
\end{equation*}
or, equivalently,

\begin{align*}
    \mathbf{D}^{\partial} &= \begin{bmatrix}
    \mathbf{N}_x & \mathbf{N}_y
    \end{bmatrix}\ {^{\g}\mathbf{V}^{\partial, k}}\\
    &= \begin{bmatrix}
    \mathbf{N}_x & \mathbf{N}_y
    \end{bmatrix}\ \begin{bmatrix} {^{p}\mathbf{V}^{\partial,k}} & \mathbf{O} \in \mathbb{R}^{N^{E,e}(k+1) \times n_k}\\
    \mathbf{O} \in \mathbb{R}^{N^{E,e}(k+1) \times n_k} & {^{p}\mathbf{V}^{\partial,k}} \end{bmatrix}\ \begin{bmatrix}
    \left(\mathbf{T}^{\nabla,k}\right)^T & \left(\mathbf{T}^{\bot,k}\right)^T
    \end{bmatrix},
\end{align*}
where $\mathbf{N}_{\ast} \in \R^{N^{E,e}(k+1) \times N^{E,e}(k+1) }$, for all $\ast \in \{x,y\}$, is the diagonal matrix whose diagonal entries are the $\ast$-component of the normal vectors at the edges of $E$, properly arranged.

Finally, the stability matrix can be computed as

\begin{equation*}
    \mathbf{K}^a_S = \overline{\K}\vert E \vert\left(\mathbf{I} - \mathbf{D}\mathbf{\Pi^0_k}\right)^T\left(\mathbf{I} - \mathbf{D}\mathbf{\Pi^0_k}\right).
\end{equation*}

\subsection{The advection term \texorpdfstring{$- \left(\bbeta p_h, \PPi^0_k \vv_h\right)_{0,E}$}{}}

The matrix $\mathbf{T}^{\beta} \in \R^{N^{\dof}_E \times n_k}$ corresponding to the local term

\begin{equation*}
    - \left(\bbeta p_h, \vv_h\right)_{0,E},
\end{equation*}
$\forall i=1,\dots,N^{\dof}_E$, $\forall \alpha=1,\dots,n_k$, is defined as

 \begin{align*}
    \mathbf{T}^{\beta}_{i \alpha} &= - \int_E  \bbeta \cdot \PPi^0_k \vvarphi_i\ p_{\alpha}^k
    = - \sum_{I=1}^{2n_k} \left(\mathbf{\Pi^0_k}\right)_{Ii} \int_E \bbeta \cdot \g_I^k \ p_{\alpha}^k.
 \end{align*}
Given $\bbeta = \begin{bmatrix}
 
 \beta_x, \beta_y
 \end{bmatrix}^T$, in matrix form $\mathbf{T}^{\beta}$ reads
 
 \begin{equation*}
    \mathbf{T}^{\bbeta} = \left(\mathbf{\Pi^0_k}\right)^T \left(\ggVk\right)^T \begin{bmatrix}
    \mathbf{W}^{Q,\beta_x} & \mathbf{O}\\
    \mathbf{O} & \mathbf{W}^{Q,\beta_y}
    \end{bmatrix} \begin{bmatrix}
    {^{p}\mathbf{V}^k}\\ {^{p}\mathbf{V}^k}
    \end{bmatrix},
 \end{equation*}
 where $\mathbf{W}^{Q,\beta_{\ast}} \in \R^{N^Q \times N^Q}$ $\forall \ast \in \{x,y\}$ is the diagonal matrix whose non zero entries are defined as
 
 \begin{equation*}
     (\mathbf{W}^{Q,\bbeta_{\ast}})_{ii} = \beta_{\ast}(\x_i^{Q}) w_i^{Q}, \quad \forall i = 1,\dots,N^Q.
 \end{equation*}

\subsection{The reaction term \texorpdfstring{$\left(\gamma p_h,q_h\right)_{0,E}$}{}}

We define the local matrix $\mathbf{H}^{\gamma} \in \R^{n_k\times n_k}$ as the matrix that collects the terms

\begin{equation*}
   \mathbf{H}^{\gamma}_{\alpha \beta} = \int_E \gamma p^k_{\alpha} p^k_{\beta}\quad \forall \alpha,\beta =1,\dots,n_k. 
\end{equation*}
Thus,

\begin{equation*}
   \mathbf{H}^{\gamma} = \left({^{p}\mathbf{V}^k}\right)^T \mathbf{W}^{Q,\gamma}\  {^{p}\mathbf{V}^k},
\end{equation*}
where 

\begin{equation*}
    \mathbf{W}^{Q,\gamma}_{i j} = \gamma(\x_i^Q) w_j^Q \delta_{ij}, \quad \forall i,j = 1,\dots,N^Q.
\end{equation*}

\bigskip
\bigskip

In conclusion, the local system matrix $\mathbf{K}^E \in \R^{N^{\dof}_E \times N^{\dof}_E}$ for the discrete problem (\ref{eq:discreteMixedVariationalFomrulation}) is given by

\begin{equation*}
    \mathbf{K}^E = \begin{bmatrix}
    \mathbf{K}^{a}_C + \mathbf{K}^{a}_S & -\mathbf{W}^T + \mathbf{T}^{\beta}\\
    \mathbf{W} & \mathbf{H}^{\gamma}
    \end{bmatrix}.
\end{equation*}

\section{Numerical experiments}\label{sec:numericalExp}

In this section, we propose some numerical experiments to show the performance of the aforementioned approaches in terms of the following error norms:

\begin{equation}
    p_{\text{err}} = \left(\sum_{E \in \mathcal{T}_h} \| p - p_h \|_{0,E}^2\right)^{\frac{1}{2}},
    \label{eq:L2pressure}
\end{equation}
\begin{equation}
    \uu_{\text{err}} = \left(\sum_{E \in \mathcal{T}_h} \| \uu - \PPi^0_k \uu_h \|_{0,E}^2\right)^{\frac{1}{2}},
    \label{eq:L2velocity}
\end{equation}
\begin{equation}
    p_{I,\text{err}} = \left(\sum_{E \in \mathcal{T}_h} \| p_I - p_h \|_{0,E}^2\right)^{\frac{1}{2}},
    \label{eq:superconvergence}
\end{equation}
where the interpolant $p_I \in \mathcal{Q}_h$  of $p$ is computed locally as

\begin{equation*}
        p_I = \sum_{\alpha = 1}^{n_k} \bm{c}_{\alpha} p^{k}_{\alpha}.
\end{equation*}
The vector $\bm{c}\in \R^{n_k}$ of coefficients is the solution of the following linear least squares problem

\begin{equation*}
       \min \| {^{p}\mathbf{V}^{k}} \bm{c} - \bm{y} \|,
\end{equation*}
where $\bm{y} \in \R^{N^Q}$ is the vector whose entries are the evaluation of $p$ at the internal quadrature points $\{\x^Q_i\}_{i=1}^{N^Q}$.

We will also analyze the condition number of the main local matrices that allow assembling the local system matrix, namely $\mathbf{G}^k$, $\mathbf{W}$, $\mathbf{B}$, $\mathbf{\PPi^0_k}$ and $\mathbf{D}$. The condition number of a matrix is computed as the ratio between its largest and its smallest singular value.

For the first two numerical examples, we consider problem \eqref{eq:continuousProblem} on the unit square domain $\Omega = (0,1)^2$ with

\begin{equation*}
    \D = \begin{bmatrix}
    y^2+1 & -xy\\
    -xy & x^2+1
    \end{bmatrix}, \quad \bb = \begin{bmatrix}
    x\\
    y
    \end{bmatrix}, \quad \gamma = x^2+y^3.
\end{equation*}

The forcing term, the Dirichlet and the Neumann boundary conditions are set in such a way the exact solution is
\begin{equation}
    p(x,y) = x^2y + \sin(2\pi x) \sin(2\pi y) + 2,
    \label{eq:exactSolution}
\end{equation}
and the set $\Gamma_N$ coincides with the edge of the domain of interest on the $x$-axis.

In particular, in the first experiment (Test1), we validate the proposed methods showing the order of convergence of the three monomial, partial and full orthogonal approaches against rising polynomial degree $k$, for various mesh refinement levels, on meshes made of square elements. Then, in the second test (Test2), we show the performance of the approaches in terms of matrix condition number and error convergence trends in presence of collapsing polygons in the meshes.

Finally, the last example (Test3) proposes the application of the method to flow simulations in Discrete Fracture Networks (DFNs). DFN simulations, are, indeed, a typical example where highly distorted mesh elements might appear in the process to generate a conforming mesh \cite{DFN}. A simple network is here considered, for which an analytic solution is known, in order to compare the convergences curve of the three different approaches.

\subsection{Test1: Convergence rates}
In this first experiment, we validate the new approaches by showing that the computational orders of convergence match the theoretical ones for different values of the polynomial degree $k$ on regular meshes. At this aim, convergence tests are performed on a sequence of squared meshes that decompose the domain of interest in 25, 100, 400 and 1600 identical squares, respectively.

\begin{figure}[htbp]
	\centering
	\subfigure[\label{fig:ErroreL2Pressure_ret10x10}] 
	{\includegraphics[width=.3\textwidth, height = .17\textheight]{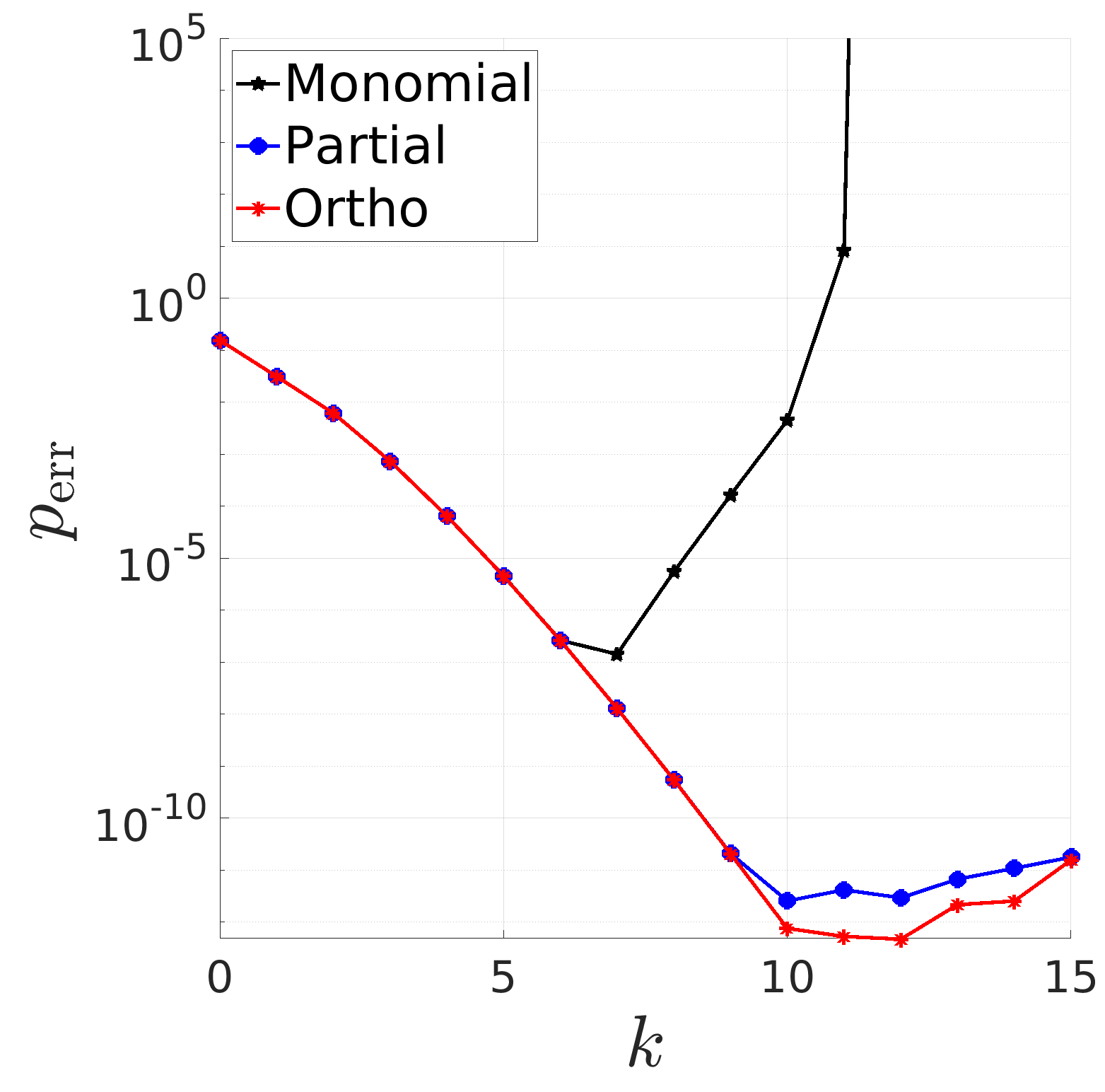}}
	\subfigure[\label{fig:ErroreL2Velocity_ret10x10}]
	{\includegraphics[width=.3\textwidth, height = .17\textheight]{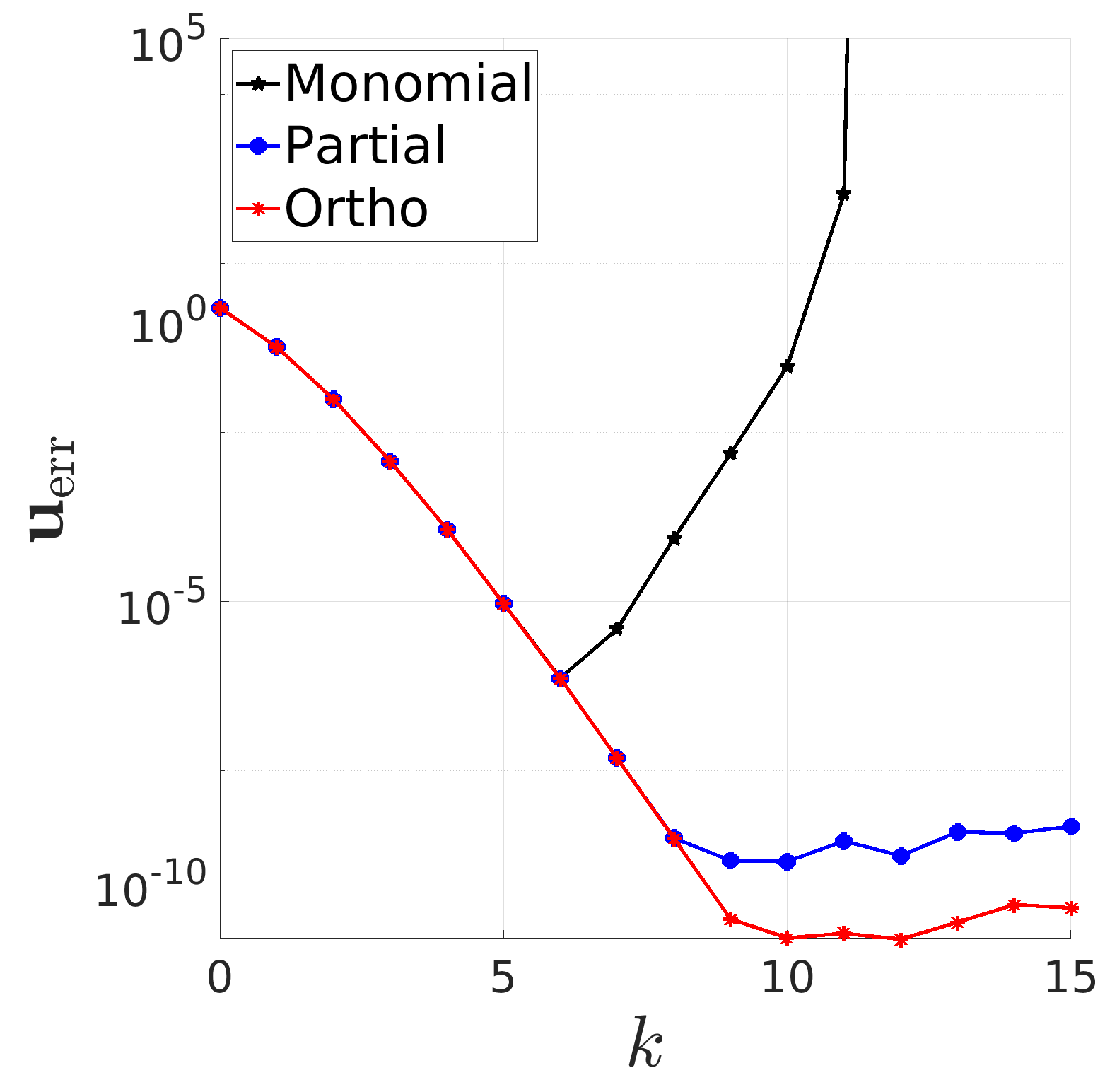}}
	\subfigure[\label{fig:SuperConvergence_ret10x10}] 
	{\includegraphics[width=.3\textwidth, height = .17\textheight]{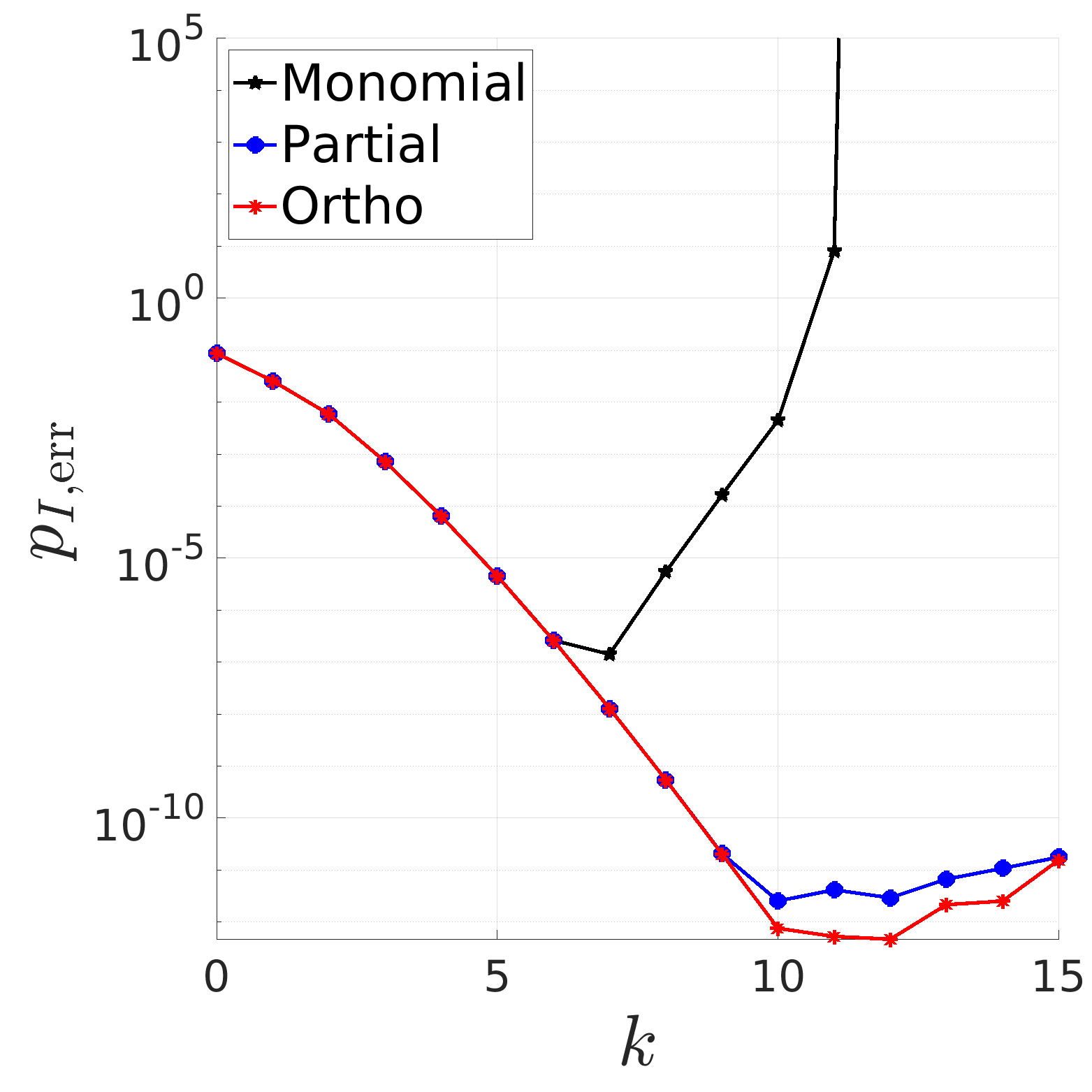}}
	\centering
	\subfigure[\label{fig:ErroreL2Pressure_ret20x20}] 
	{\includegraphics[width=.3\textwidth, height = .17\textheight]{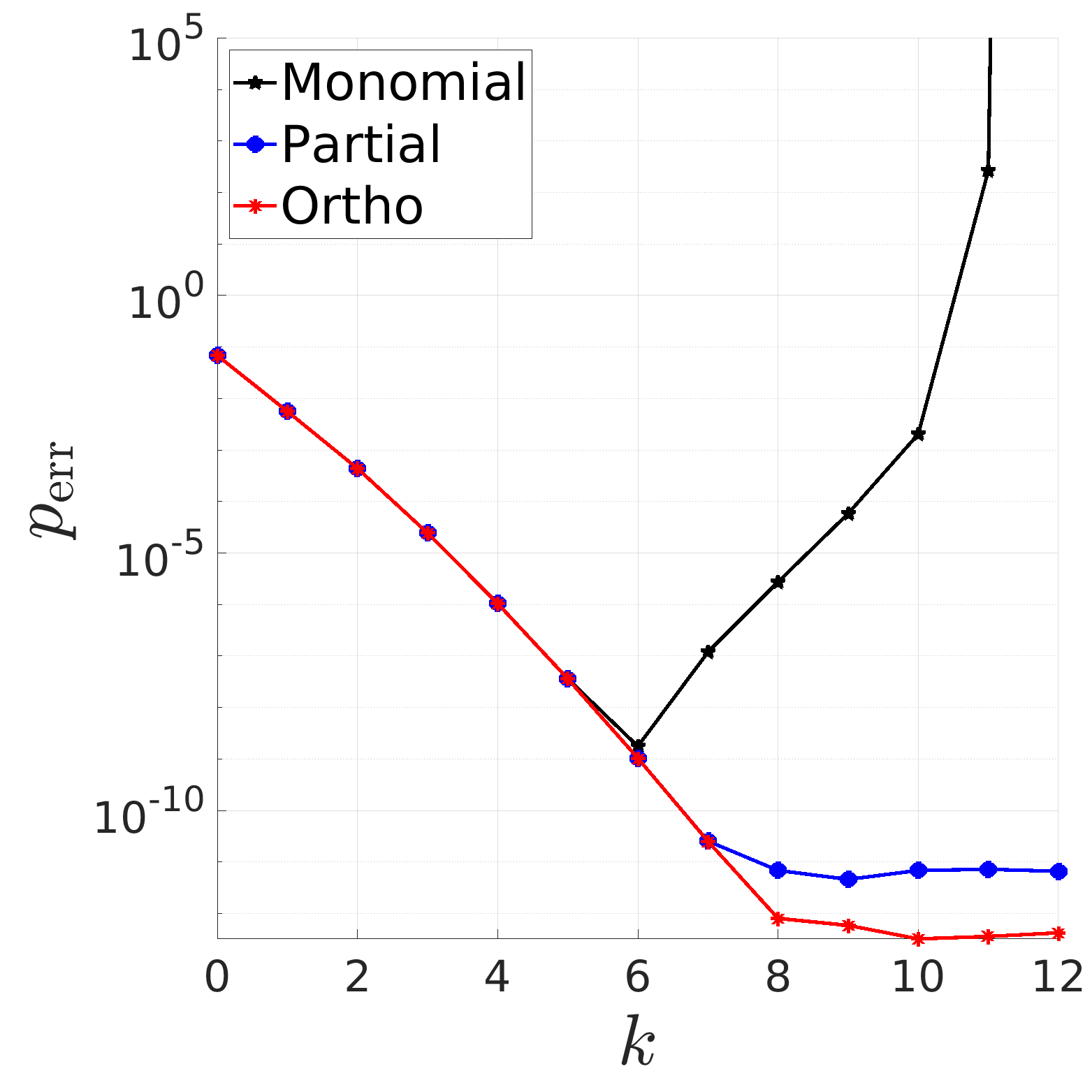}}
	\subfigure[\label{fig:ErroreL2Velocity_ret20x20}]
	{\includegraphics[width=.3\textwidth, height = .17\textheight]{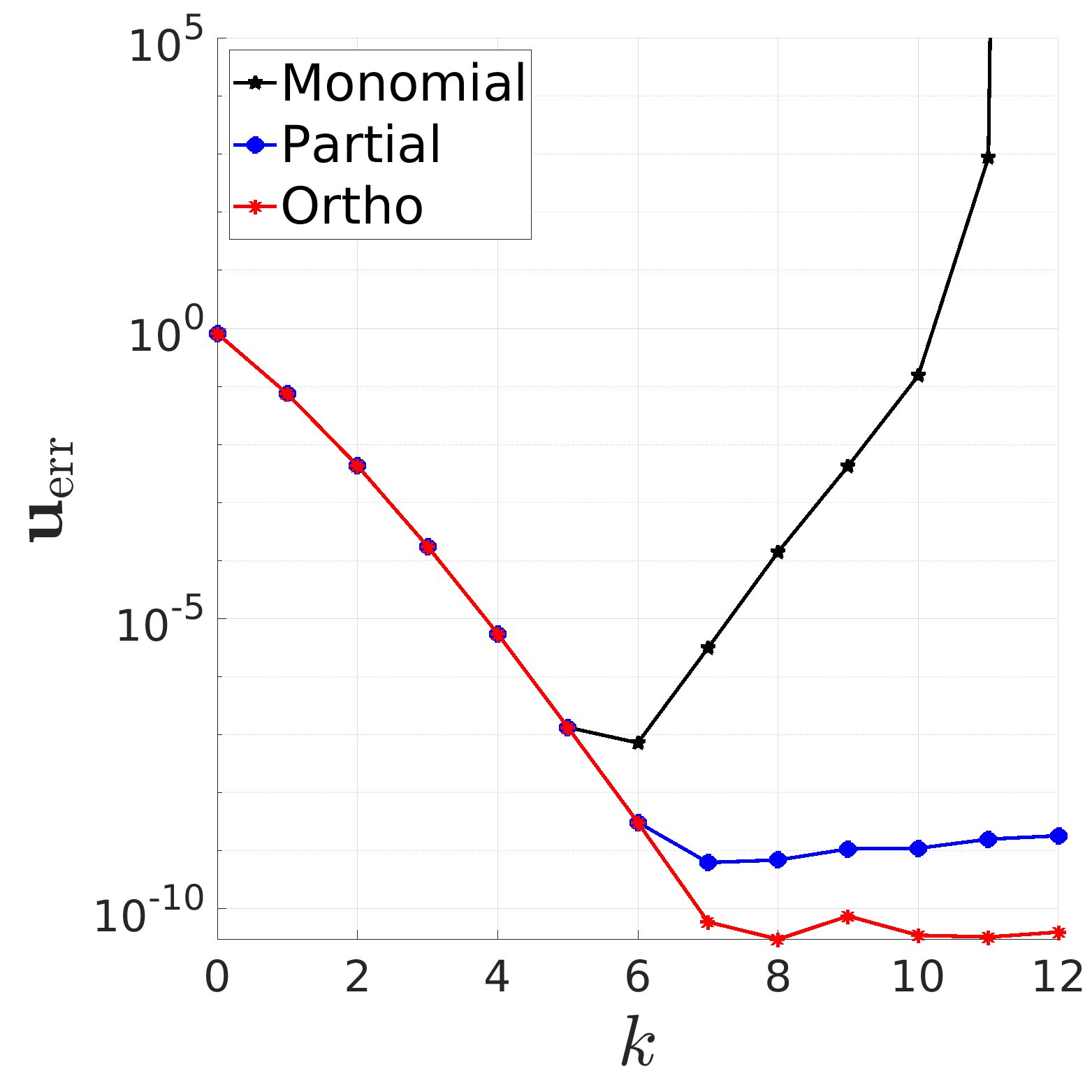}}
	\subfigure[\label{fig:SuperConvergence_ret20x20}] 
	{\includegraphics[width=.3\textwidth, height = .17\textheight]{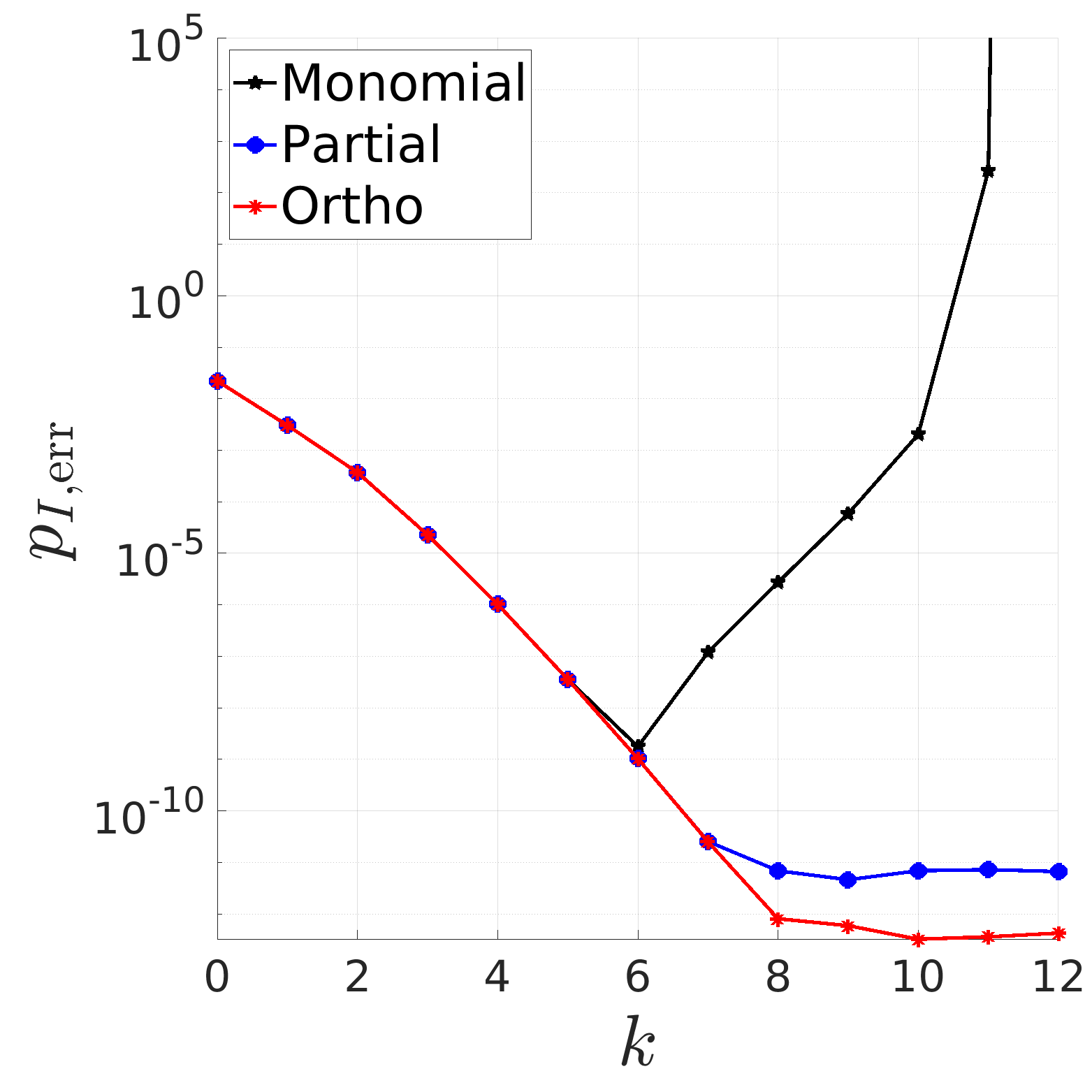}}
	\subfigure[\label{fig:ErroreL2Pressure_ret40x40}] 
	{\includegraphics[width=.3\textwidth, height = .17\textheight]{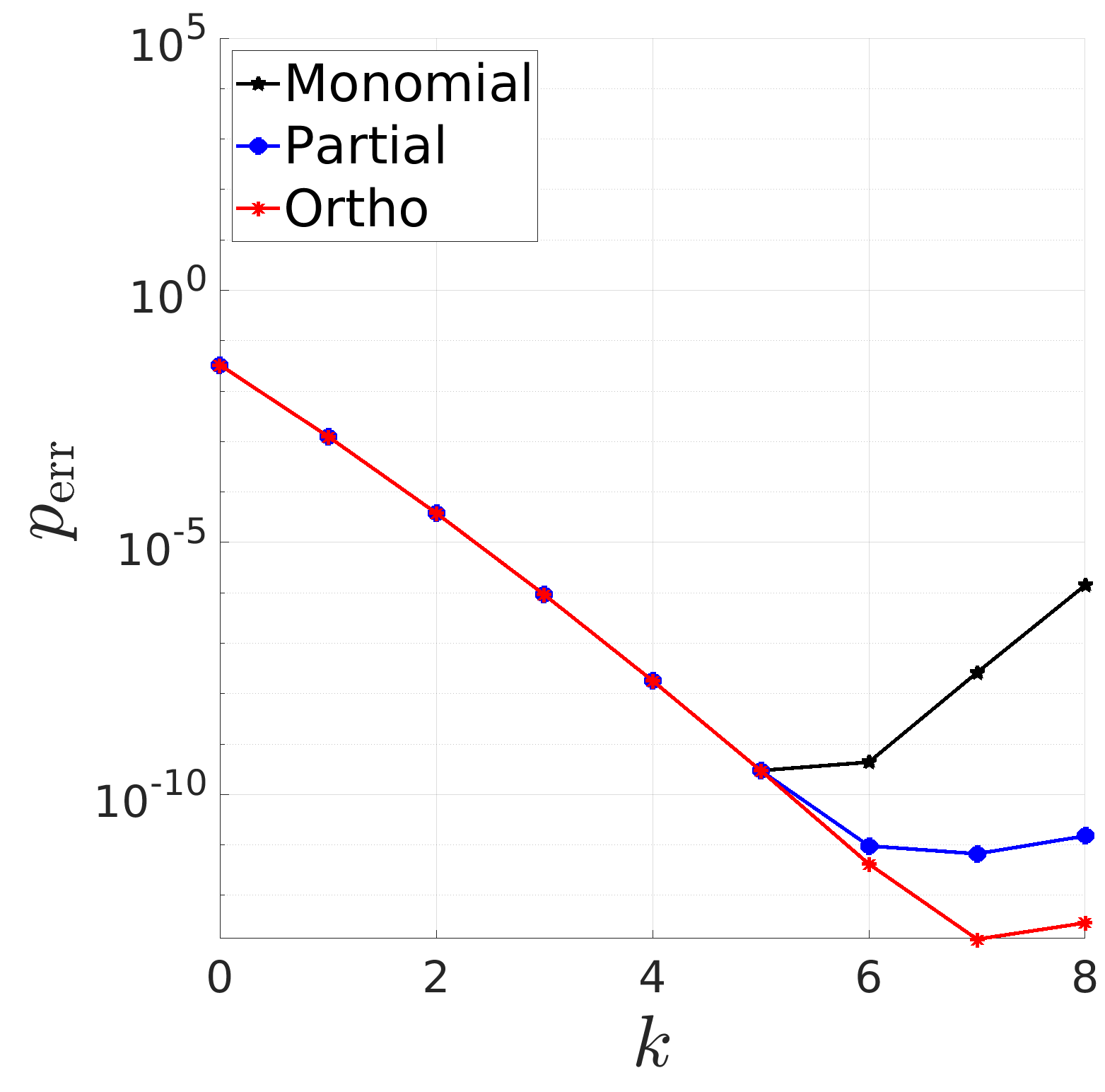}}
	\subfigure[\label{fig:ErroreL2Velocity_ret40x40}]
	{\includegraphics[width=.3\textwidth, height = .17\textheight]{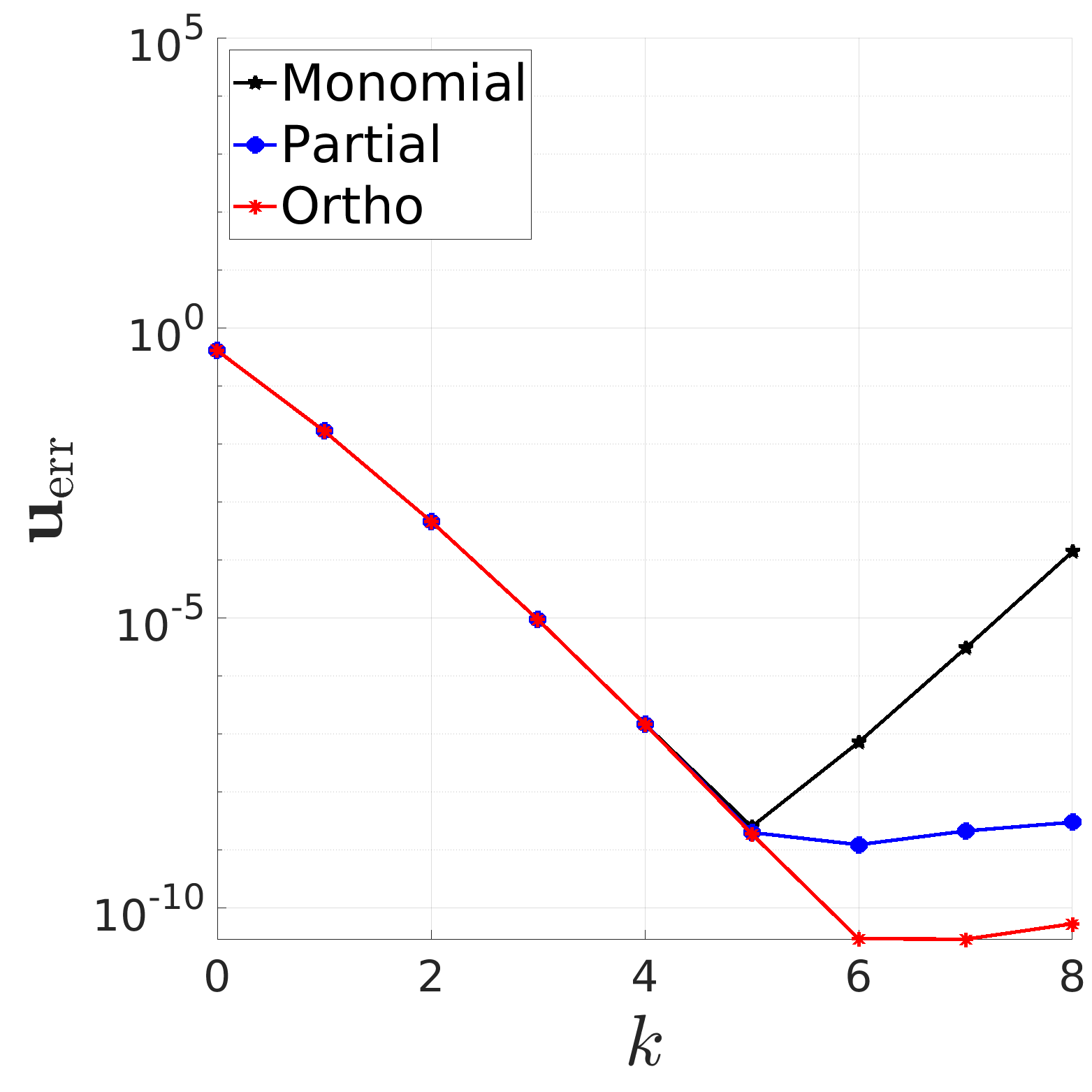}}
	\subfigure[\label{fig:SuperConvergence_ret40x40}] 
	{\includegraphics[width=.3\textwidth, height = .17\textheight]{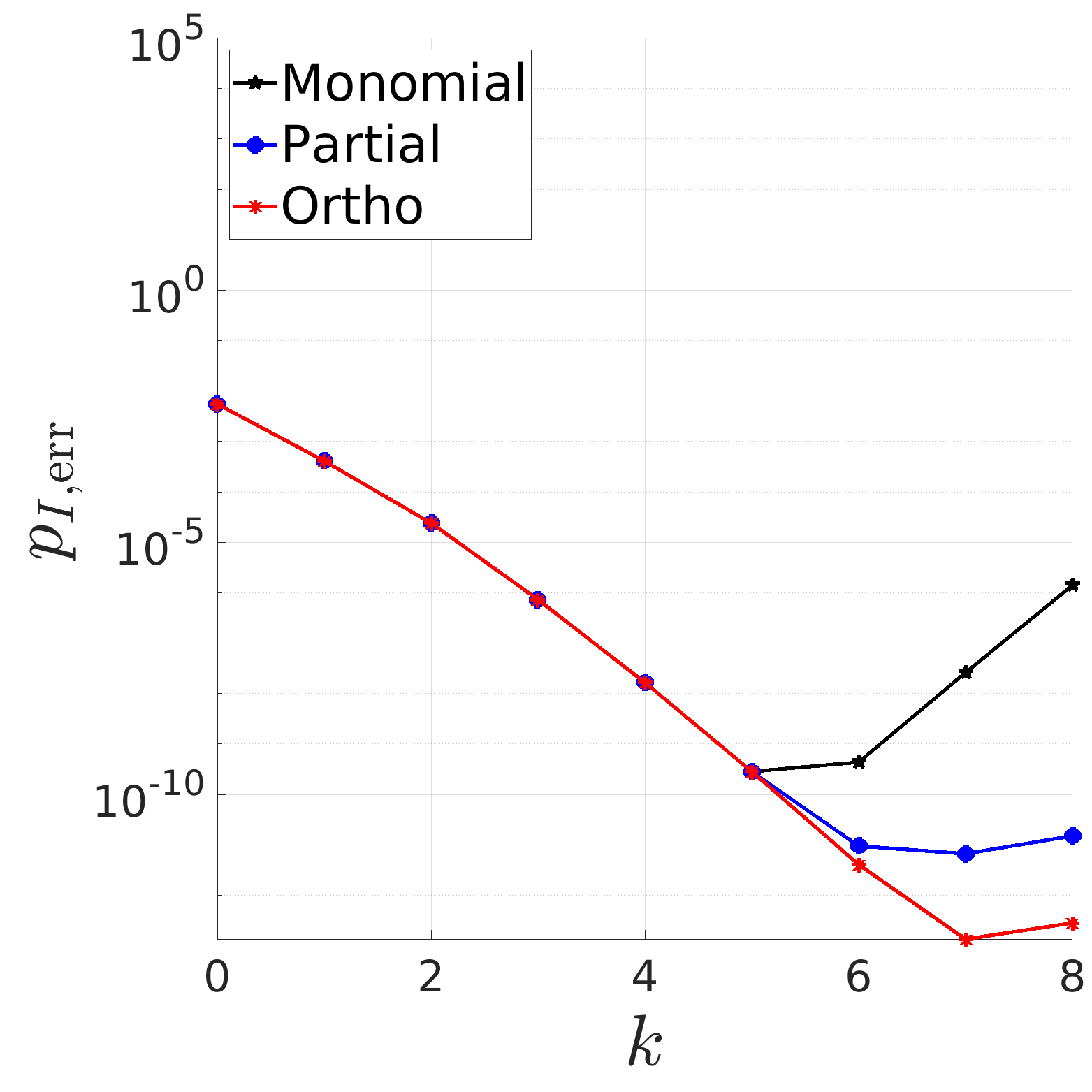}}
    \caption{Test1: Behaviour of errors \eqref{eq:L2pressure}, \eqref{eq:L2velocity} and \eqref{eq:superconvergence} at varying $k$ on square meshes. Pictures on each row refer to a different mesh refinement level: 100, 400 and 1600 element meshes from top to bottom.}
    \label{fig:errors_firstExp}
\end{figure}

Obtained convergence rates are shown in Table \ref{tab:Mixed_SquareConv}, while convergence curves of the errors are shown in Figure \ref{fig:errors_firstExp}. In this figure, each row reports the graphs of the three error norms at varying $k$ on a fixed mesh. The three rows correspond to the last three considered refinement levels, i.e. the 100, 400 and 1600 element meshes, respectively. In the figure, the upper bound of $y$-axis is fixed to $1.0e+5$. 

The obtained results show that the rates of convergence related to all three approaches match the theoretical ones up to polynomial order $k=5$. 
For higher orders, errors in the monomial approach start to increase due to ill-conditioning, while both the partial and the full orthonormal approaches still provide the expected results in terms of errors, up to stagnation due to finite precision arithmetic. We highlight that, in the partial and full orthonormal approach, the local projection matrices and the global system matrix are well-conditioned, as we will see in the next tests.

\begin{table}[]\footnotesize
\centering
\resizebox{\textwidth}{!}{
\begin{tabular}{lccccccccll}
\hline
                  & \textbf{$k$}         & \textbf{0}                 & \textbf{1}                 & \textbf{2}                 & \textbf{3}                 & \textbf{4}                 & \textbf{5}                 & \textbf{6}                 & \multicolumn{1}{c}{\textbf{7}} & \multicolumn{1}{c}{\textbf{8}} \\ \hline
\textbf{Monomial} & $\pp_{\text{err}}$   & 1.1956                     & 2.4760                     & 3.7464                     & 4.8565                     & 5.9153                     & 6.9456                     & 5.8756                     & 2.4067                         & 1.5749                         \\
                  & $\uu_{\text{err}}$   & \multicolumn{1}{l}{0.9947} & \multicolumn{1}{l}{2.0937} & \multicolumn{1}{l}{3.1263} & \multicolumn{1}{l}{4.0877} & \multicolumn{1}{l}{5.1110} & \multicolumn{1}{l}{5.9812} & \multicolumn{1}{l}{3.2093} & 0.3267                         & -0.0167                        \\
\textbf{}         & $\pp_{I,\text{err}}$ & 1.9417                     & 3.0056                     & 3.9580                     & 4.9588                     & 5.9529                     & 6.9680                     & 5.8761                     & 2.4071                         & 1.5749                         \\ \hline
\textbf{Partial}  & $\pp_{\text{err}}$   & 1.1956                     & 2.4760                     & 3.7464                     & 4.8565                     & 5.9153                     & 6.9456                     & 7.6064                     & 6.8582                         & 5.1747                         \\
                  & $\uu_{\text{err}}$   & \multicolumn{1}{l}{0.9947} & \multicolumn{1}{l}{2.0937} & \multicolumn{1}{l}{3.1263} & \multicolumn{1}{l}{4.0877} & \multicolumn{1}{l}{5.1110} & \multicolumn{1}{l}{6.0812} & \multicolumn{1}{l}{5.4221} & 3.8693                         & 2.1663                         \\
\textbf{}         & $\pp_{I,\text{err}}$ & 1.9417                     & 3.0056                     & 3.9580                     & 4.9588                     & 5.9529                     & 6.9680                     & 7.6098                     & 6.8590                         & 5.1751                         \\ \hline
\textbf{Ortho}    & $\pp_{\text{err}}$   & 1.1956                     & 2.4760                     & 3.7464                     & 4.8565                     & 5.9153                     & 6.9456                     & 7.9637                     & 8.5476                         & 7.2031                         \\
                  & $\uu_{\text{err}}$   & \multicolumn{1}{l}{0.9947} & \multicolumn{1}{l}{2.0937} & \multicolumn{1}{l}{3.1263} & \multicolumn{1}{l}{4.0877} & \multicolumn{1}{l}{5.1110} & \multicolumn{1}{l}{6.1058} & \multicolumn{1}{l}{7.0385} & 6.0692                         & 4.3635                         \\
\textbf{}         & $\pp_{I,\text{err}}$ & 1.9417                     & 3.0056                     & 3.9580                     & 4.9588                     & 5.9529                     & 6.9680                     & 7.9773                     & 8.5484                         & 7.2035                         \\ \hline
\end{tabular}}
\caption{Test1: convergence rates on squared mesh.}
\label{tab:Mixed_SquareConv}
\end{table}

\subsection{Test2: Collapsing polygons}

Now, we analyze the condition number of local matrices in case of collapsing polygons. For this purpose, we use three meshes of rectangular elements, with elements having aspect ratios of 10, 50, and 100, respectively. We remark that the aspect ratio of an element is here defined as the ratio between the maximum and the minimum length of its edges. These meshes are built starting from a mesh made up of 100 identical squares, then subdivided into rectangles with length fixed to the value 0.1, equal to that of the original square, and a height computed according to the desired aspect ratio value. 

\begin{figure}[htbp]
	\centering
	\subfigure[\label{fig:condG_ret10x100}] 
	{\includegraphics[width=.3\textwidth, height = .17\textheight]{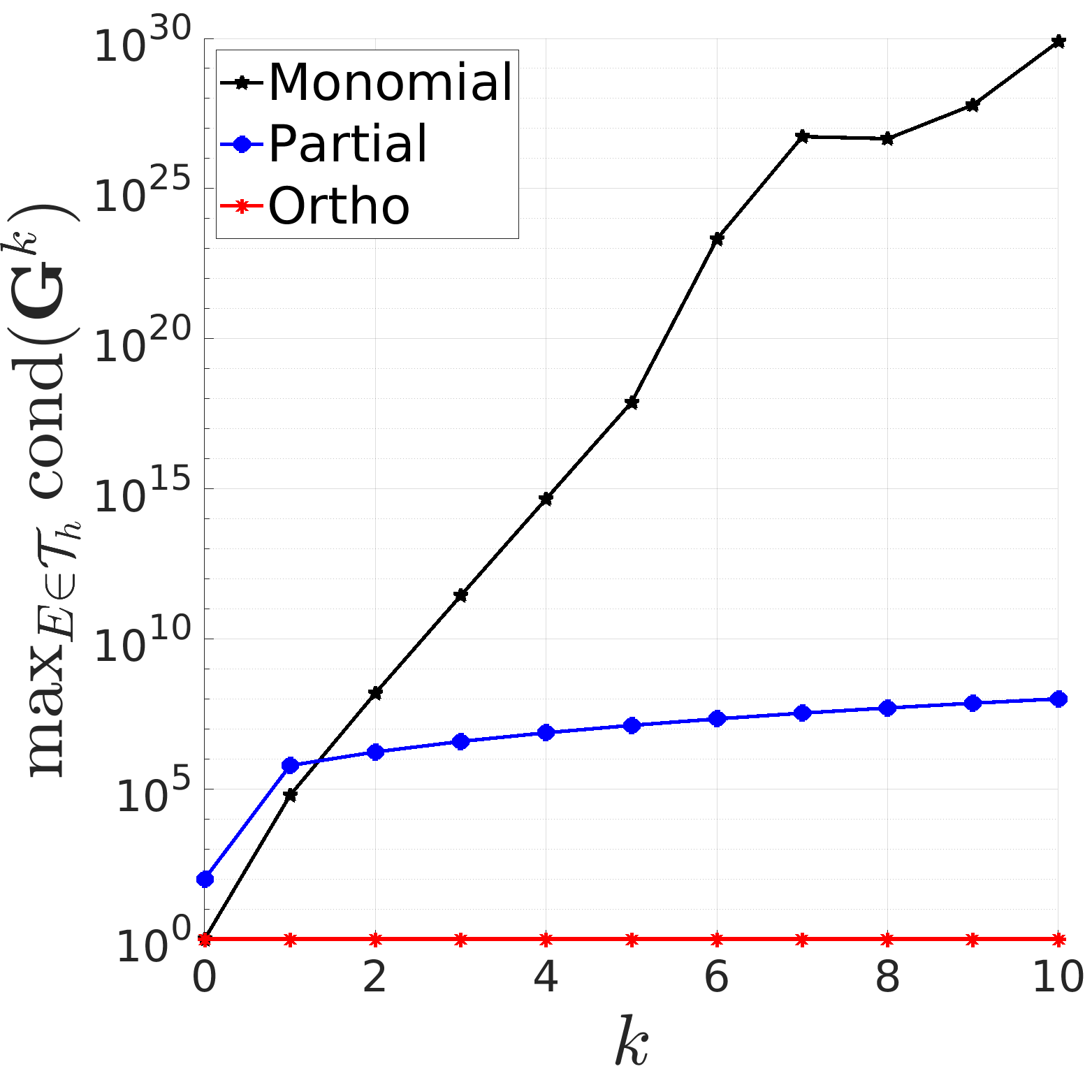}}
	\subfigure[\label{fig:condW_ret10x100}]
	{\includegraphics[width=.3\textwidth, height = .17\textheight]{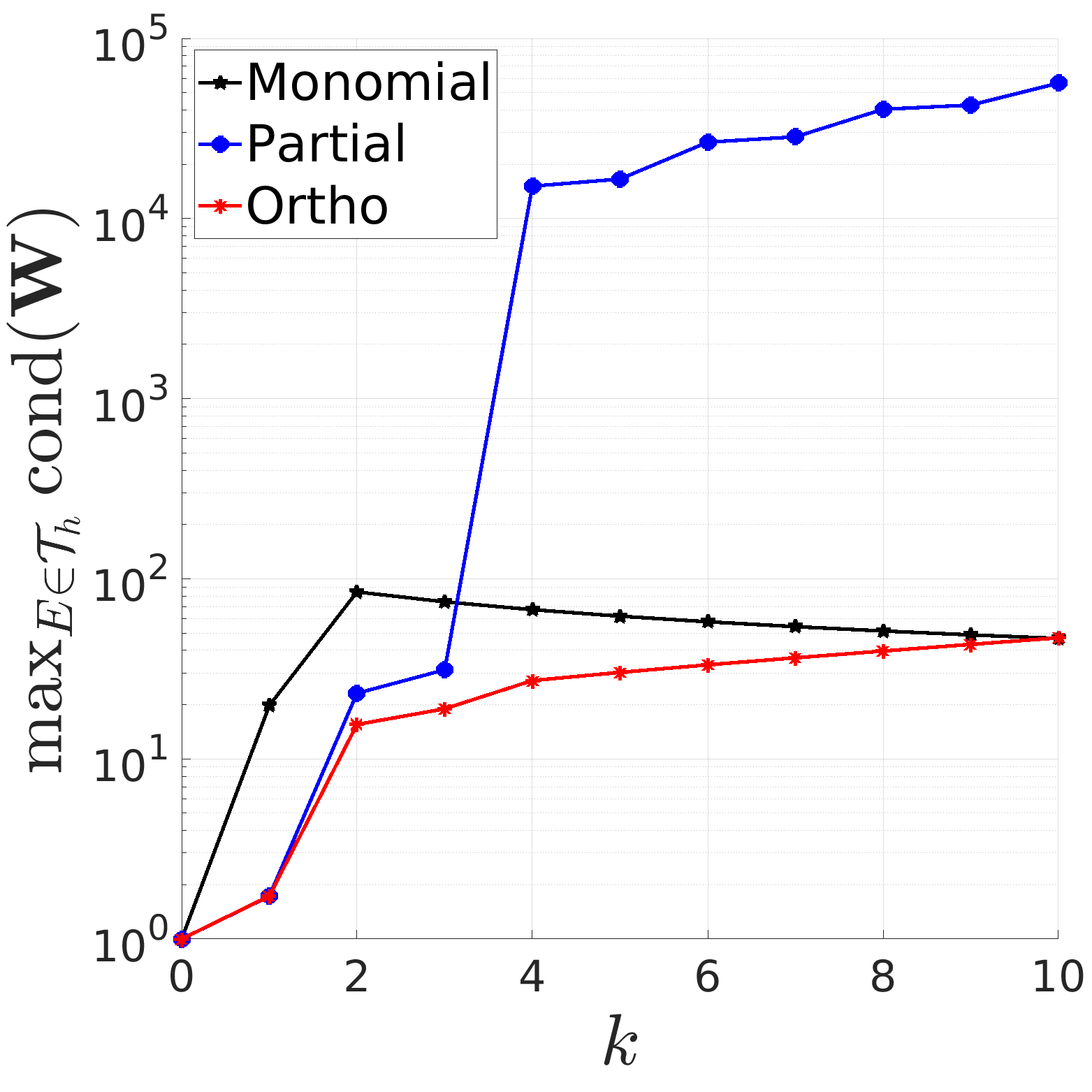}}
	\subfigure[\label{fig:condB_ret10x100}] 
	{\includegraphics[width=.3\textwidth, height = .17\textheight]{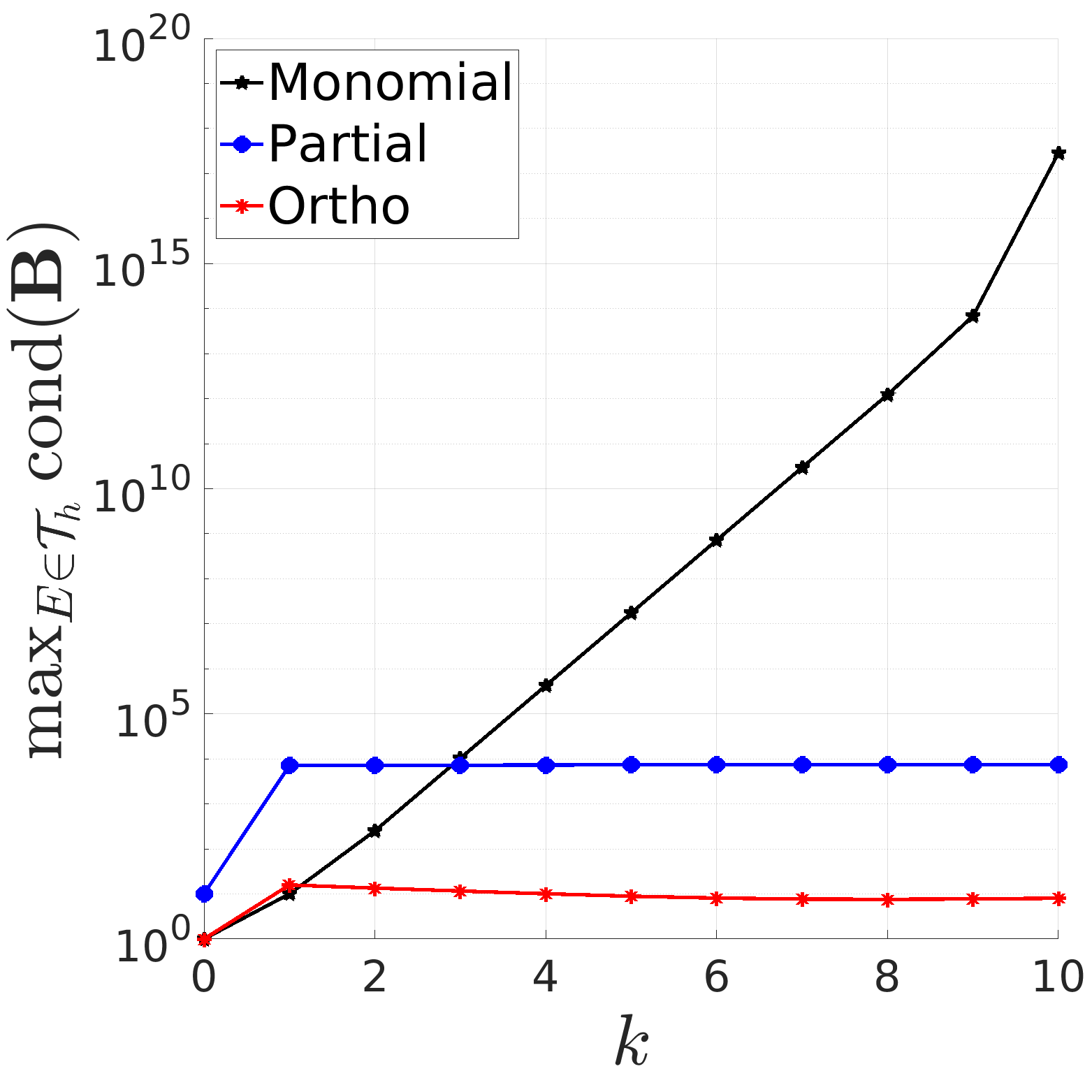}}
	\subfigure[\label{fig:Pi0k_ret10x100}] 
	{\includegraphics[width=.3\textwidth, height = .17\textheight]{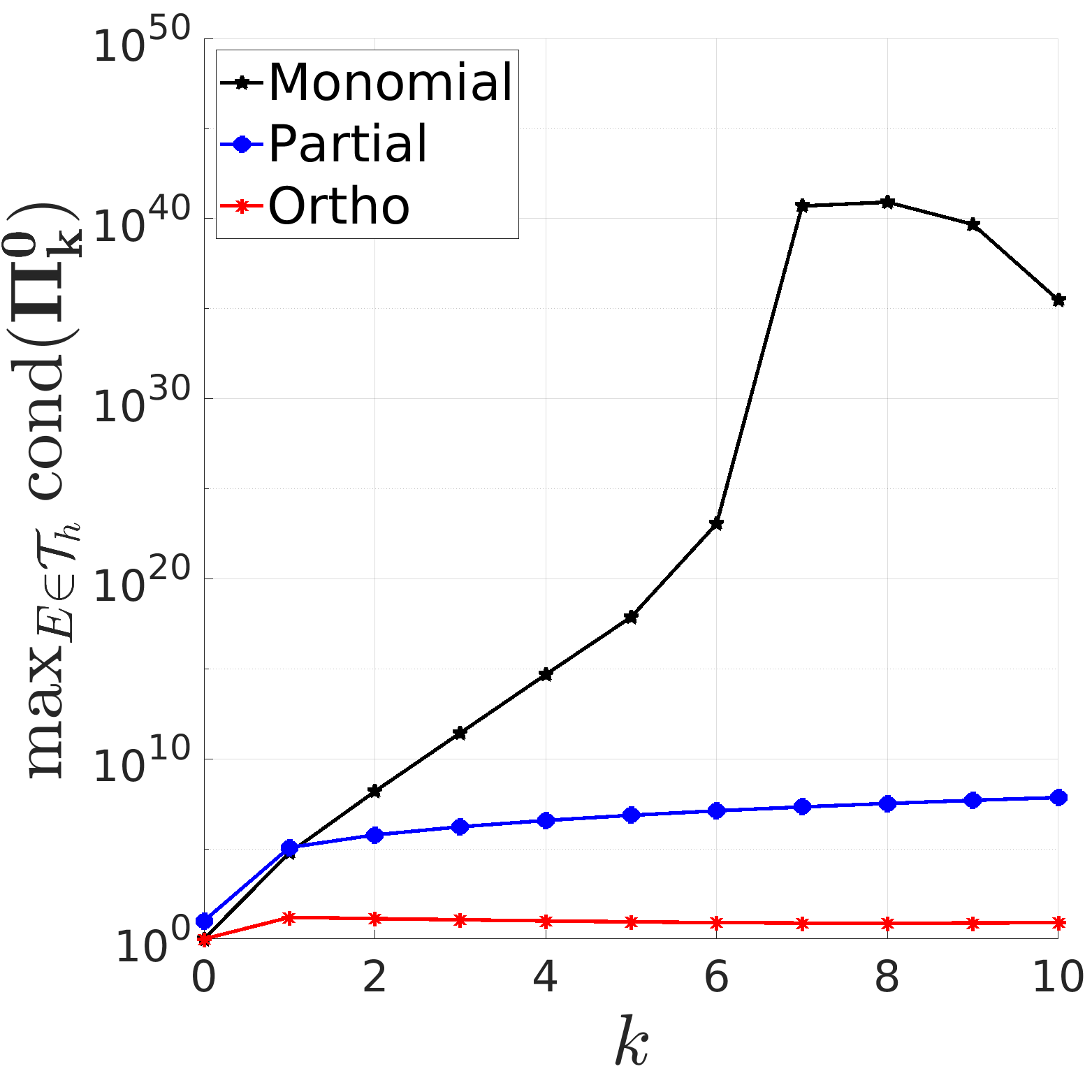}}
	\subfigure[\label{fig:condD_ret10x100}]
	{\includegraphics[width=.3\textwidth, height = .17\textheight]{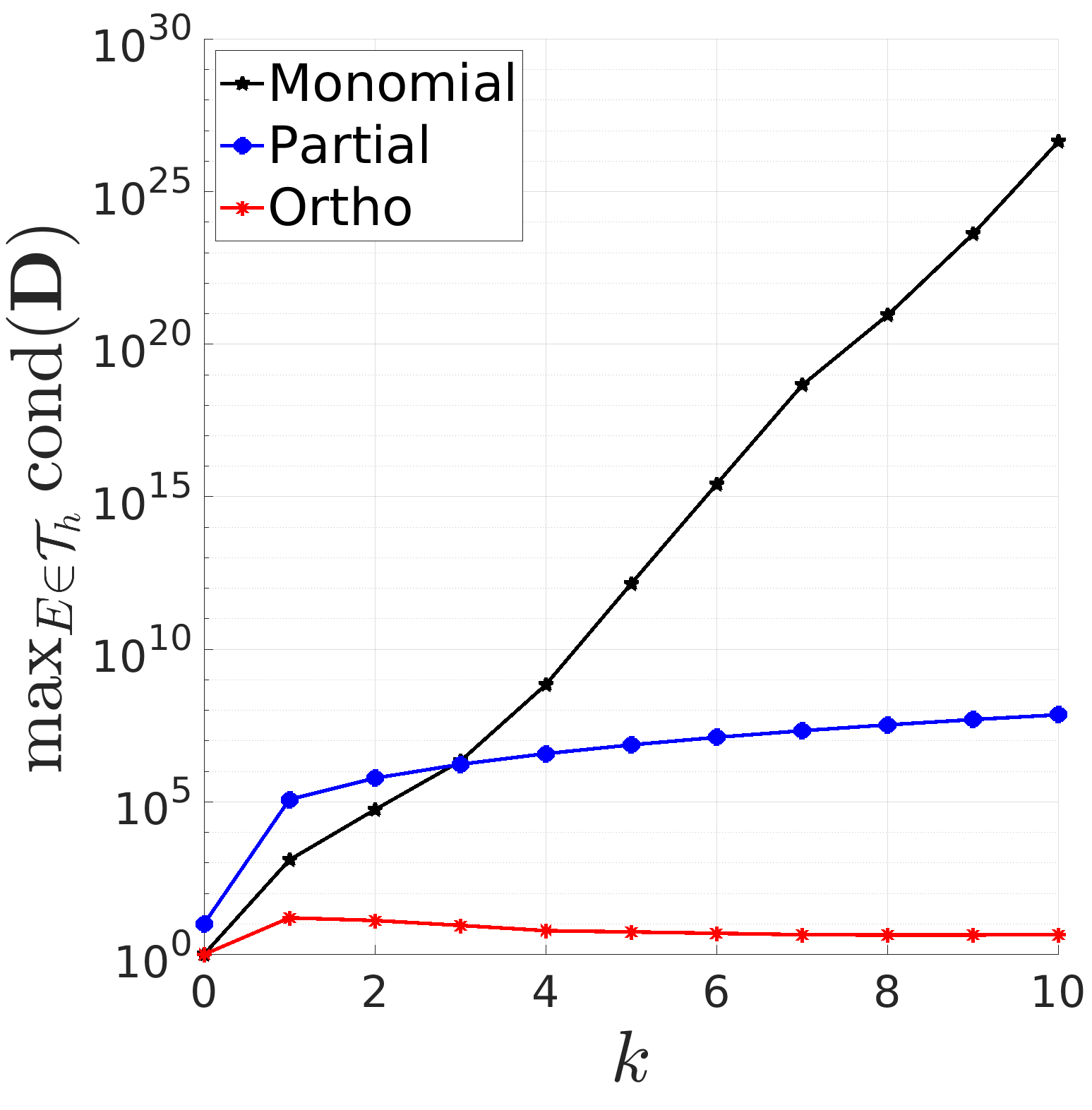}}
    \caption{Test2: Maximum condition number of local matrices among elements, at varying $k$. Mesh  with aspect ratio 10.}
    \label{fig:cond_ret10x100}
\end{figure}

\begin{figure}[htbp]
	\centering
	\subfigure[\label{fig:condG_ret10x500}] 
	{\includegraphics[width=.3\textwidth, height = .17\textheight]{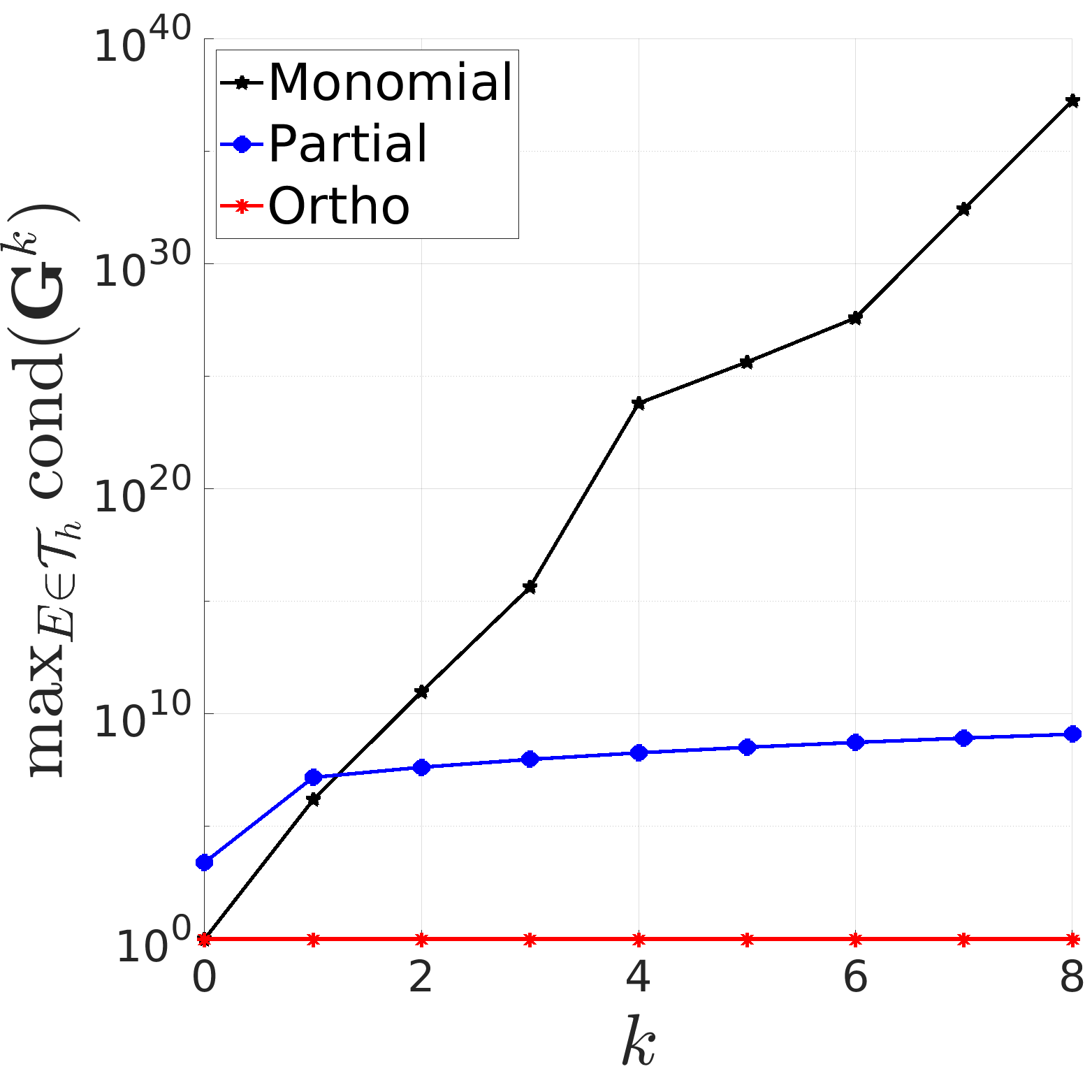}}
	\subfigure[\label{fig:condW_ret10x500}]
	{\includegraphics[width=.3\textwidth, height = .17\textheight]{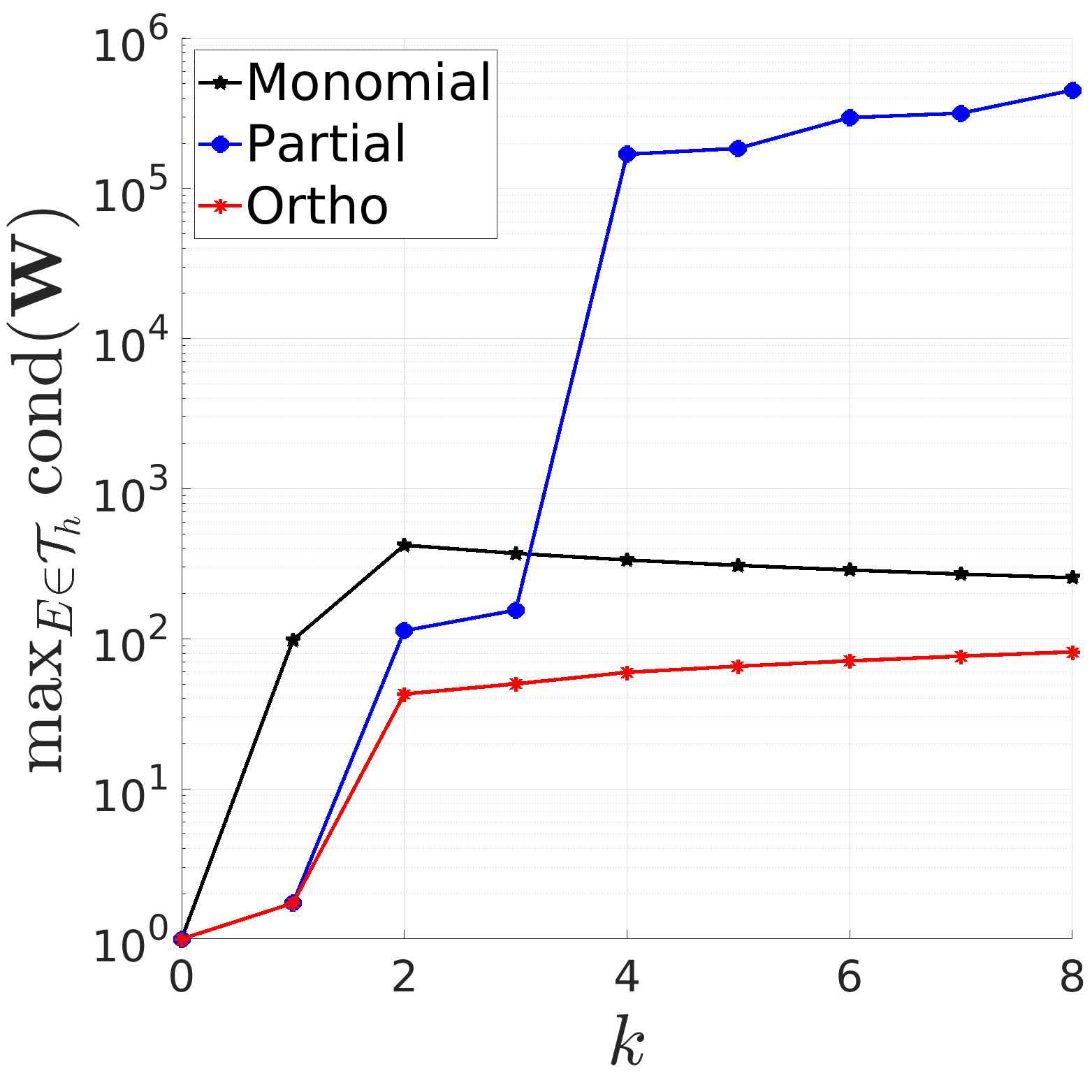}}
	\subfigure[\label{fig:condB_ret10x500}] 
	{\includegraphics[width=.3\textwidth, height = .17\textheight]{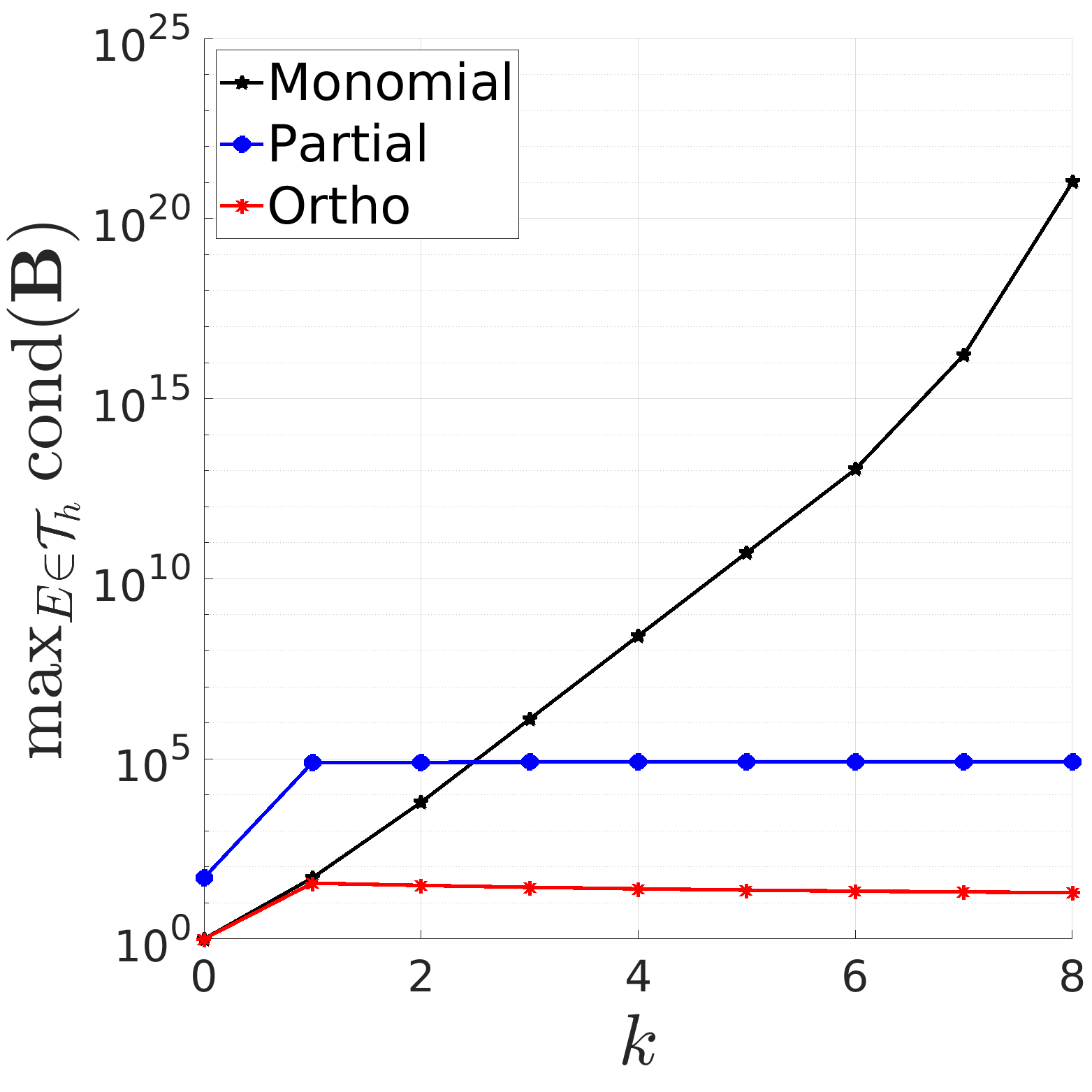}}
	\subfigure[\label{fig:Pi0k_ret10x500}] 
	{\includegraphics[width=.3\textwidth, height = .17\textheight]{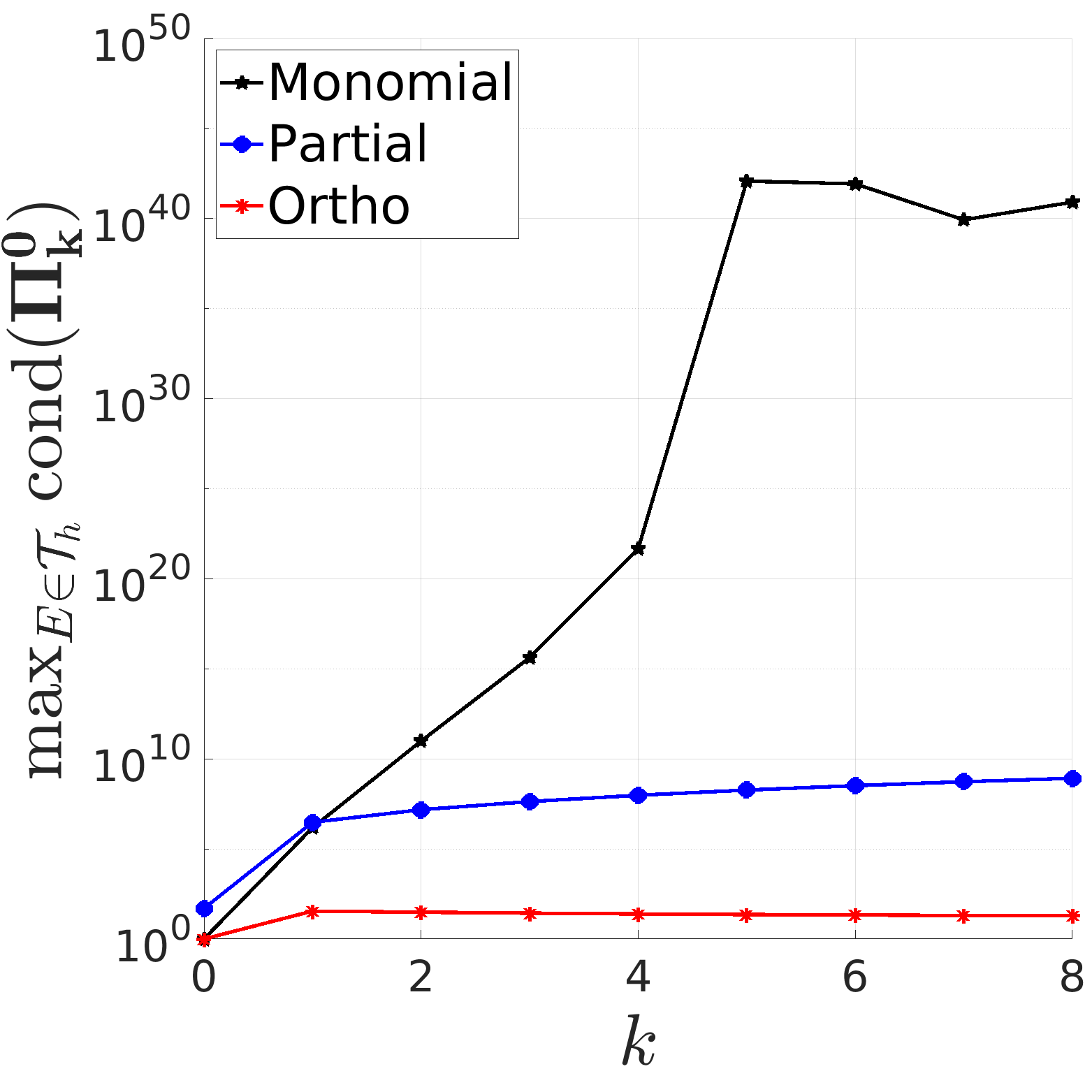}}
	\subfigure[\label{fig:condD_ret10x500}]
	{\includegraphics[width=.3\textwidth, height = .17\textheight]{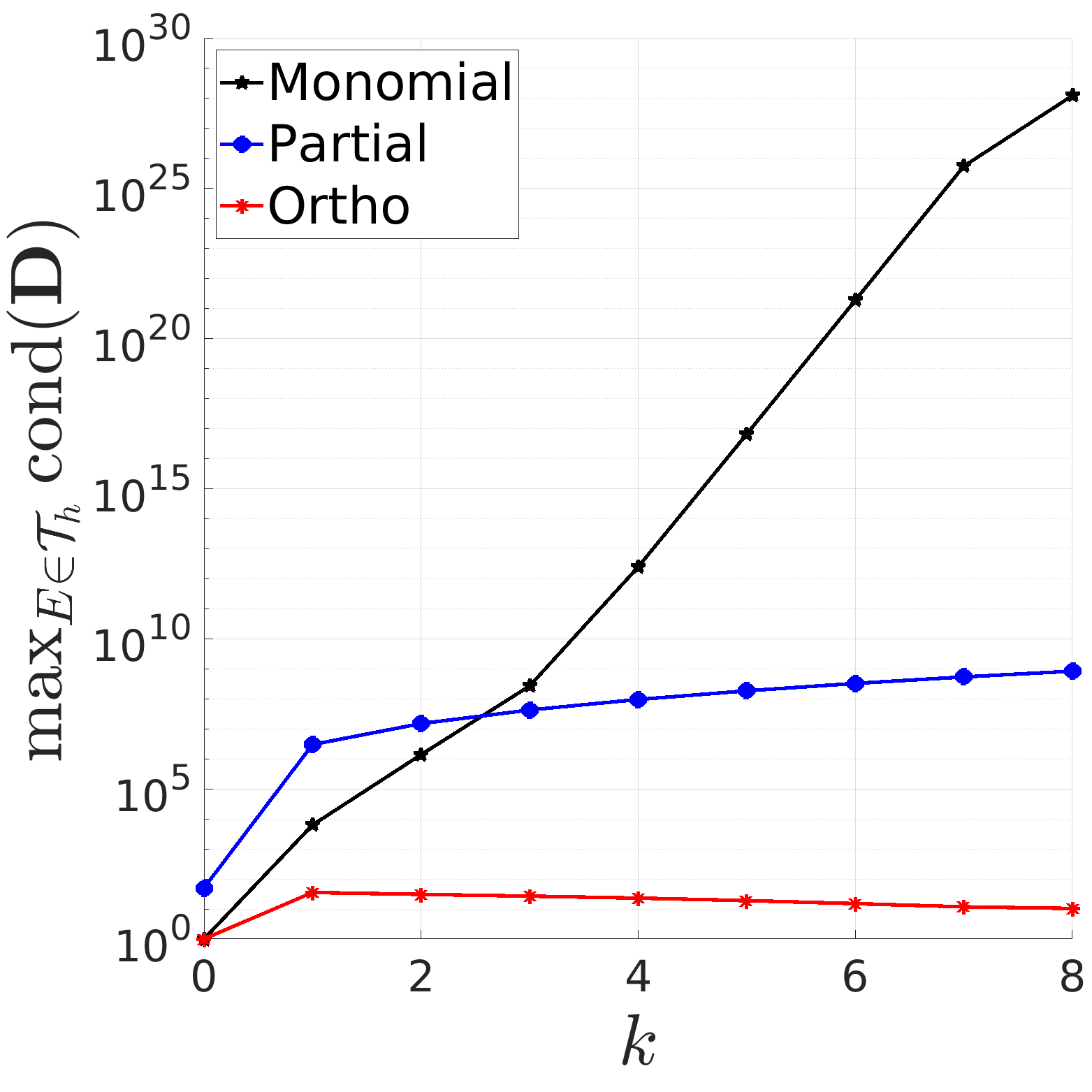}}
    \caption{Test2: Maximum condition number of local matrices among elements at varying $k$. Mesh with aspect ratio 50.}
    \label{fig:cond_ret10x500}
\end{figure}

\begin{figure}[htbp]
	\centering
	\subfigure[\label{fig:condG_ret10x1000}] 
	{\includegraphics[width=.3\textwidth, height = .17\textheight]{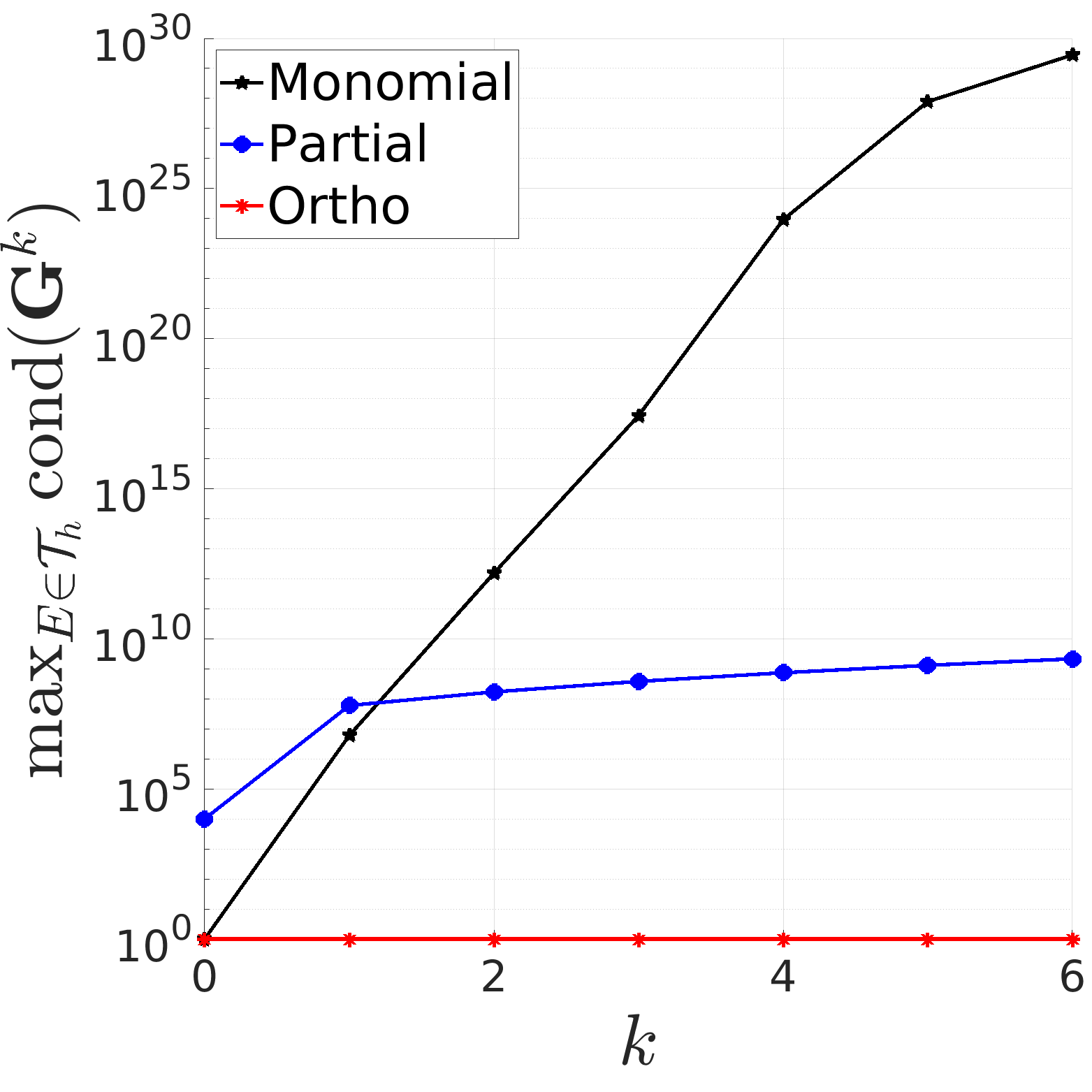}}
	\subfigure[\label{fig:condW_ret10x1000}]
	{\includegraphics[width=.3\textwidth, height = .17\textheight]{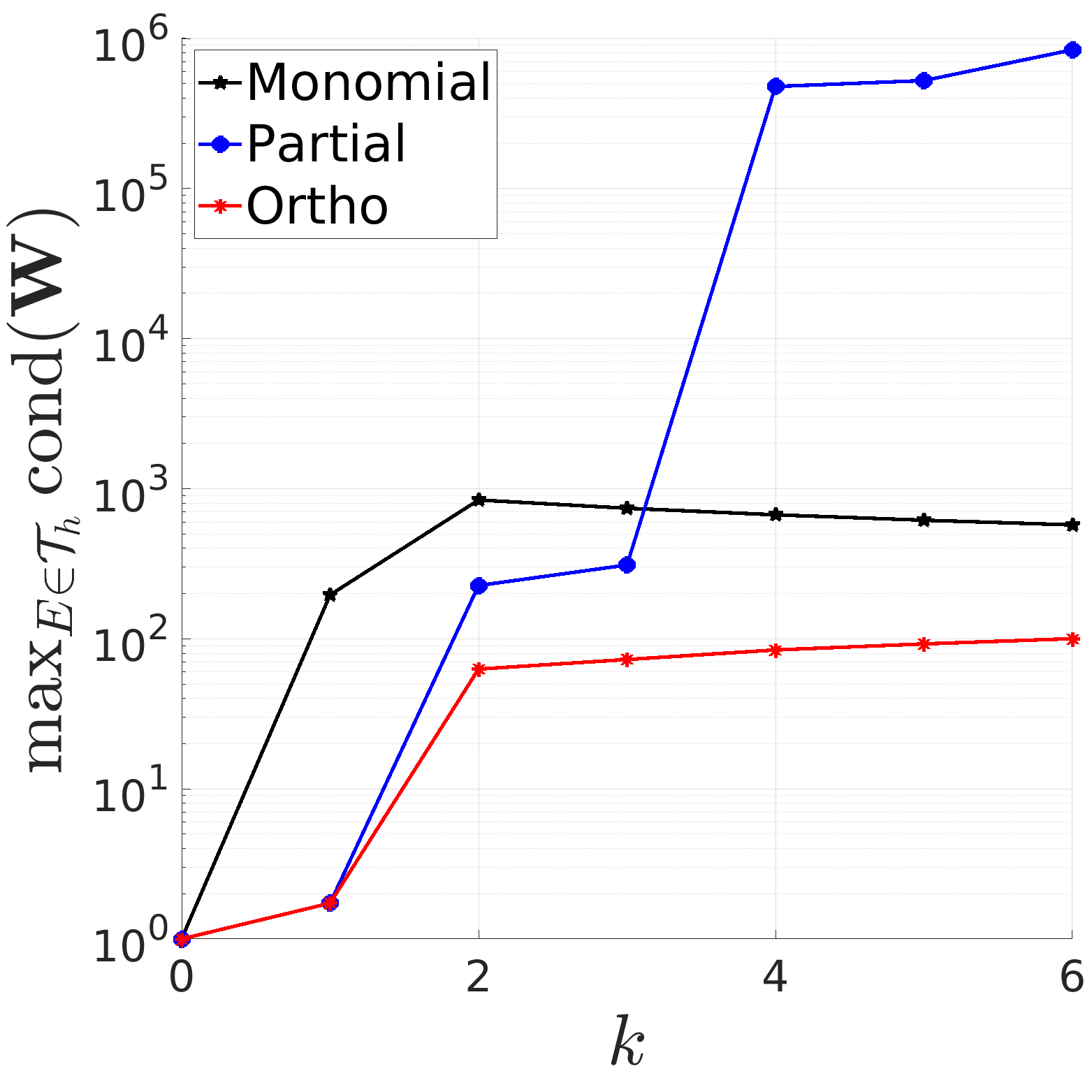}}
	\subfigure[\label{fig:condB_ret10x1000}] 
	{\includegraphics[width=.3\textwidth, height = .17\textheight]{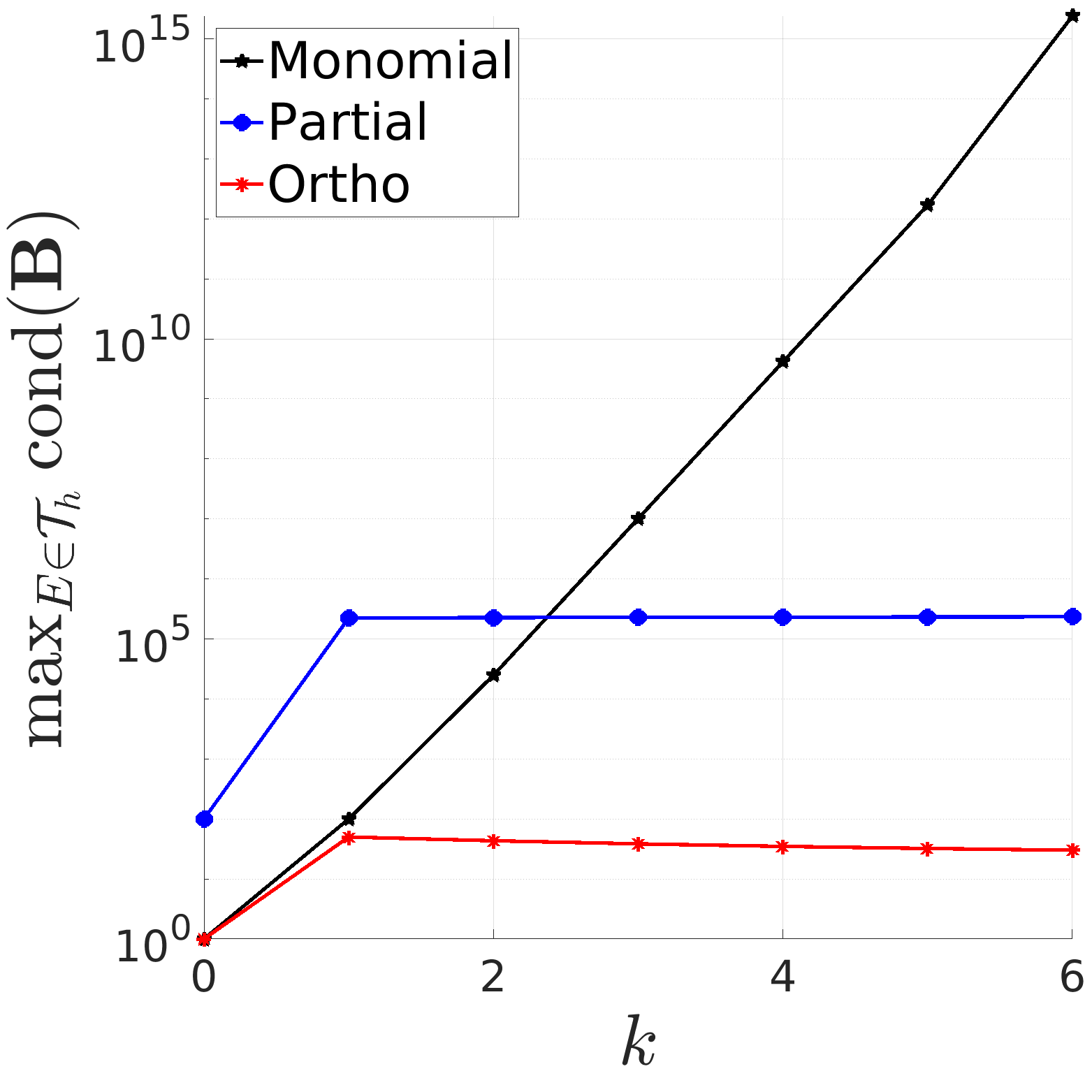}}
	\subfigure[\label{fig:Pi0k_ret10x1000}] 
	{\includegraphics[width=.3\textwidth, height = .17\textheight]{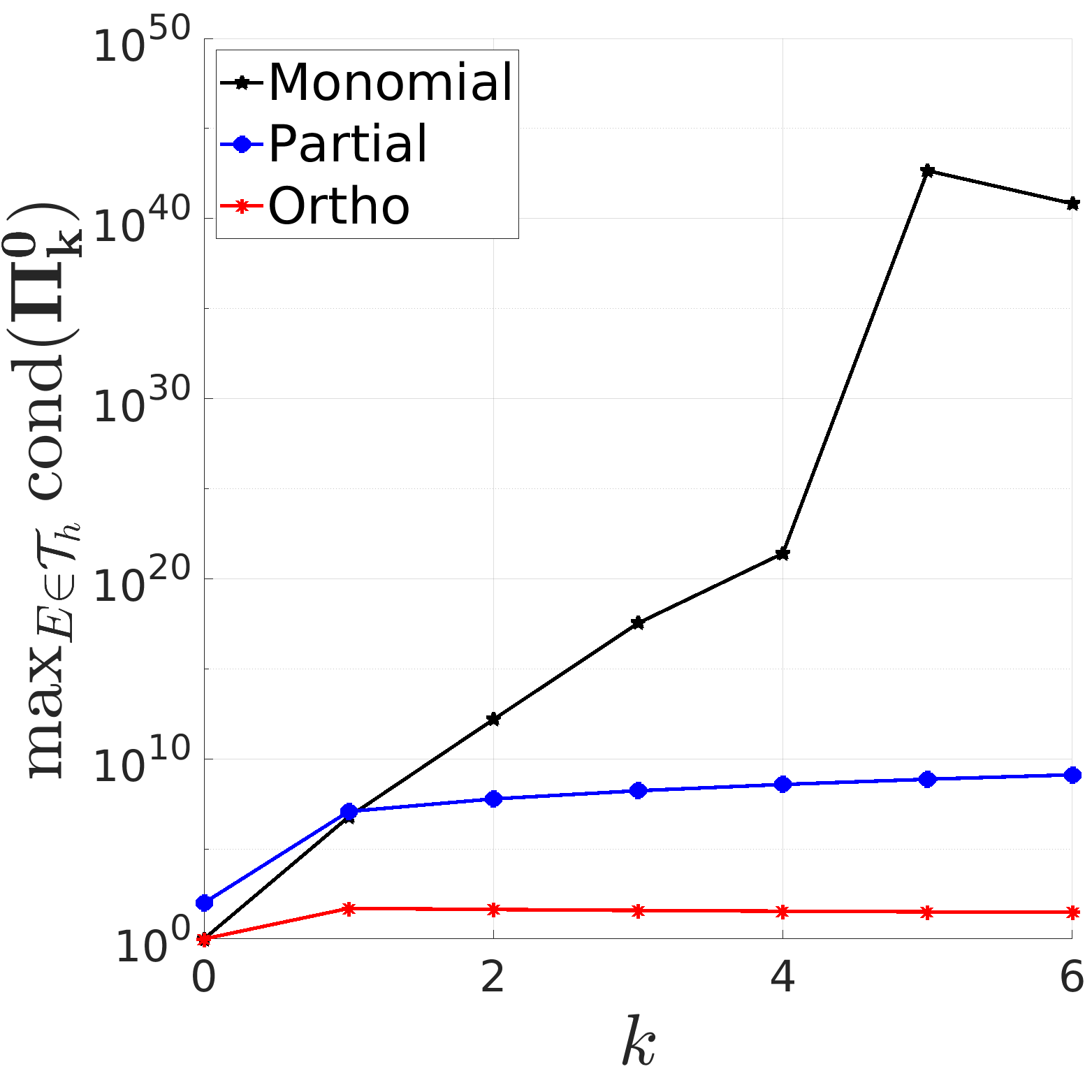}}
	\subfigure[\label{fig:condD_ret10x1000}]
	{\includegraphics[width=.3\textwidth, height = .17\textheight]{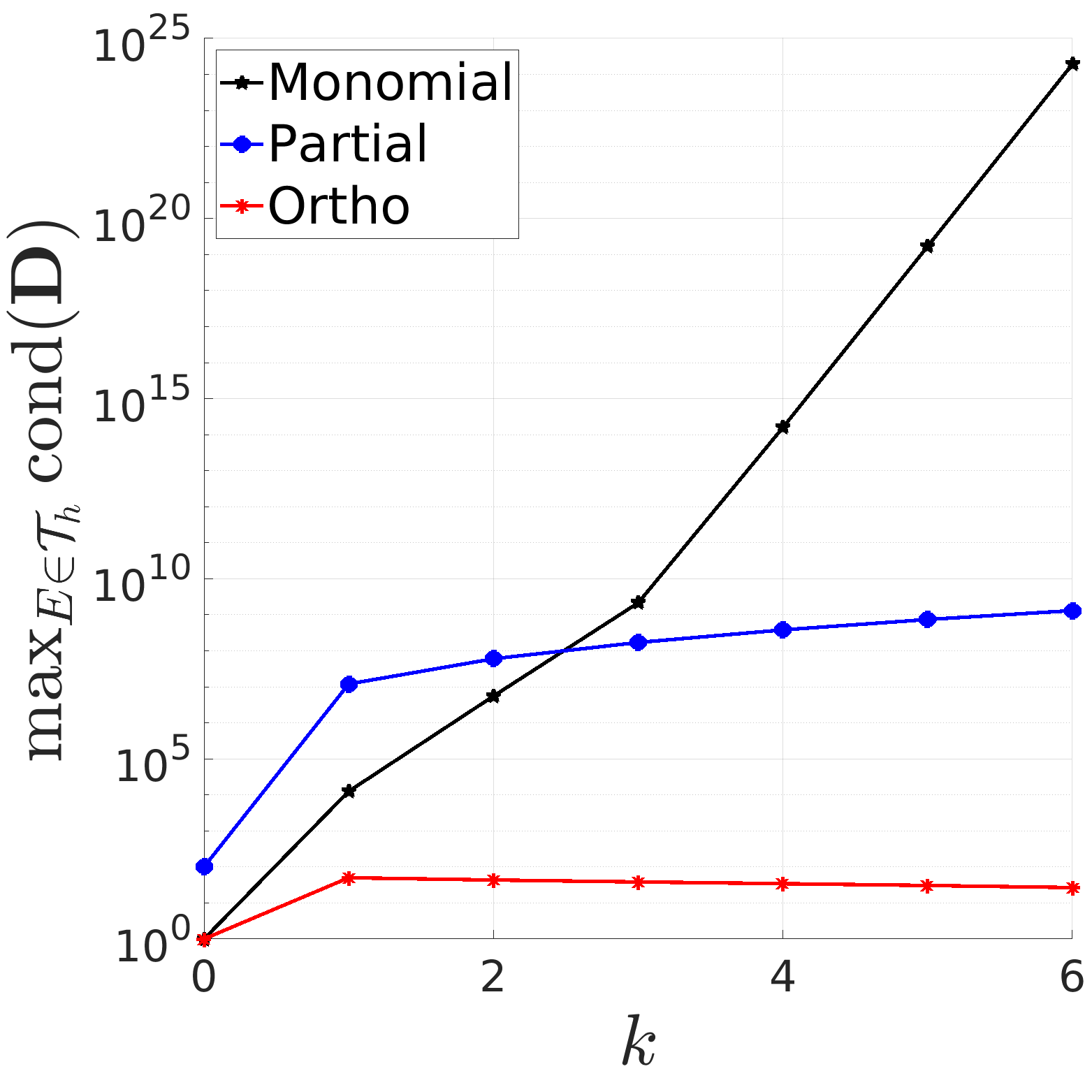}}
    \caption{Test2: Maximum condition number of local matrices among elements at varying $k$. Mesh with aspect ratio 100.}
    \label{fig:cond_ret10x1000}
\end{figure}

Figures \ref{fig:cond_ret10x100}, \ref{fig:cond_ret10x500} and \ref{fig:cond_ret10x1000} report the highest value among mesh elements of the condition number of elemental matrices $\mathbf{G}^k$, $\mathbf{W}$, $\mathbf{B}$, $\mathbf{\PPi^0_k}$ and $\mathbf{D}$, on the three considered meshes. 
As shown in these figures, mass matrix $\mathbf{G}^k$, resulting from the full orthonormal approach, is perfectly well-conditioned, since it corresponds to the identity matrix, up to machine precision. Also, an algebraic growth of the condition number of the mass matrix obtained with the partial orthonormal approach is observed, whereas conditioning grows exponentially for the monomial one. 
More generally, we can observe that the condition numbers of the local matrices vary in a very limited range in the partial and full orthonormal approaches, as $k$ increases. Instead, the condition number of local matrices grows exponentially in the monomial approach, with the only exception of matrices $\mathbf{W}$, which however have good conditioning in all cases.  

Conditioning of matrices obtained with the partial-orthonormal bases appear to be slightly more affected by increasing aspect ratio values than the corresponding matrices with the full-orthonormal approach, however, still showing much lower values than the matrices given by the monomial basis.

\begin{figure}[htbp]
	\centering
	\subfigure[\label{fig:ErrorL2Pressure_ret10x100}] 
	{\includegraphics[width=.3\textwidth, height = .17\textheight]{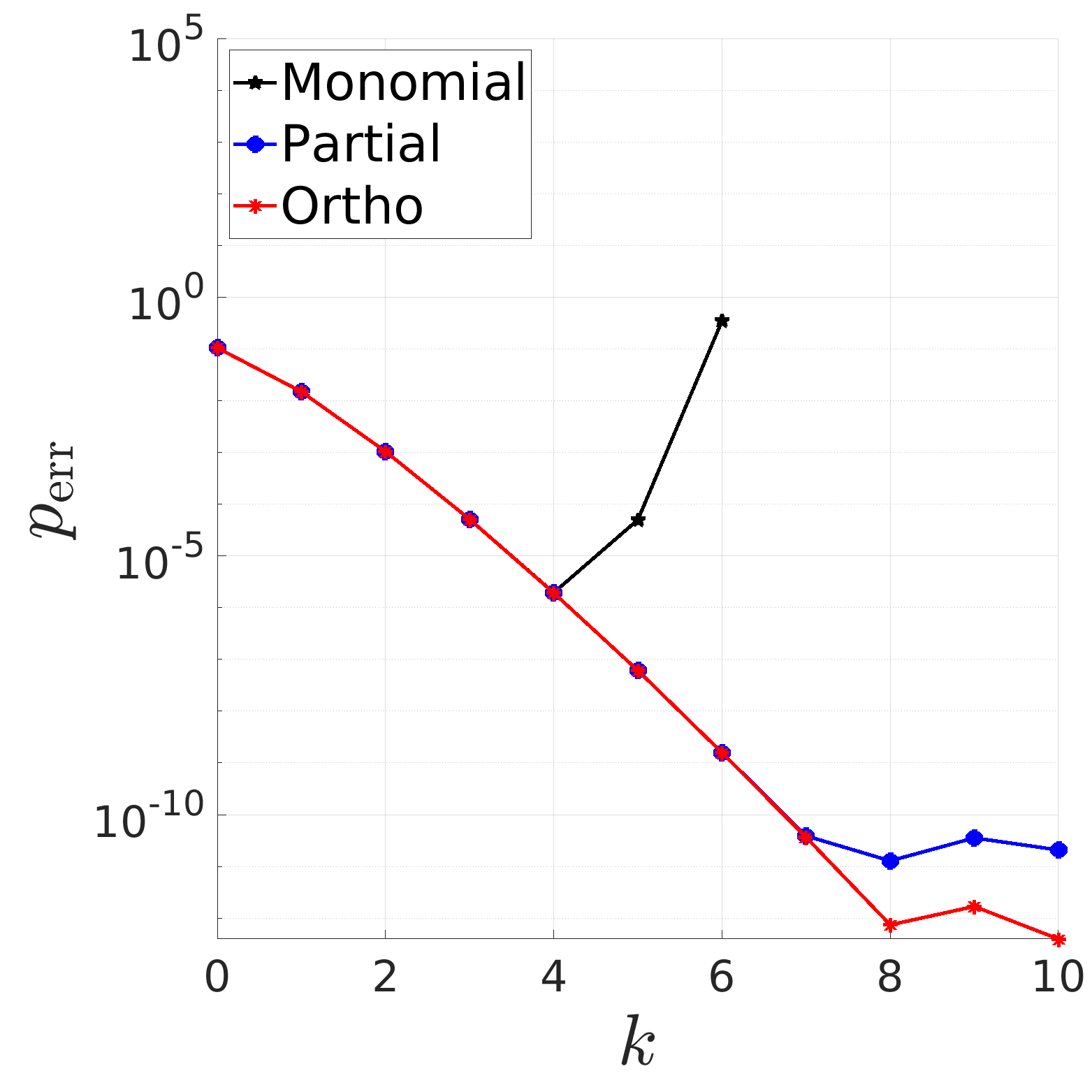}}
	\subfigure[\label{fig:ErrorL2Velocity_ret10x100}] 
	{\includegraphics[width=.3\textwidth, height = .17\textheight]{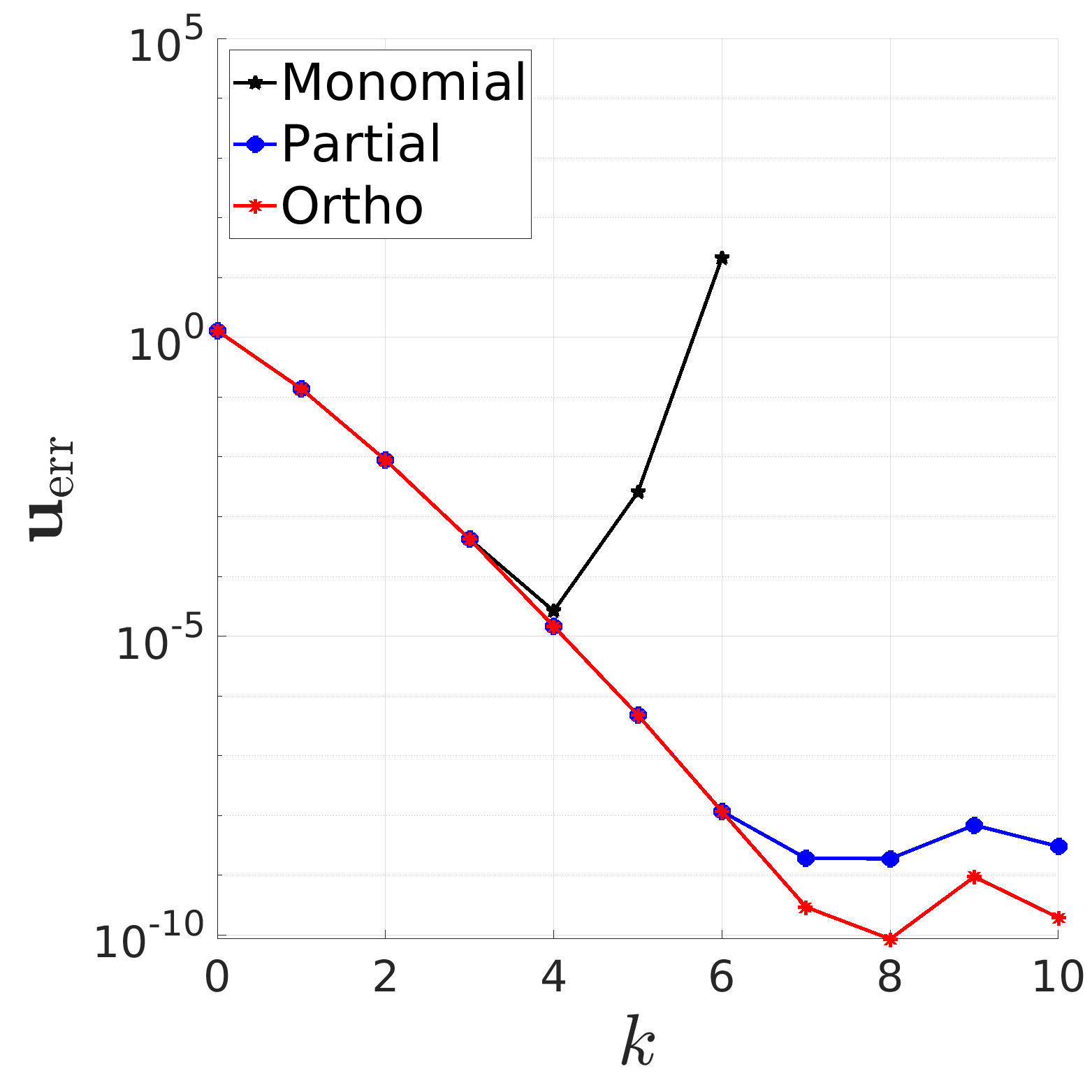}}
	\subfigure[\label{fig:SuperConvergence_ret10x100}] 
	{\includegraphics[width=.3\textwidth, height = .17\textheight]{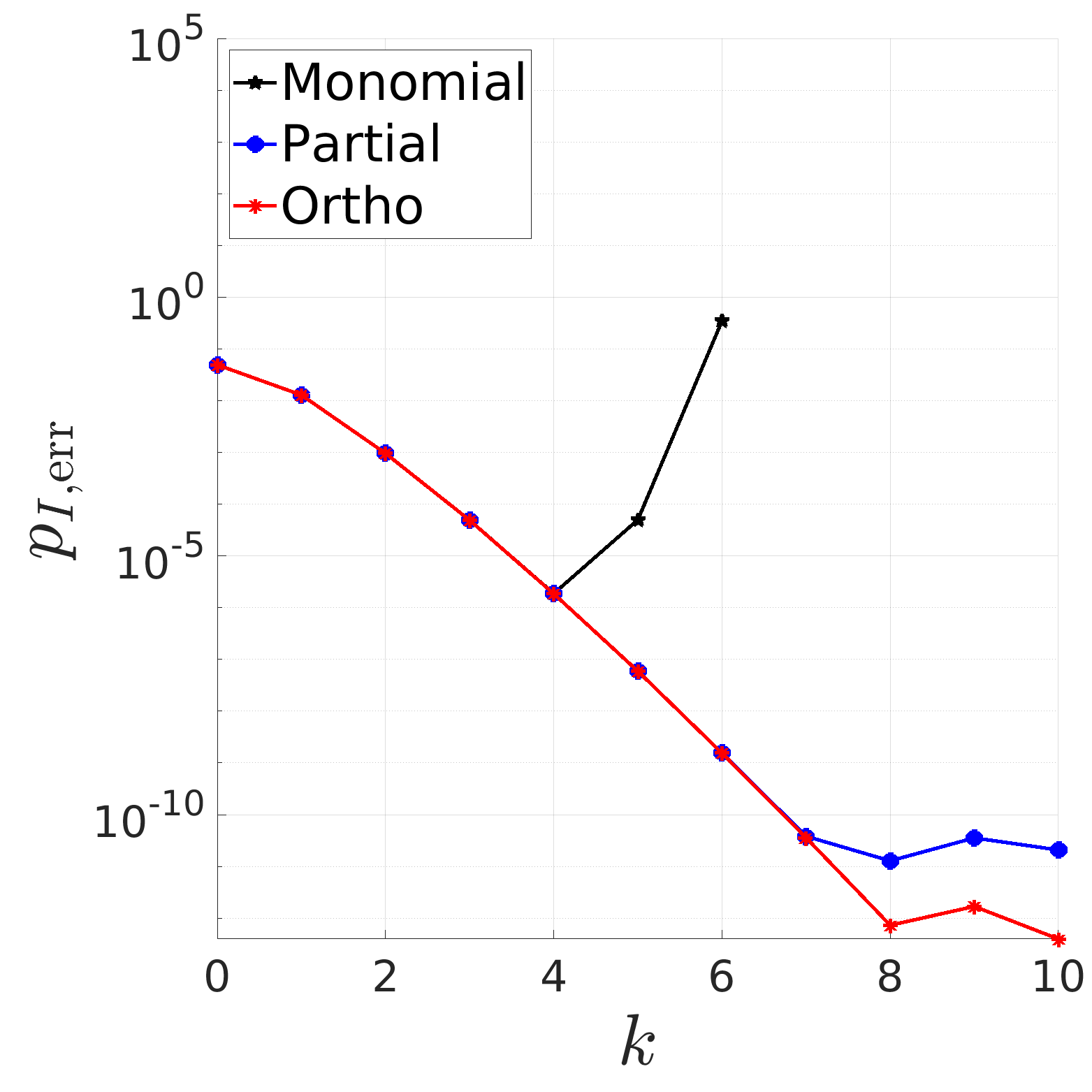}}
	\subfigure[\label{fig:ErrorL2Pressure_ret10x500}] 
	{\includegraphics[width=.3\textwidth, height = .17\textheight]{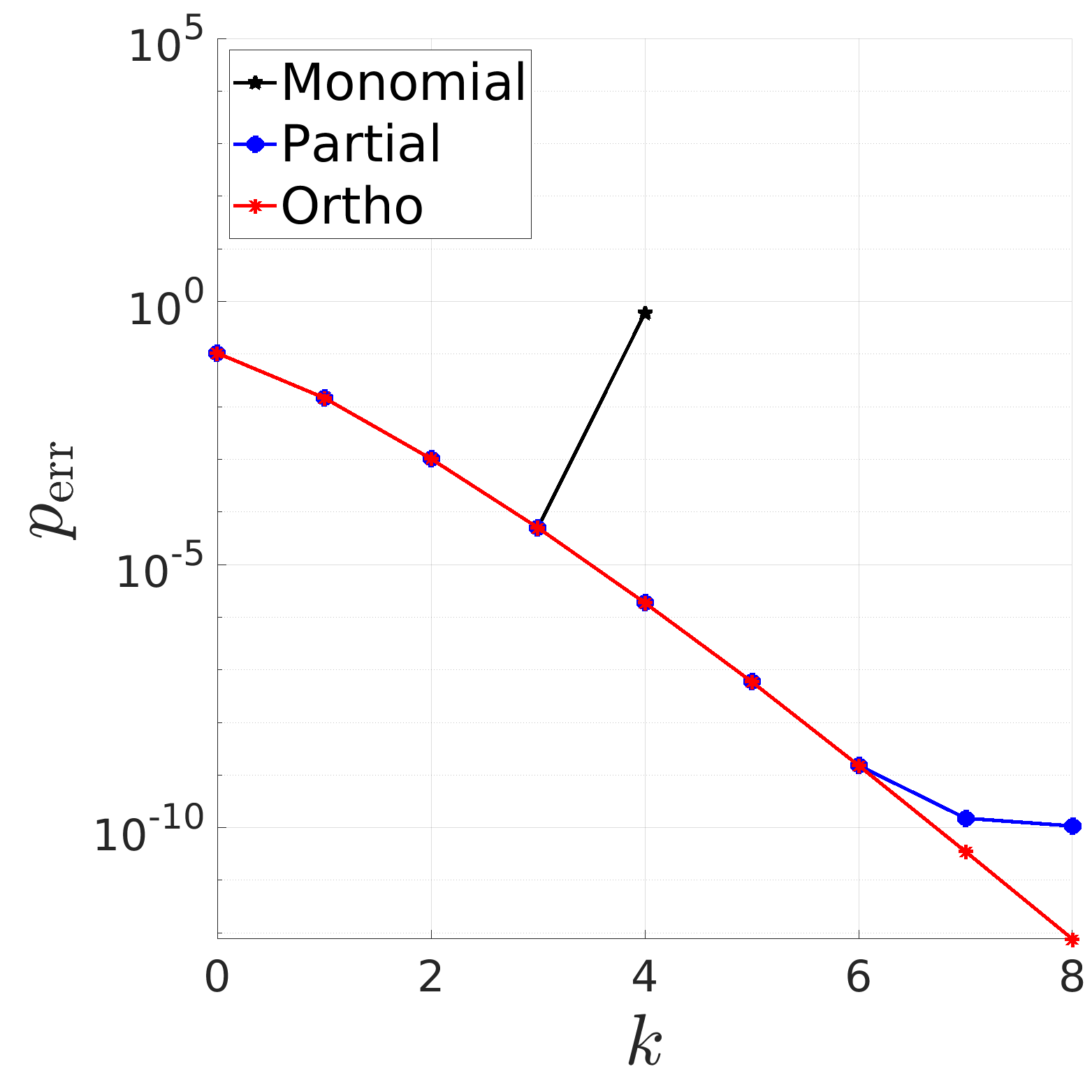}}
	\subfigure[\label{fig:ErrorL2Velocity_ret10x500}] 
	{\includegraphics[width=.3\textwidth, height = .17\textheight]{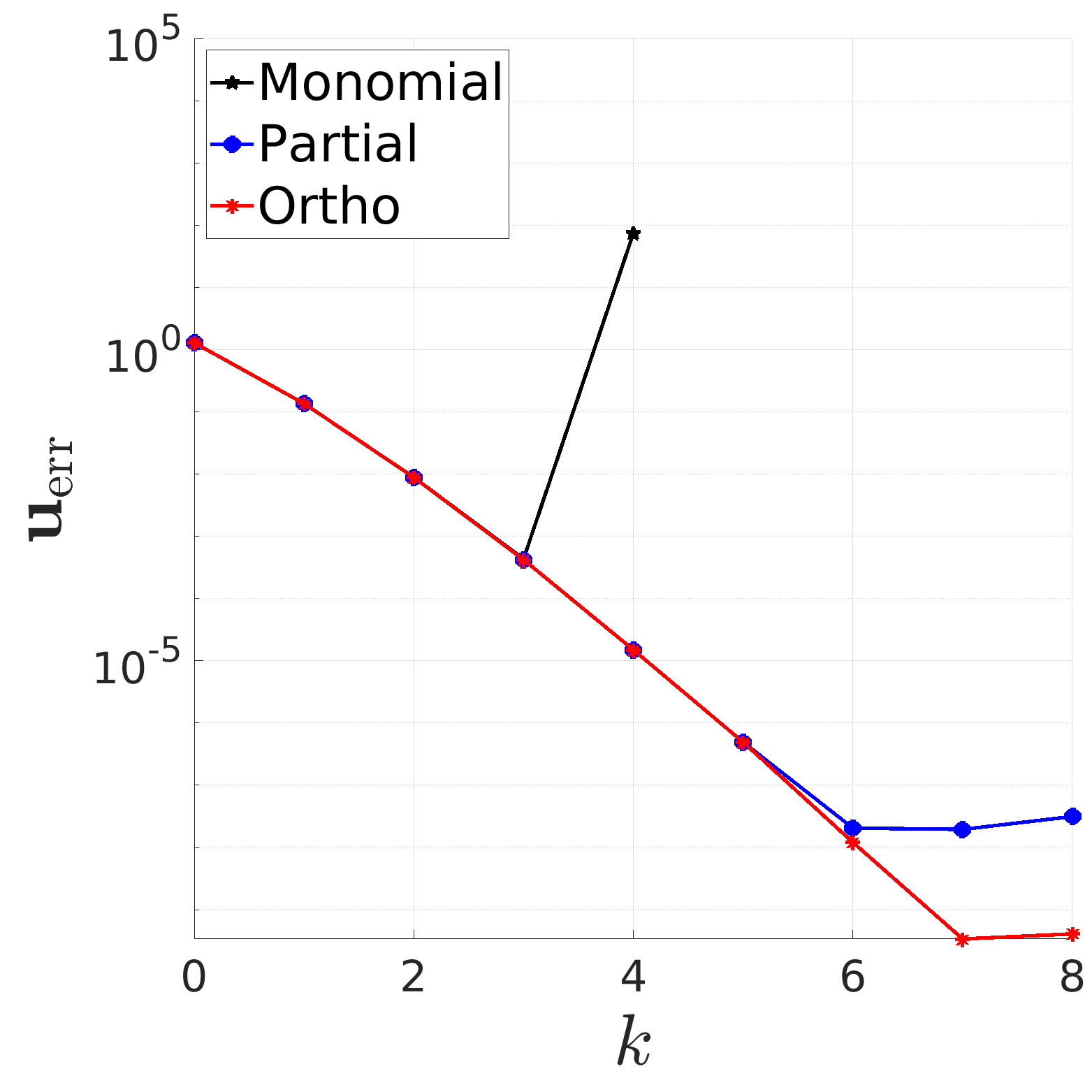}}
	\subfigure[\label{fig:SuperConvergence_ret10x500}] 
	{\includegraphics[width=.3\textwidth, height = .17\textheight]{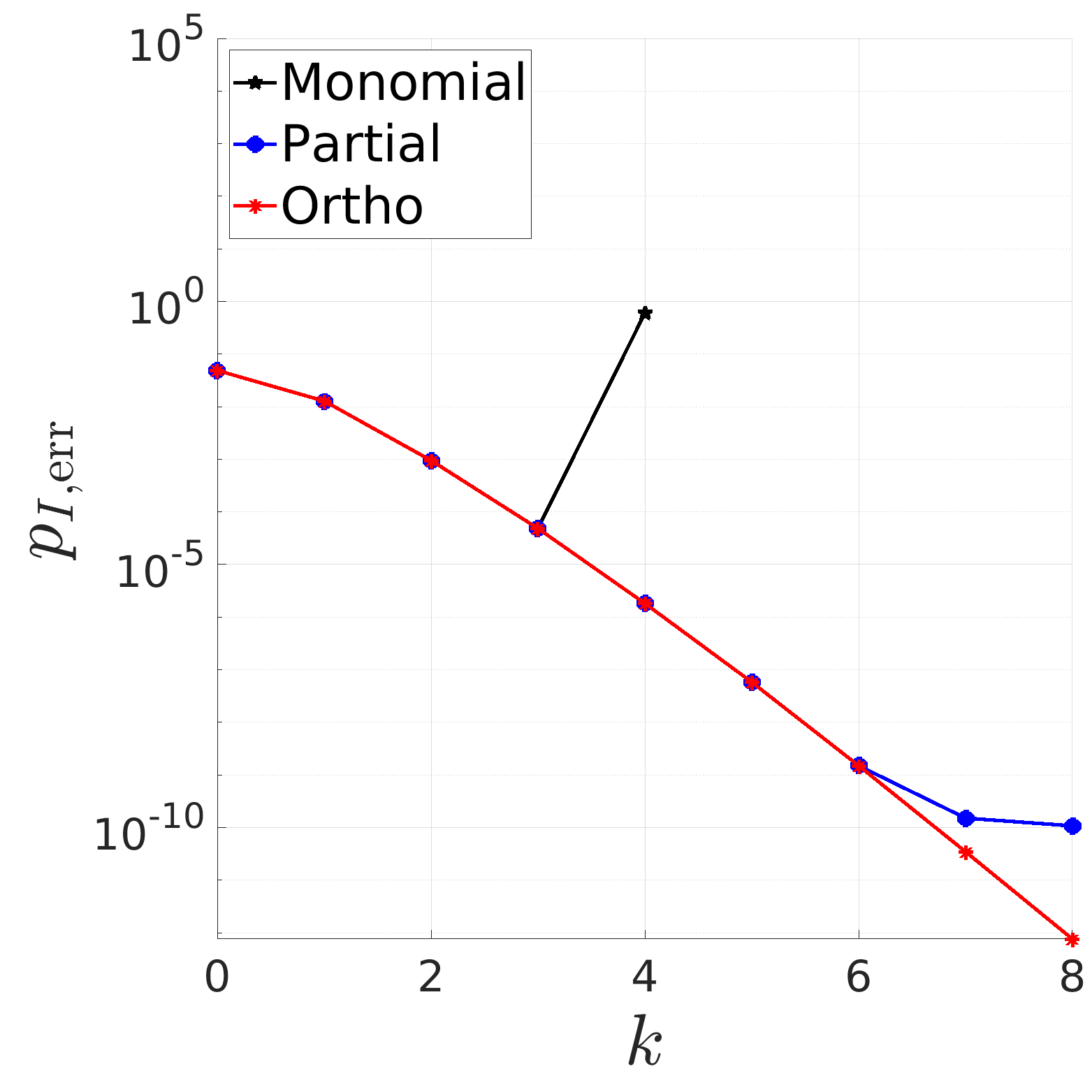}}
	\subfigure[\label{fig:ErrorL2Pressure_ret10x1000}] 
	{\includegraphics[width=.3\textwidth, height = .17\textheight]{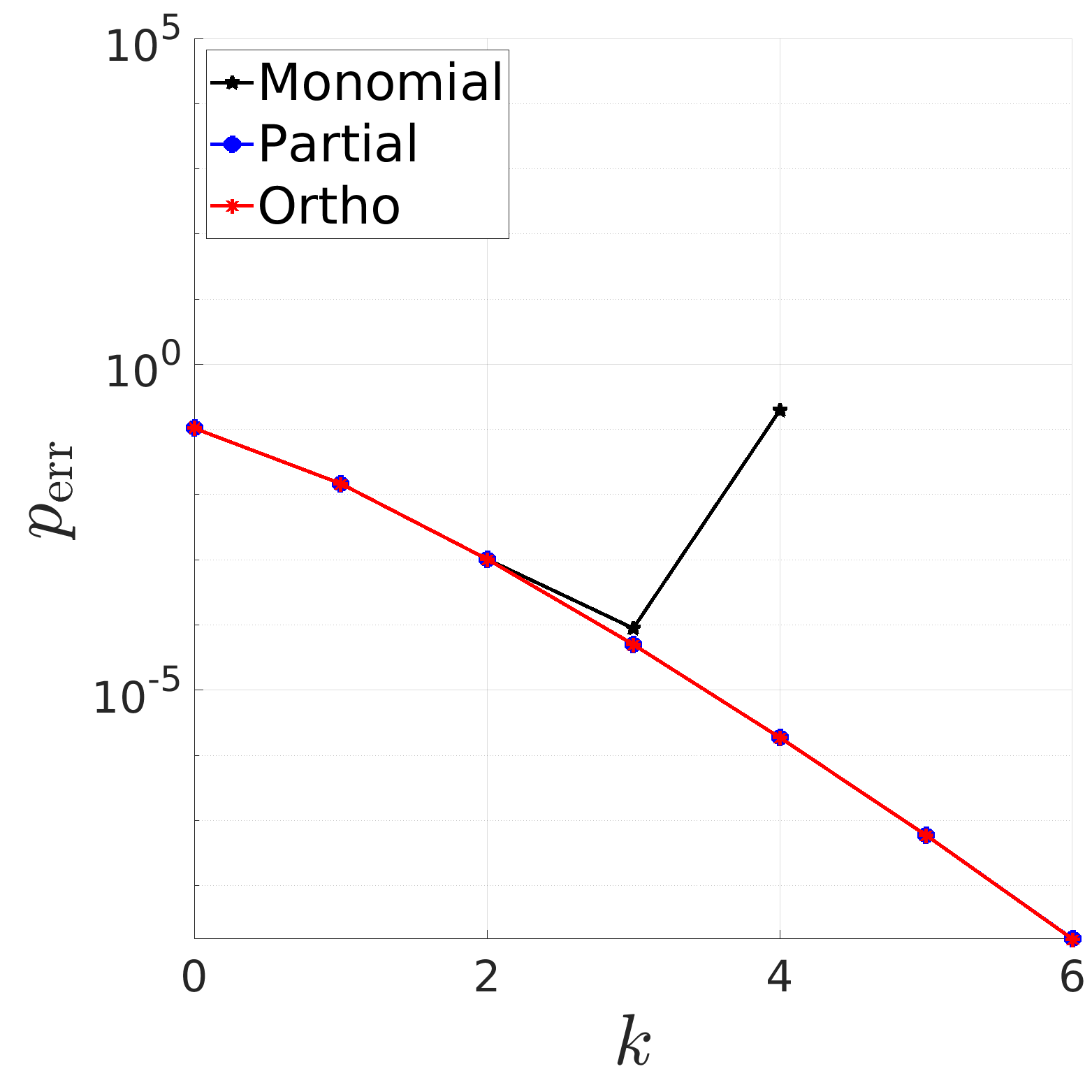}}
	\subfigure[\label{fig:ErrorL2Velocity_ret10x1000}] 
	{\includegraphics[width=.3\textwidth, height = .17\textheight]{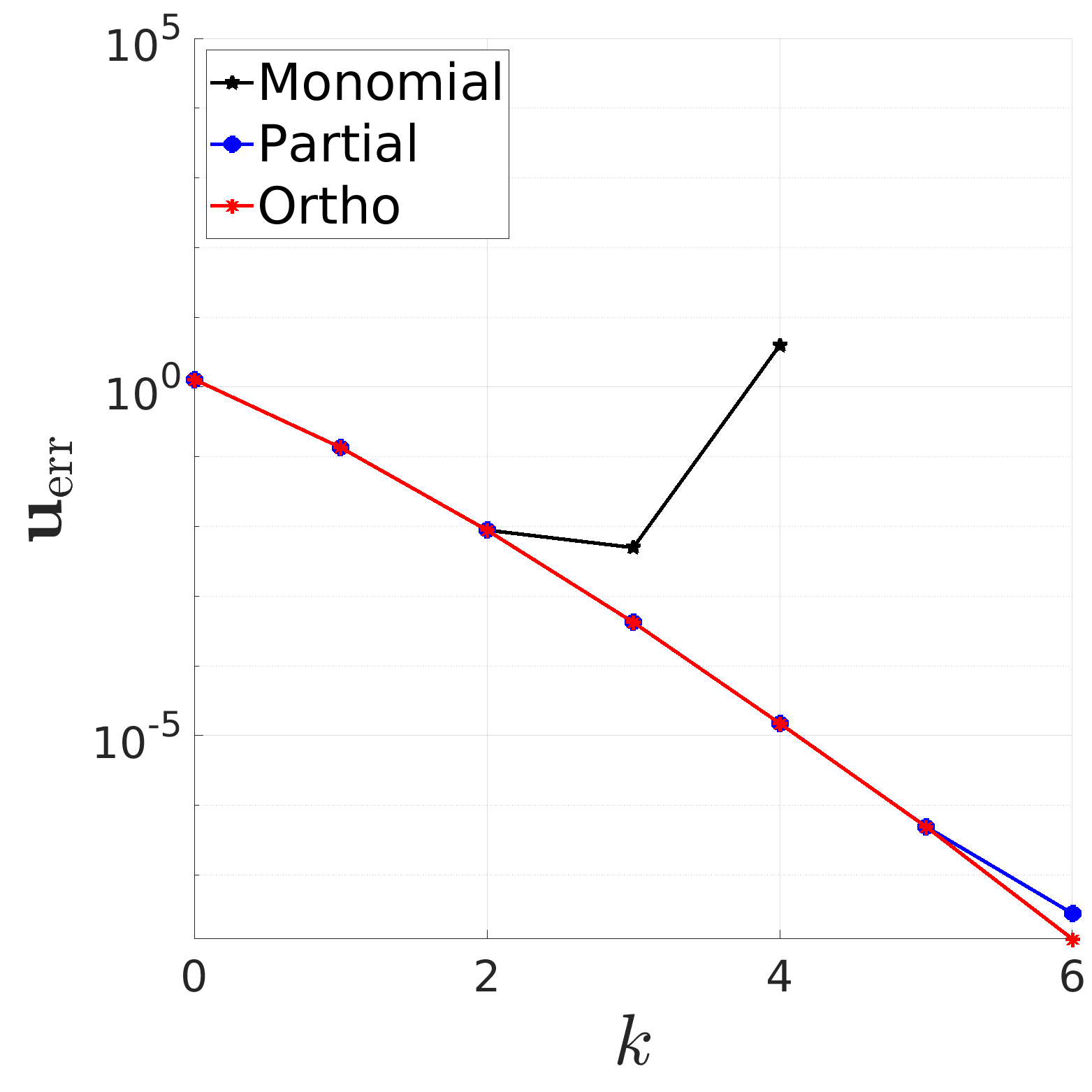}}
	\subfigure[\label{fig:SuperConvergence_ret10x1000}] 
	{\includegraphics[width=.3\textwidth, height = .17\textheight]{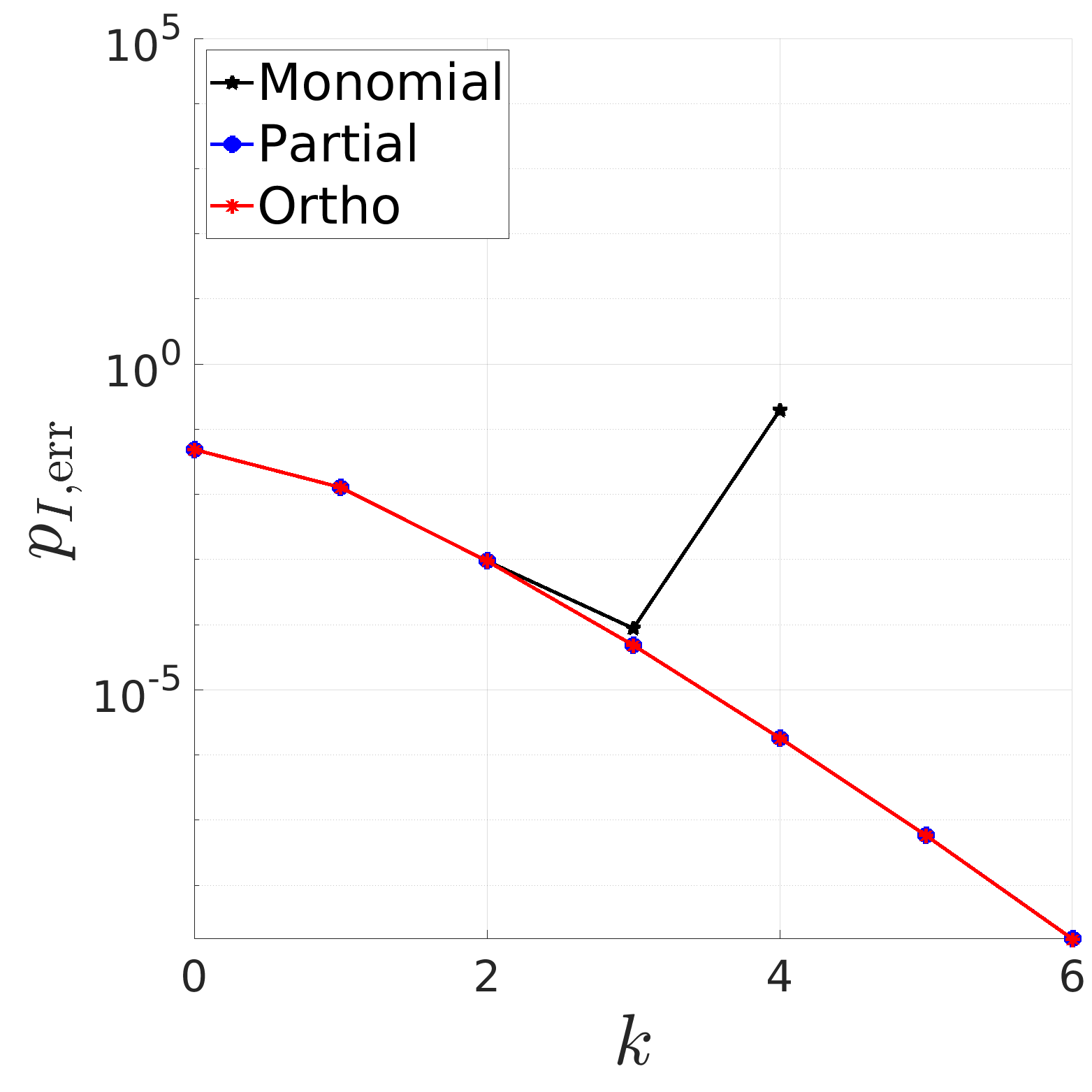}}
    \caption{Test2: Behaviour of errors \eqref{eq:L2pressure}, \eqref{eq:L2velocity} and \eqref{eq:superconvergence} at varying $k$ on rectangular meshes. Each row represents a different mesh: 10, 50, 100 from top to bottom.}
    \label{fig:errors_secondExp}
\end{figure}

Figure \ref{fig:errors_secondExp} reports convergence curves of the error against growing polynomial accuracy $k$, for the three considered meshes: top row for the mesh with aspect ratio 10, middle row for aspect ratio 50 and bottom row for aspect ratio 100.
At low values of $k$, the curves corresponding to the monomial approach are well overlapped with the curves of partial and full-orthonormal approaches. Further, the maximum value of polynomial accuracy $k$ for which the monomial approach provides errors in line with the other approaches reduces as the aspect ratio of mesh elements increases. Finally, error curves relating to the monomial approach are interrupted at values of $k\leq 6$ due to failure of linear algebra libraries (the \textit{SparseLU} solver of Eigen) in computing a solution due to ill-conditioning.
Small differences are instead noticed between the curves obtained with the partial and full orthonormal approaches for all the considered values of $k$.

\subsection{Test3: Simulations in Discrete Fracture Networks}
In this last example, the application of the proposed orthonormal bases to Discrete Fracture Network (DFN) problems is presented. Discrete Fracture Networks are obtained as the union of planar polygonal domains with arbitrary orientations in the 3D space and are used to model the fractures in a porous medium \cite{DFNbook}. As fracture thickness is typically orders of magnitude smaller than the other dimensions, fractures are geometrically reduced to 2D domains, and suitable equations, averaged across fracture thickness, are derived to describe the phenomena occurring in such domains \cite{MJR2005}. Interface equations are then added at fracture intersections. A major complexity in DFN simulations consists in the generation of a conforming mesh, for realistic configurations characterized by intricate networks with a large number of fractures and fracture intersections. A possible strategy for DFN meshing is proposed, e.g. in \cite{DFN}, based on the use of mixed virtual elements: first a triangular mesh is constructed on each fracture domain, independently of the intersections with the other fractures; then, the elements of these meshes are cut according to the interfaces, and hanging nodes are added where needed, to obtain a fully conforming mesh of the whole domain (see \cite{DFN} for more details). Highly elongated elements are likely to be generated in this process, such that the standard choice of VEM basis function might yield badly conditioned problems for high order approximations \cite{CurvedMVEMAppl}.

An advection-diffusion-reaction problem on a network made up by the union of three fractures $F_i$ with three intersections $\Gamma_i$, is here considered, namely

\begin{align*}
    F_1 &= \{(x,y,z) \in \mathbb{R}^3: -1 \leq x \leq 1, -1 \leq y \leq 1, z = 0\},\\
    F_2 &= \{(x,y,z) \in \mathbb{R}^3: -1 \leq z \leq 1, -1 \leq x \leq 1, y = 0\},\\
    F_3 &= \{(x,y,z) \in \mathbb{R}^3: -1 \leq y \leq 1, -1 \leq z \leq 1, x = 0\},
\end{align*}
\begin{align*}
    \Gamma_1 &= \{(x,y,z) \in \mathbb{R}^3: -1 \leq x \leq 1, y = 0, z = 0\},\\
    \Gamma_2 &= \{(x,y,z) \in \mathbb{R}^3: -1 \leq y \leq 1, x = 0, z = 0\},\\
    \Gamma_3 &= \{(x,y,z) \in \mathbb{R}^3: -1 \leq z \leq 1, x = 0, y = 0\}.
\end{align*}
On each fracture, we choose 

\begin{equation*}
    \D_i(\hat{x},\hat{y}) = \begin{bmatrix}
    1 + \hat{y}^2 & -\frac{\hat{x}\hat{y}}{2}\\
    -\frac{\hat{x}\hat{y}}{2} & 1+ \hat{x}^2
    \end{bmatrix},\quad \bb_i(\hat{x},\hat{y}) = \begin{bmatrix}
    \hat{x} - \hat{y}\\ \hat{y}-1    \end{bmatrix} , \quad \gamma_i(\hat{x},\hat{y}) = \hat{x}^3 + \hat{y}.
\end{equation*}
where $(\hat{x},\hat{y})$ is a proper fracture-local reference system, and a forcing term and Neumann boundary conditions are defined in such a way the exact solution on each fracture is:

\begin{align*}
    h_1(x,y) &= -\vert x \vert (1+x)(1-x)y(1+y)(1-y) ,\\
    h_2(z,x) &= -z(1+z)(1-z)x(1+x)(1-x) ,\\
    h_3(y,z) &= y(1+y)(1-y)\vert z \vert (1+z)(1-z).
\end{align*}

\begin{figure}[htbp]
	\centering
	{\includegraphics[width=.8\textwidth, height = .5\textheight]{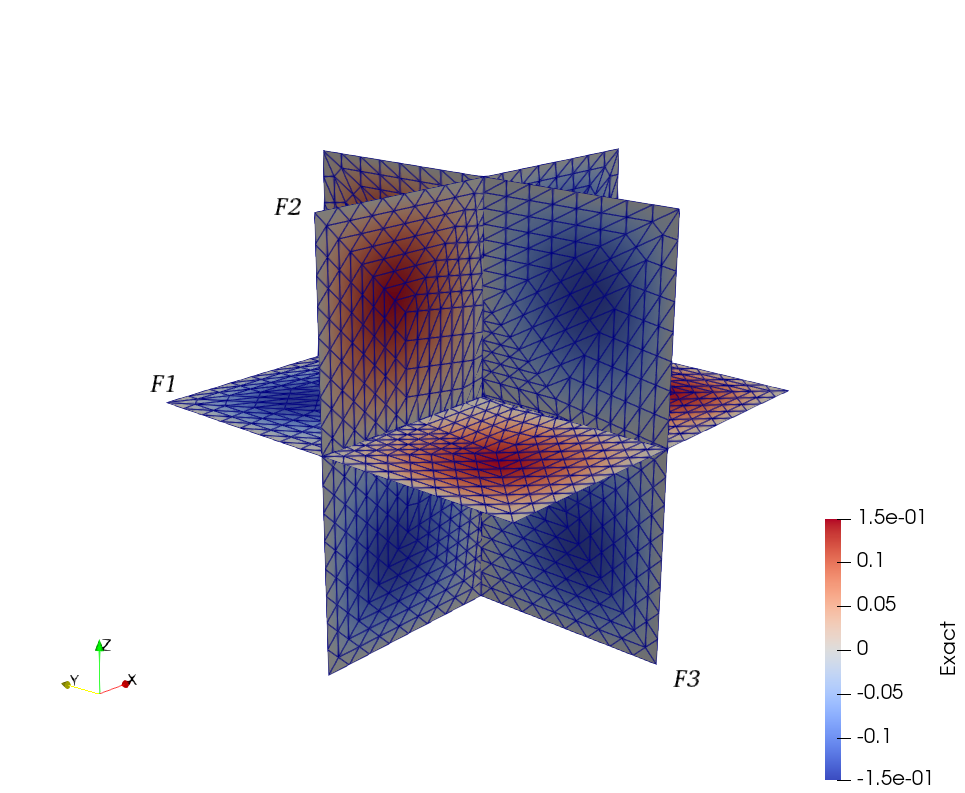}}
    \caption{Test3: Exact solution DFN benchmark problem.}
    \label{fig:sol_3expe}
\end{figure}

The same problem is considered in \cite{DFN} up to order $5$. The network and the exact solution are shown in Figure \ref{fig:sol_3expe}.

In this numerical test, we use four refinements of an initially triangular mesh, modified, as mentioned above, in such a way that the final polygonal meshes are conforming at the traces. Table~\ref{tab:DFN_mesh} reports the number, the minimum and maximum aspect ratio of mesh elements for the four refinements (R0 to R3).

\begin{table}[]\footnotesize
\centering
\begin{tabular}{cccc}
\hline
 & \textbf{Num Cells} & \textbf{Min AR} & \textbf{Max AR} \\
 \hline 
 \textbf{R0} & 75 & 1.41 & 38.43
 \\
 \textbf{R1} & 246 & 1.19 & 62.11 \\
 \textbf{R2} & 882 & 1.19 & 152.11 \\
 \textbf{R3} & 3294 & 1.14 & 215.13
\\
\hline
\end{tabular}
\caption{Test3: number of cells (Num Cells), minimum aspect ratio (Min AR) and maximum aspect ratio (Max AR) of mesh cells for the four refinement levels (R0 to R3).}
\label{tab:DFN_mesh}
\end{table}

\begin{table}[]\footnotesize
\centering
\resizebox{\textwidth}{!}{
\begin{tabular}{lccccccccl}
\hline
                  & \textbf{$k$}         & \textbf{0}                 & \textbf{1}                 & \textbf{2}                 & \textbf{3}                 & \textbf{4}                 & \textbf{5}                 & \textbf{6}                 & \multicolumn{1}{c}{\textbf{7}} \\ \hline
\textbf{Monomial} & $\pp_{\text{err}}$   & 1.1984                     & 2.3292                     & 3.5494                     & 4.6377                     & 5.6597                     & 2.8423                     & 1.8227                     & -4.7912                        \\
                  & $\uu_{\text{err}}$   & \multicolumn{1}{l}{1.0354} & \multicolumn{1}{l}{1.8430} & \multicolumn{1}{l}{2.8787} & \multicolumn{1}{l}{3.9900} & \multicolumn{1}{l}{4.9328} & \multicolumn{1}{l}{1.1135} & \multicolumn{1}{l}{0.6710} & -6.5224                        \\
\textbf{}         & $\pp_{I,\text{err}}$ & 1.7581                     & 2.6572                     & 3.8542                     & 4.8566                     & 5.8788                     & 2.8333                     & 1.8227                     & -4.7912                        \\ \hline
\textbf{Partial}  & $\pp_{\text{err}}$   & 1.1984                     & 2.3292                     & 3.5494                     & 4.6377                     & 5.6633                     & 6.8989                     & 5.7860                     & -2.8821                        \\
                  & $\uu_{\text{err}}$   & \multicolumn{1}{l}{1.0354} & \multicolumn{1}{l}{1.8430} & \multicolumn{1}{l}{2.8787} & \multicolumn{1}{l}{3.9900} & \multicolumn{1}{l}{5.0341} & \multicolumn{1}{l}{4.7002} & \multicolumn{1}{l}{2.8836} & -4.0248                        \\
\textbf{}         & $\pp_{I,\text{err}}$ & 1.7581                     & 2.6572                     & 3.8542                     & 4.8566                     & 5.8885                     & 6.9844                     & 5.7860                     & -2.8821                        \\ \hline
\textbf{Ortho}    & $\pp_{\text{err}}$   & 1.1984                     & 2.3292                     & 3.5494                     & 4.6377                     & 5.6633                     & 6.9237                     & 8.0005                     & -0.3784                        \\
                  & $\uu_{\text{err}}$   & \multicolumn{1}{l}{1.0354} & \multicolumn{1}{l}{1.8430} & \multicolumn{1}{l}{2.8787} & \multicolumn{1}{l}{3.9900} & \multicolumn{1}{l}{5.0345} & \multicolumn{1}{l}{6.0422} & \multicolumn{1}{l}{6.2590} & -2.7340                        \\
\textbf{}         & $\pp_{I,\text{err}}$ & 1.7581                     & 2.6572                     & 3.8542                     & 4.8566                     & 5.8885                     & 7.0225                     & 8.0005                     & -0.3784                        \\ \hline
\end{tabular}}
\caption{Test3: convergence rates on conforming mesh.}
\label{tab:DFN_conv}
\end{table}

Table \ref{tab:DFN_conv} shows the computed convergence rates of errors \eqref{eq:L2pressure} and \eqref{eq:superconvergence}, for all the three tested approaches.
In each sub-domain, the solution is a polynomial of degree $6$, such that, for $k \geq 6$ only errors related to floating point arithmetic computations are to be expected.

\begin{figure}[htbp]
	\centering
	\subfigure[\label{fig:condG_test9_DFN_abs_12_R0}] 
	{\includegraphics[width=.3\textwidth, height = .17\textheight]{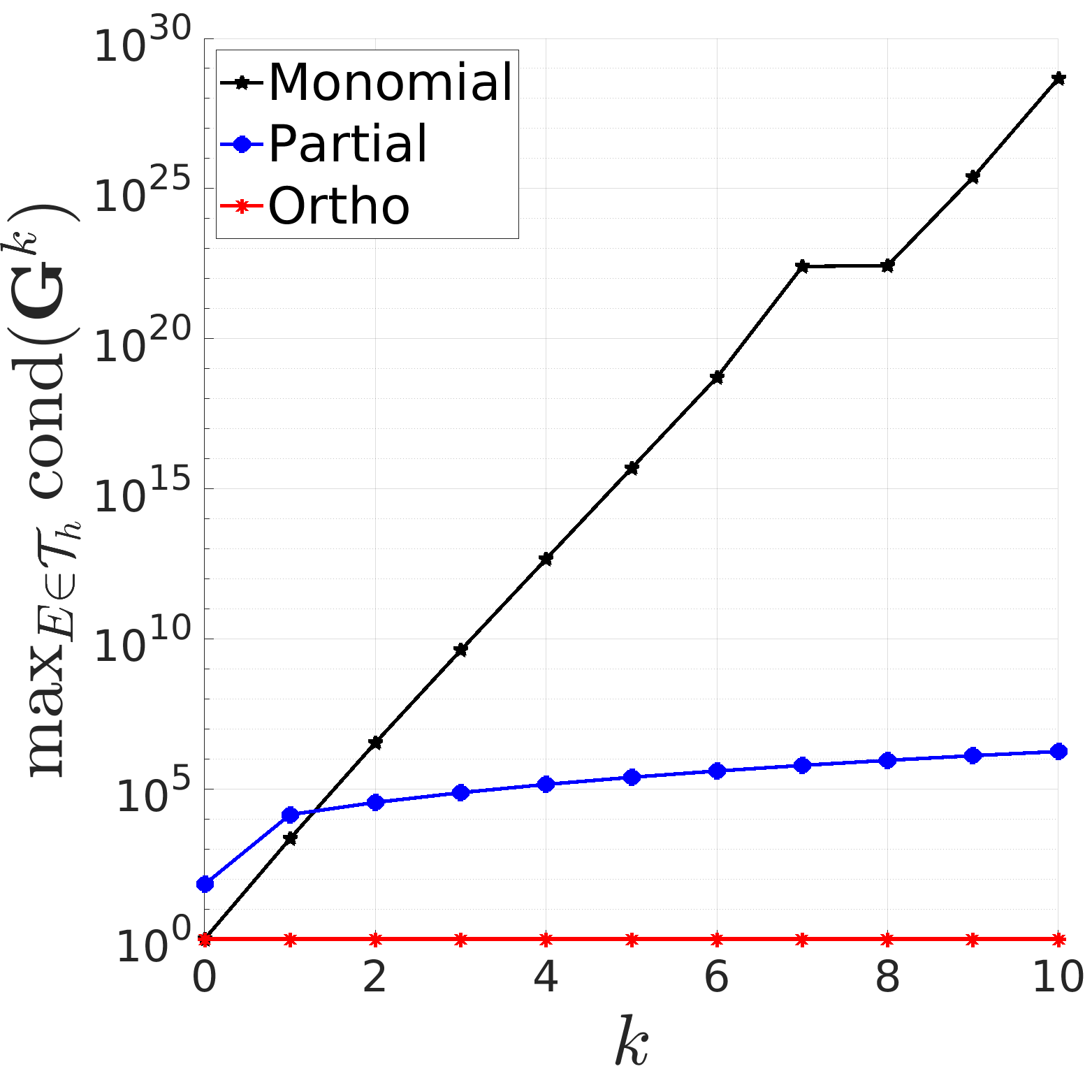}}
	\subfigure[\label{fig:condW_test9_DFN_abs_12_R0}] 
	{\includegraphics[width=.3\textwidth, height = .17\textheight]{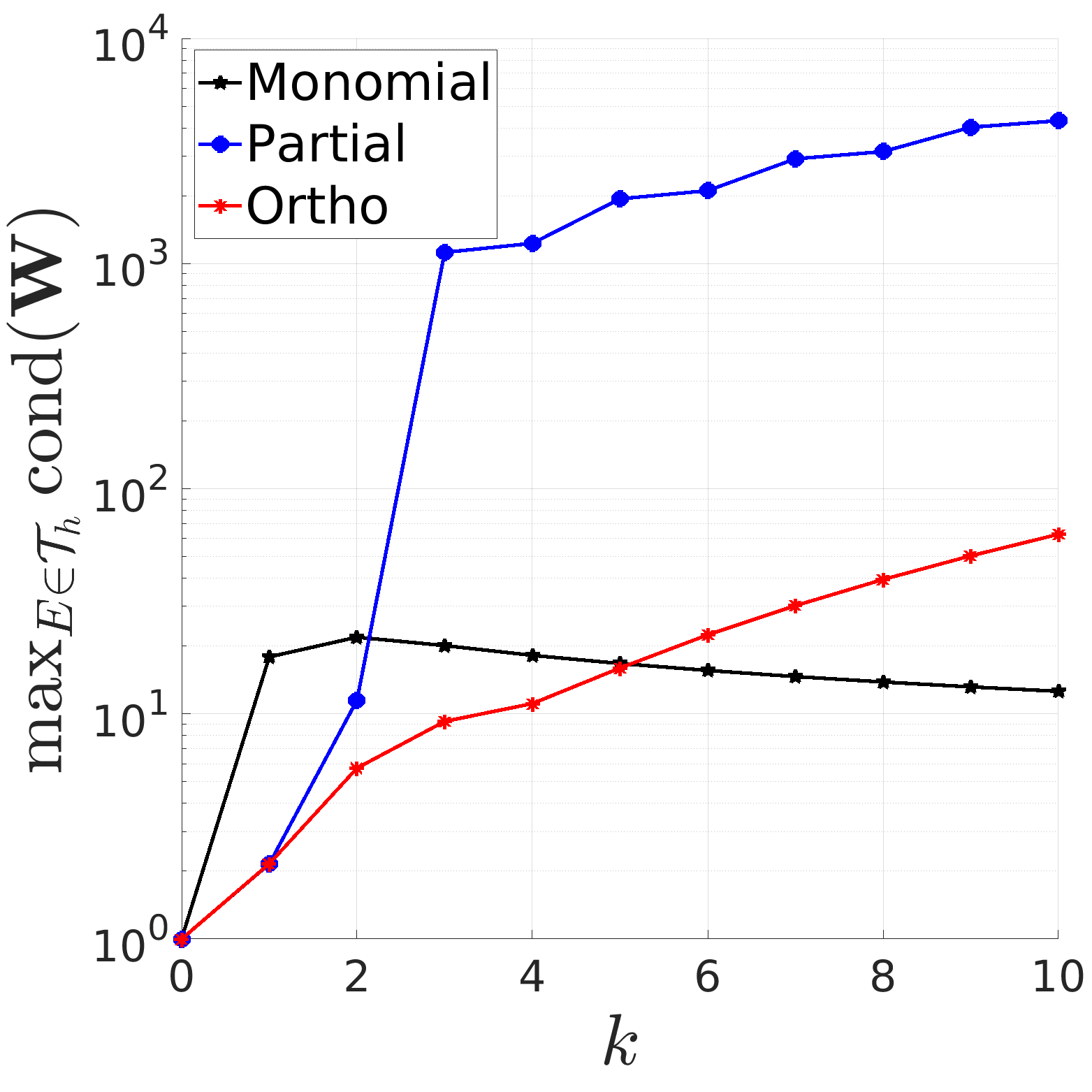}}
	\subfigure[\label{fig:condB_test9_DFN_abs_12_R0}] 
	{\includegraphics[width=.3\textwidth, height = .17\textheight]{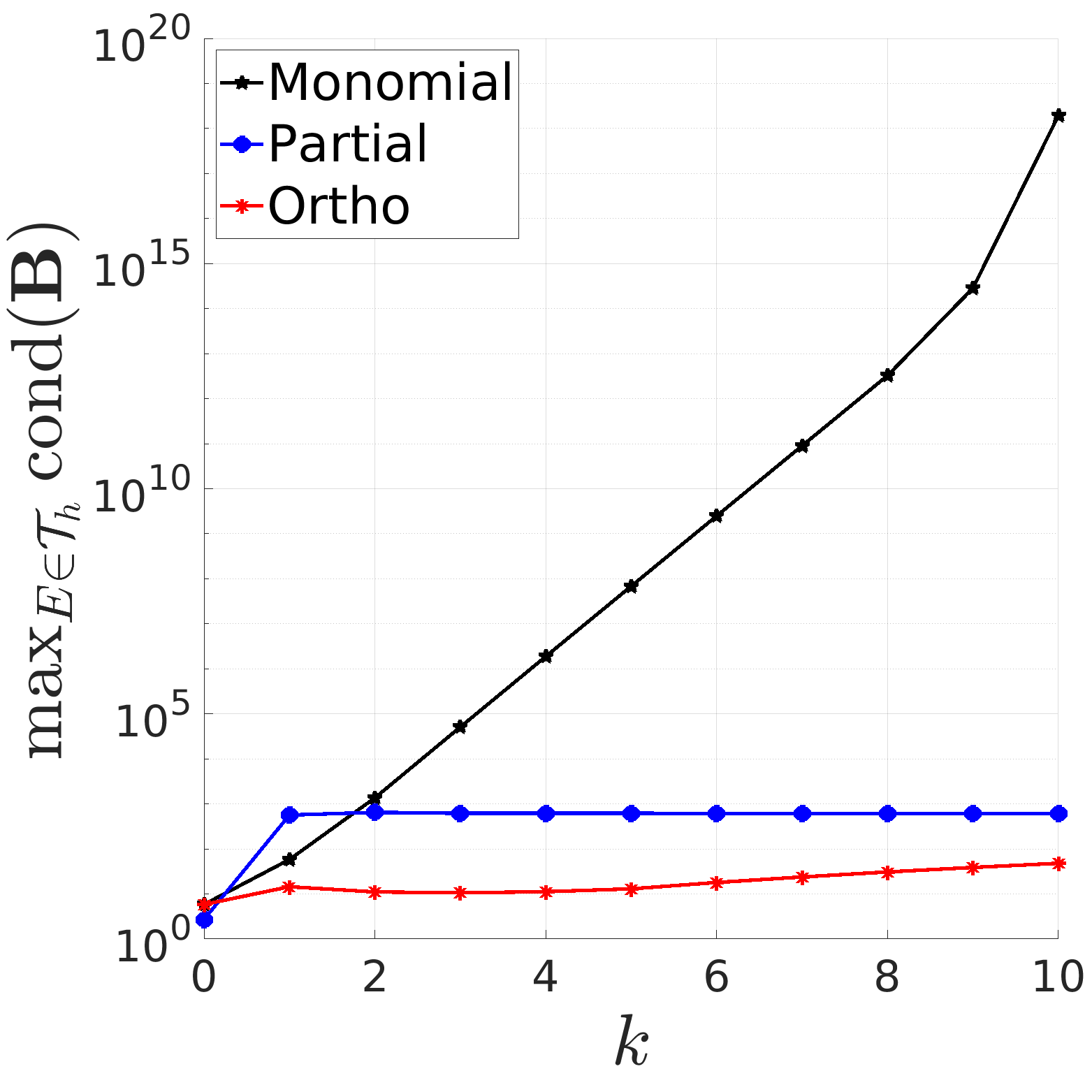}}
	\subfigure[\label{fig:Pi0k_test9_DFN_abs_12_R0}] 
	{\includegraphics[width=.3\textwidth, height = .17\textheight]{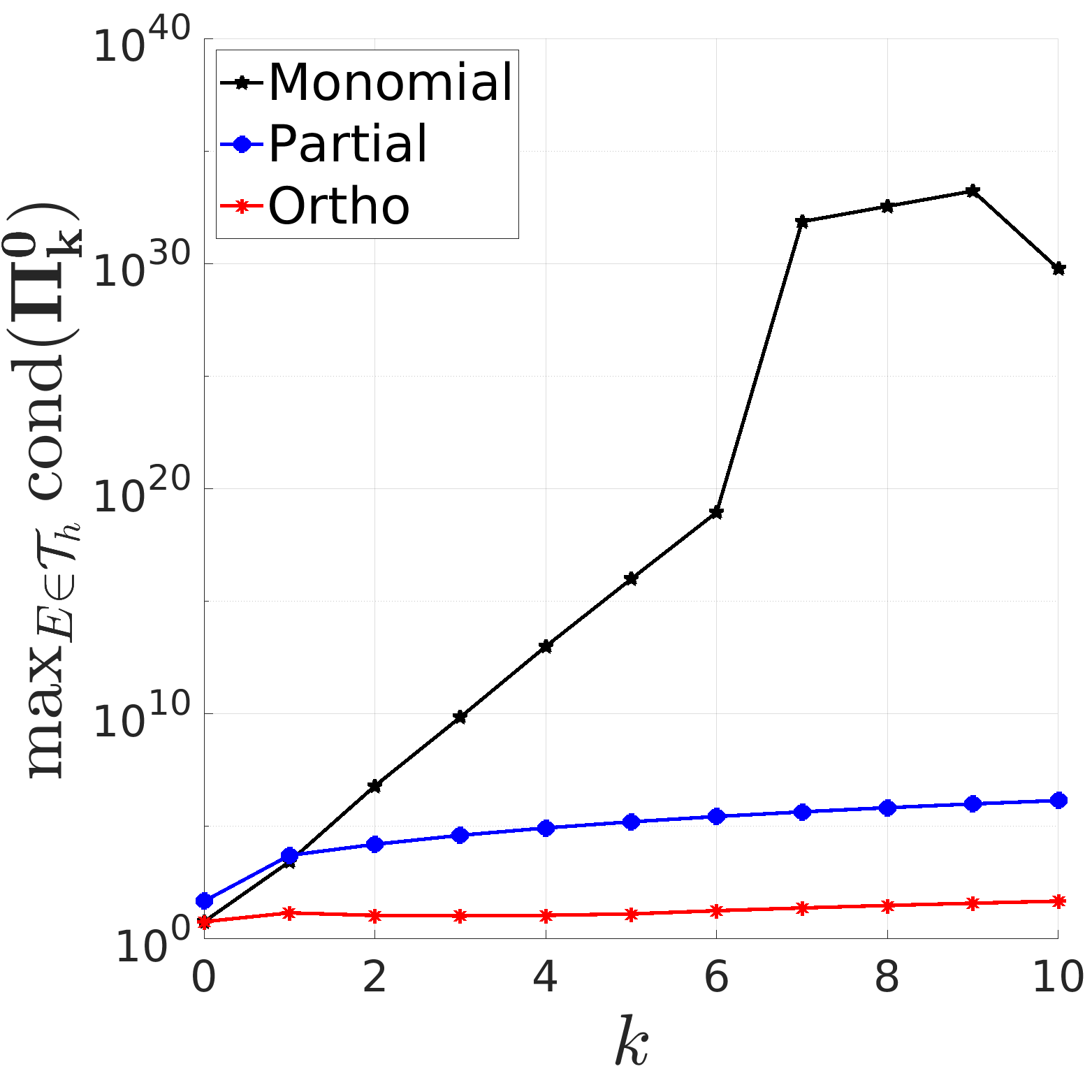}}
	\subfigure[\label{fig:condD_test9_DFN_abs_12_R0}]
	{\includegraphics[width=.3\textwidth, height = .17\textheight]{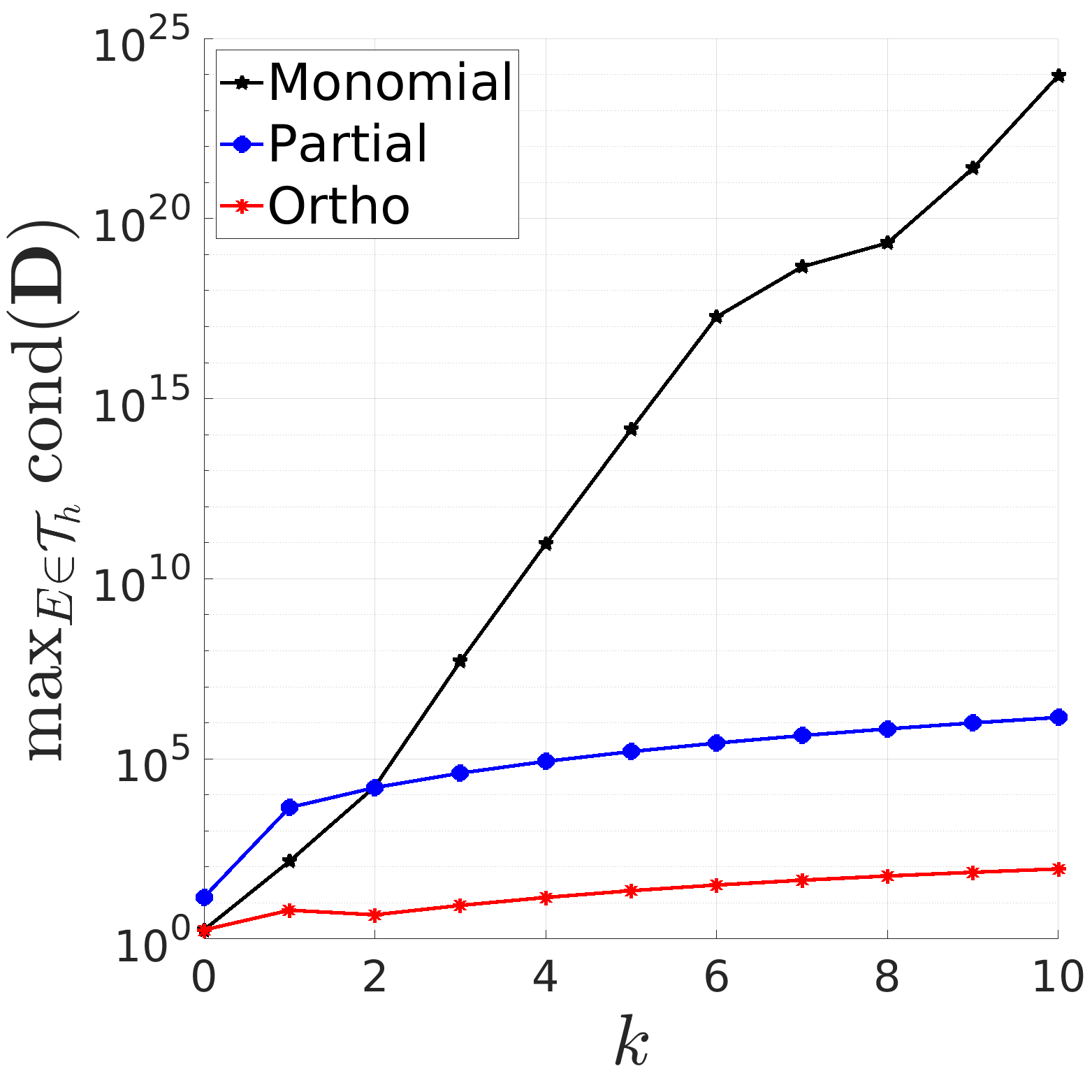}}
    \caption{Test3: Maximum condition number of local matrices among elements, at varying $k$. Coarsest mesh.}
    \label{fig:cond_first_conf}
\end{figure}

\begin{figure}[htbp]
	\centering
	\subfigure[\label{condG_test9_DFN_abs_12_R3}]
	{\includegraphics[width=.3\textwidth, height = .17\textheight]{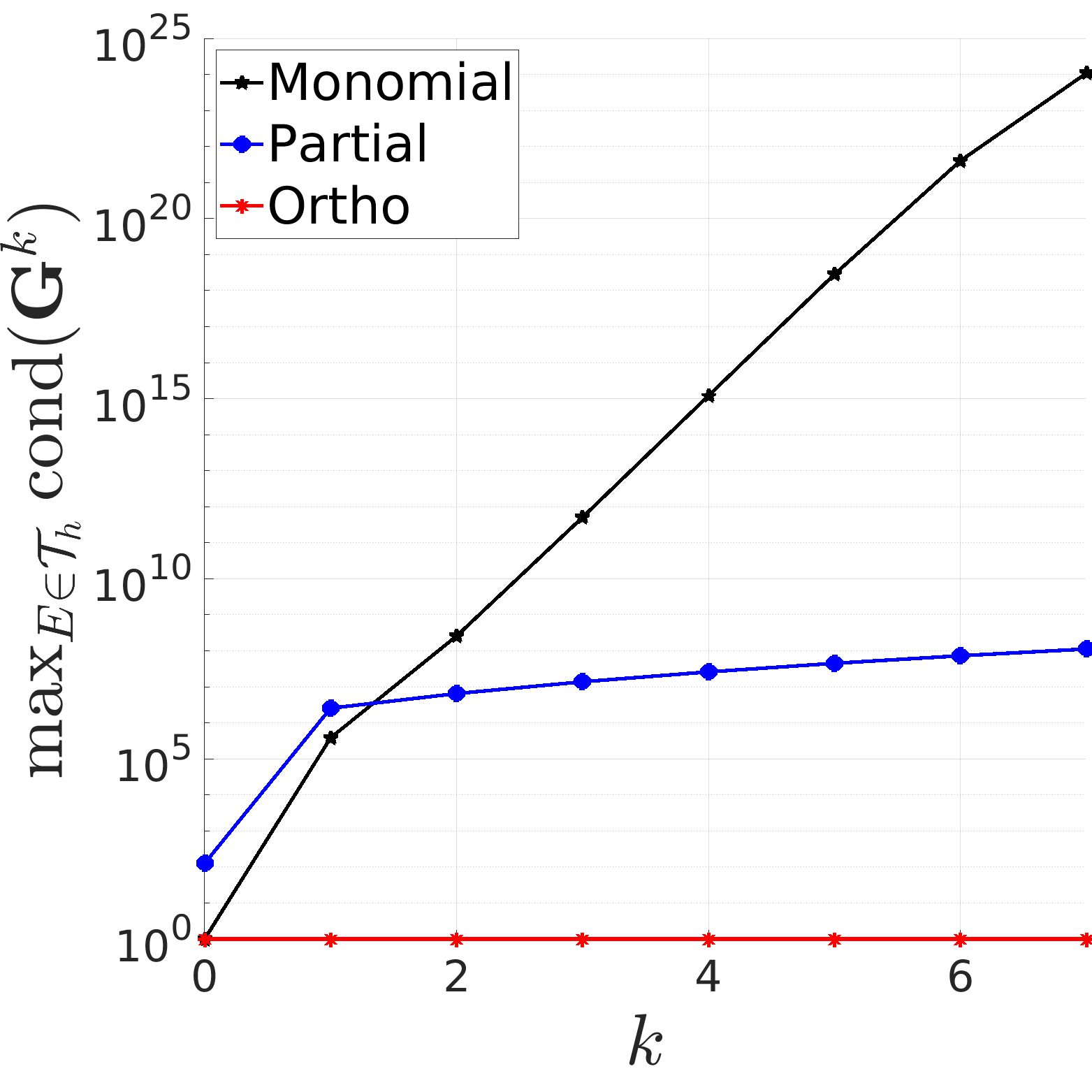}}
	\subfigure[\label{condW_test9_DFN_abs_12_R3}]
	{\includegraphics[width=.3\textwidth, height = .17\textheight]{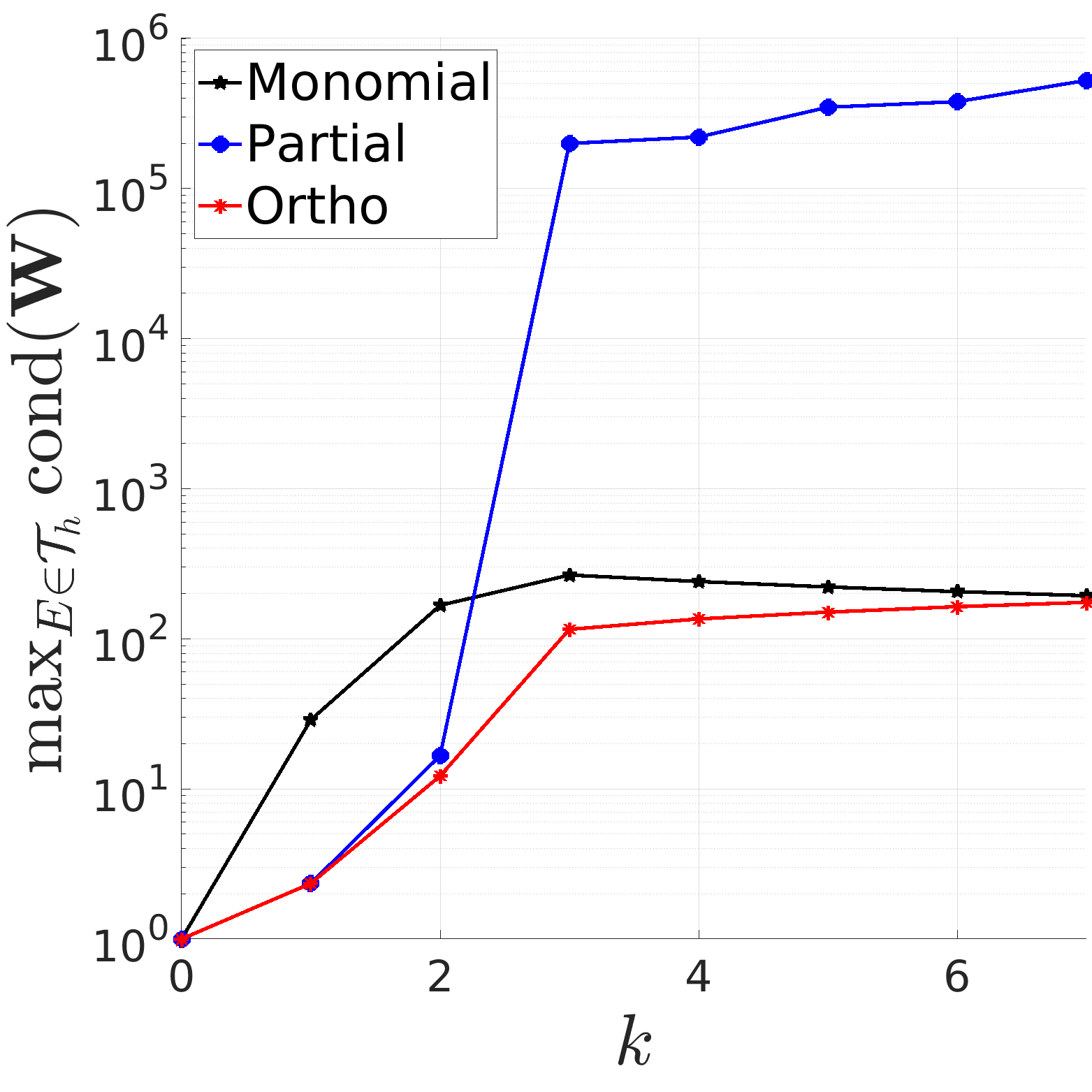}}
	\subfigure[\label{condB_test9_DFN_abs_12_R3}]
	{\includegraphics[width=.3\textwidth, height = .17\textheight]{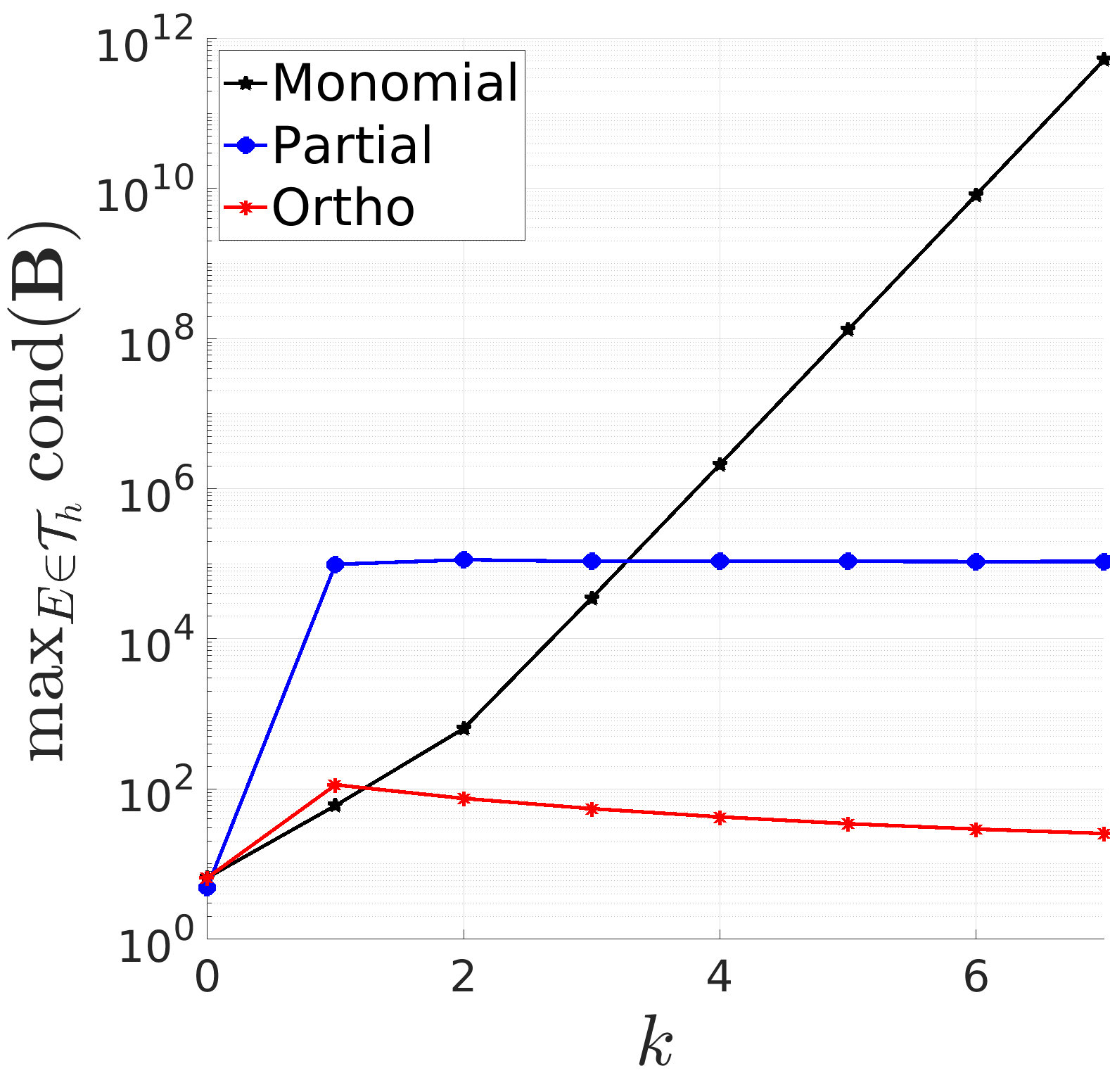}}
	\subfigure[\label{Pi0k_test9_DFN_abs_12_R3}]
	{\includegraphics[width=.3\textwidth, height = .17\textheight]{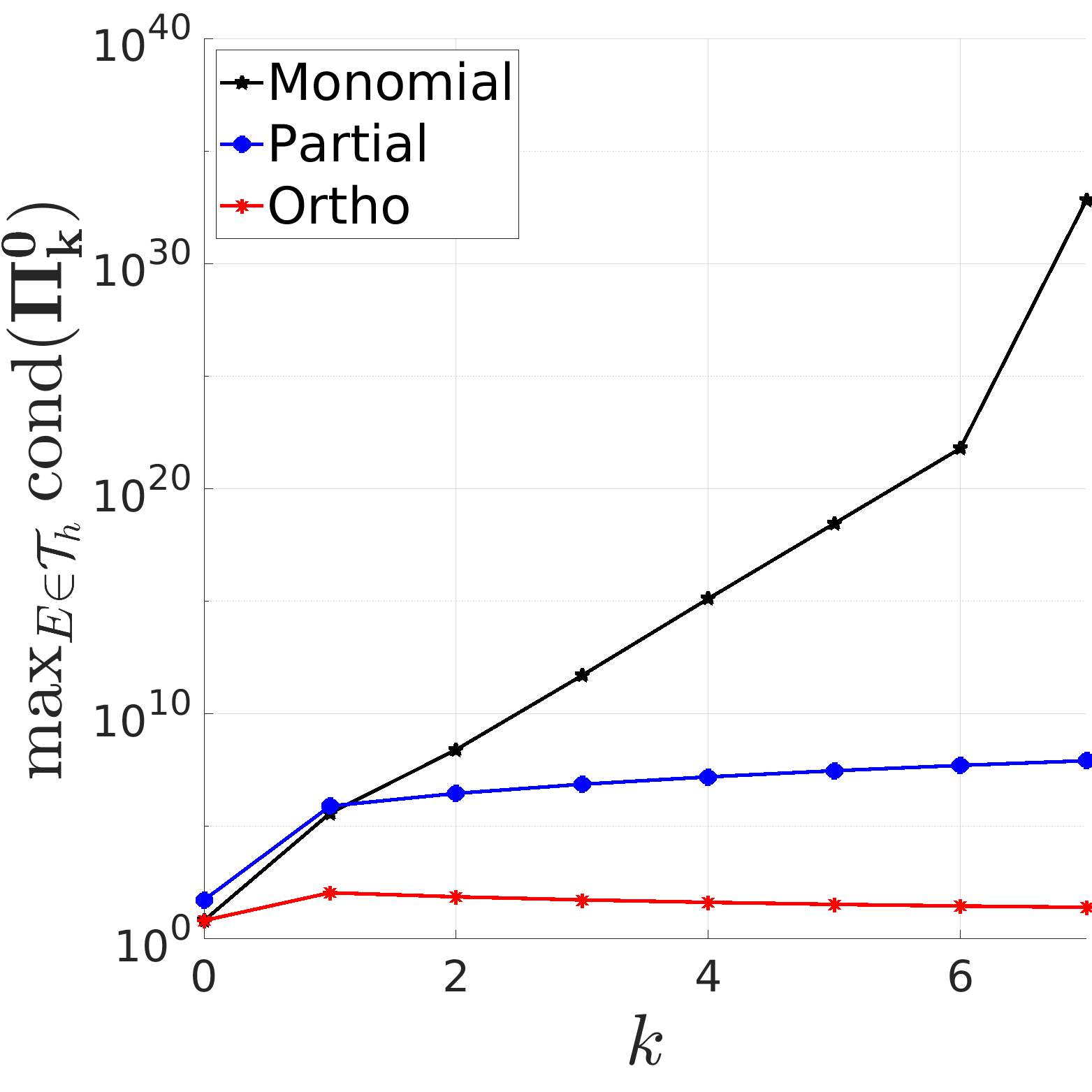}}
	\subfigure[\label{condD_test9_DFN_abs_12_R3}]
	{\includegraphics[width=.3\textwidth, height = .17\textheight]{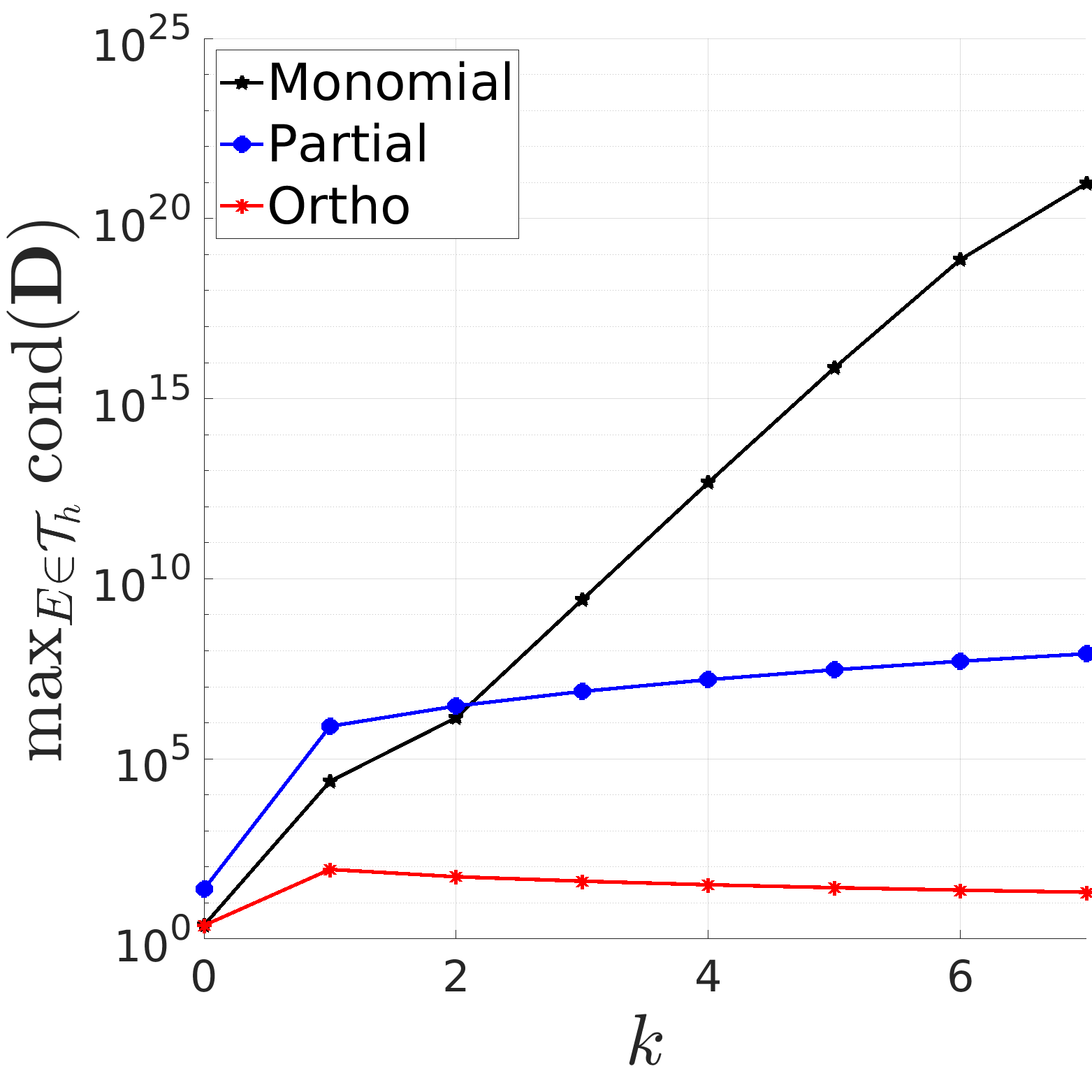}}
    \caption{Test3: Maximum condition number of local matrices among elements, at varying $k$. Finest mesh.}
    \label{fig:cond_last_conf}
\end{figure}

\begin{figure}[htbp]
	\centering
	\subfigure[\label{ErrorL2pressure_test9_DFN_abs_12_R0}]
	{\includegraphics[width=.3\textwidth, height = .17\textheight]{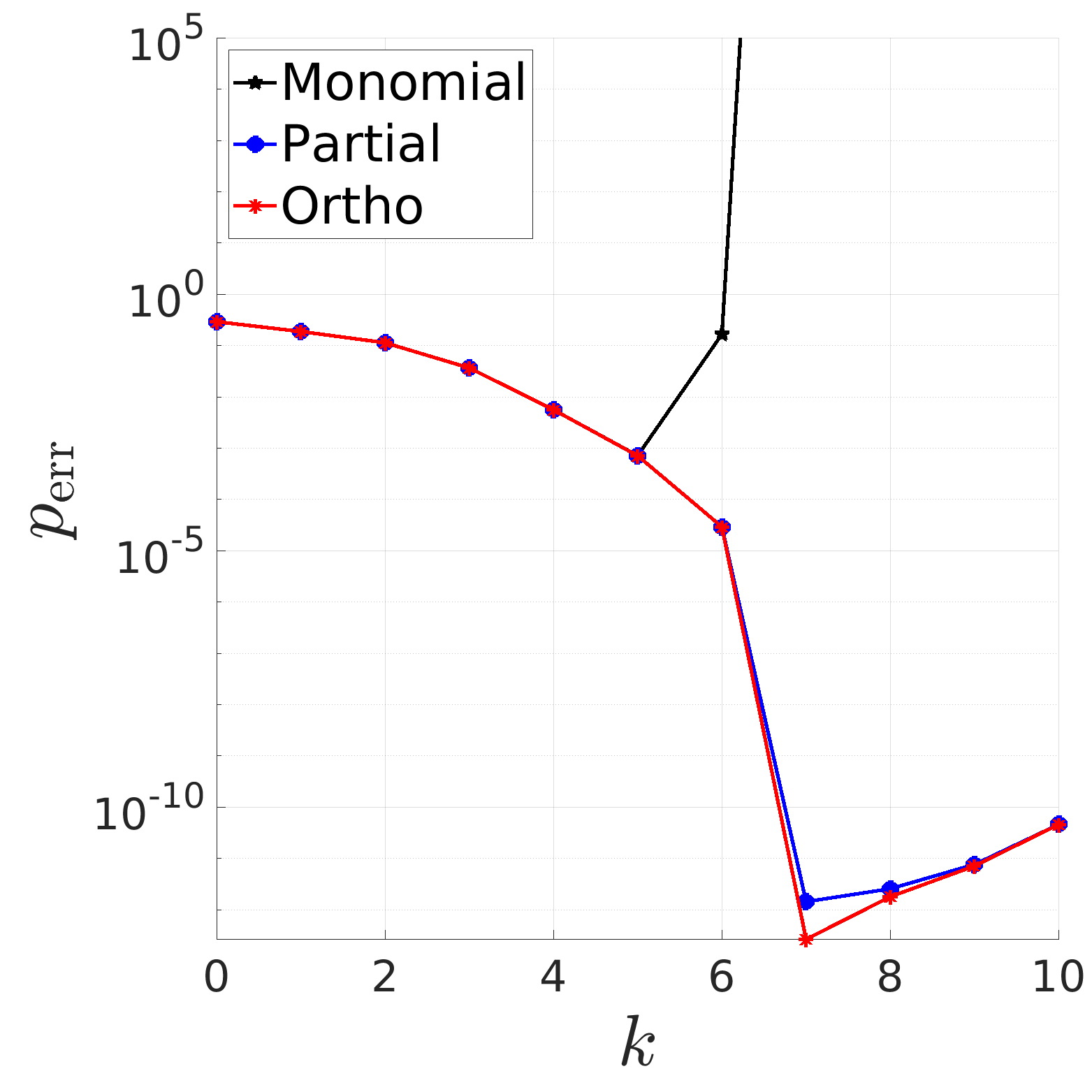}}
	\subfigure[\label{ErrorL2Velocity_test9_DFN_abs_12_R0}]
	{\includegraphics[width=.3\textwidth, height = .17\textheight]{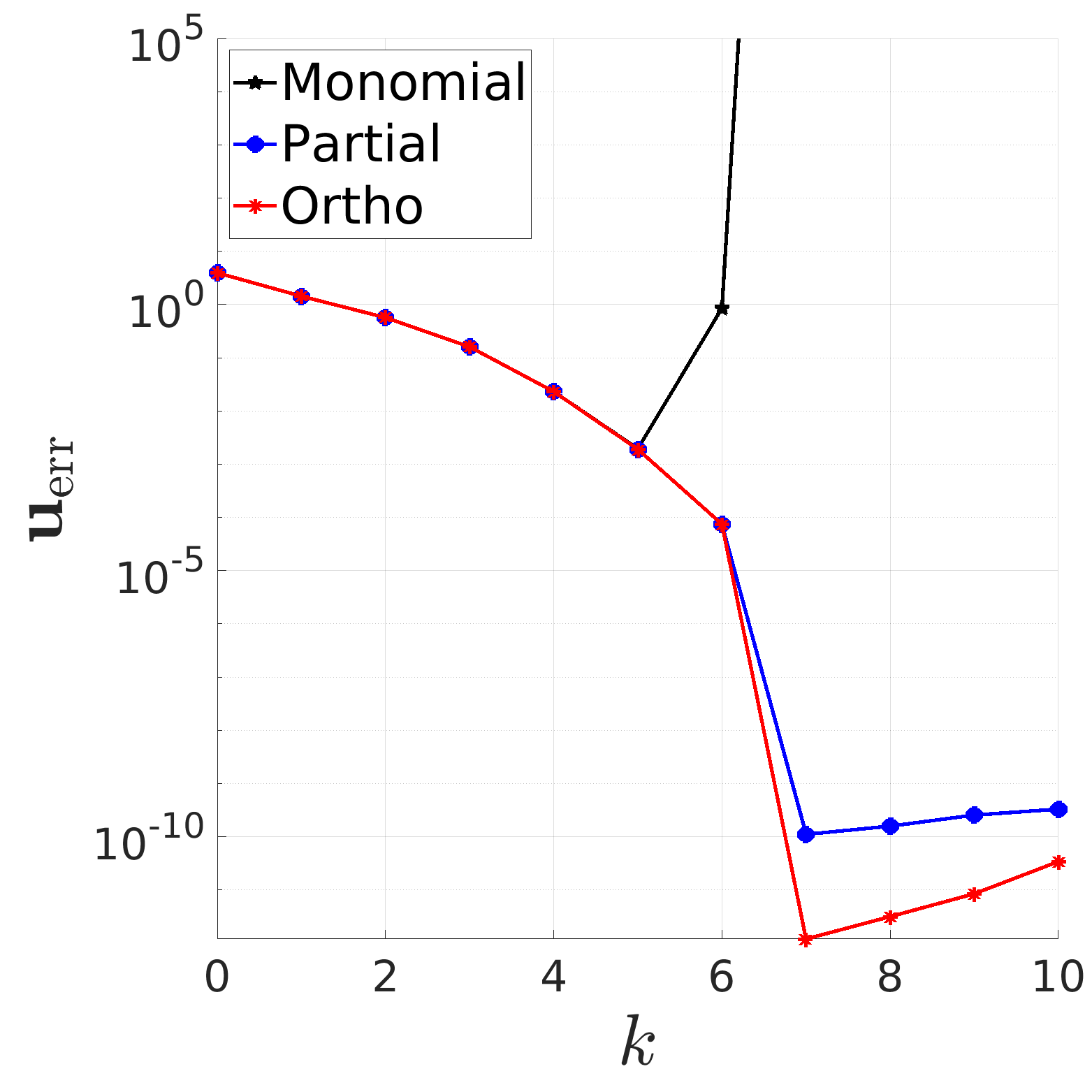}}
	\subfigure[\label{SuperConvergence_test9_DFN_abs_12_R0}]
	{\includegraphics[width=.3\textwidth, height = .17\textheight]{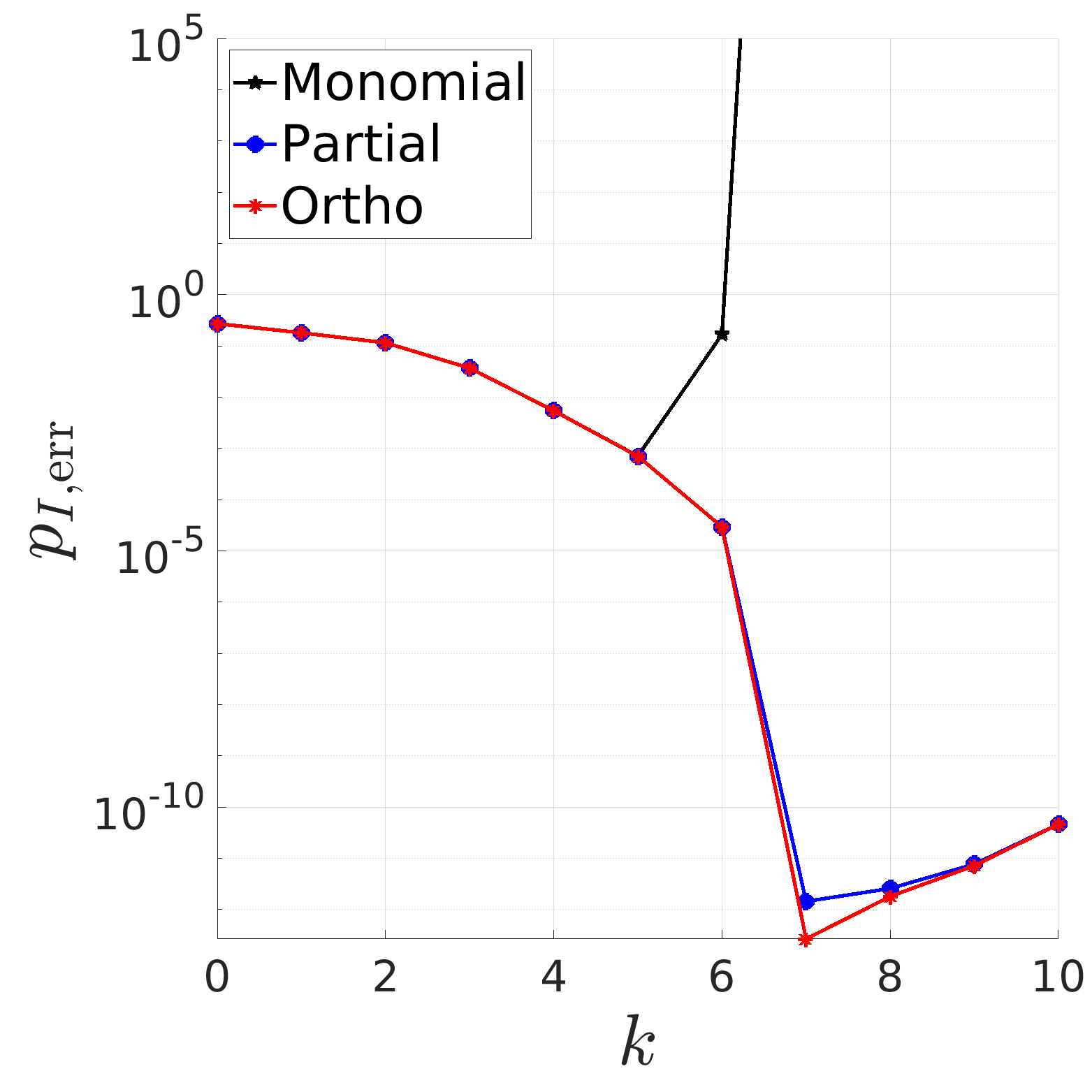}}
	\subfigure[\label{ErrorL2pressure_test9_DFN_abs_12_R1}]
	{\includegraphics[width=.3\textwidth, height = .17\textheight]{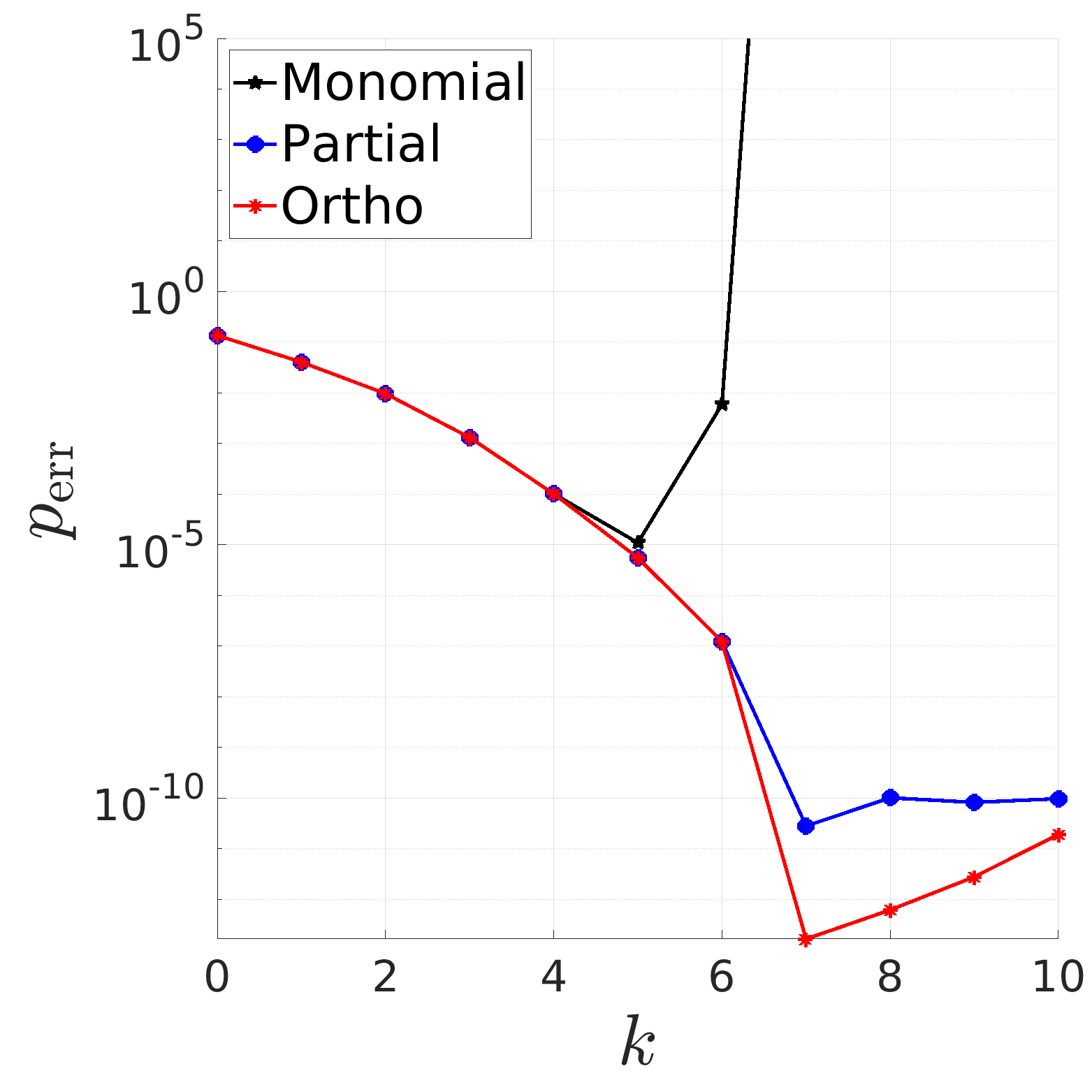}}
	\subfigure[\label{ErrorL2Velocity_test9_DFN_abs_12_R1}]
	{\includegraphics[width=.3\textwidth, height = .17\textheight]{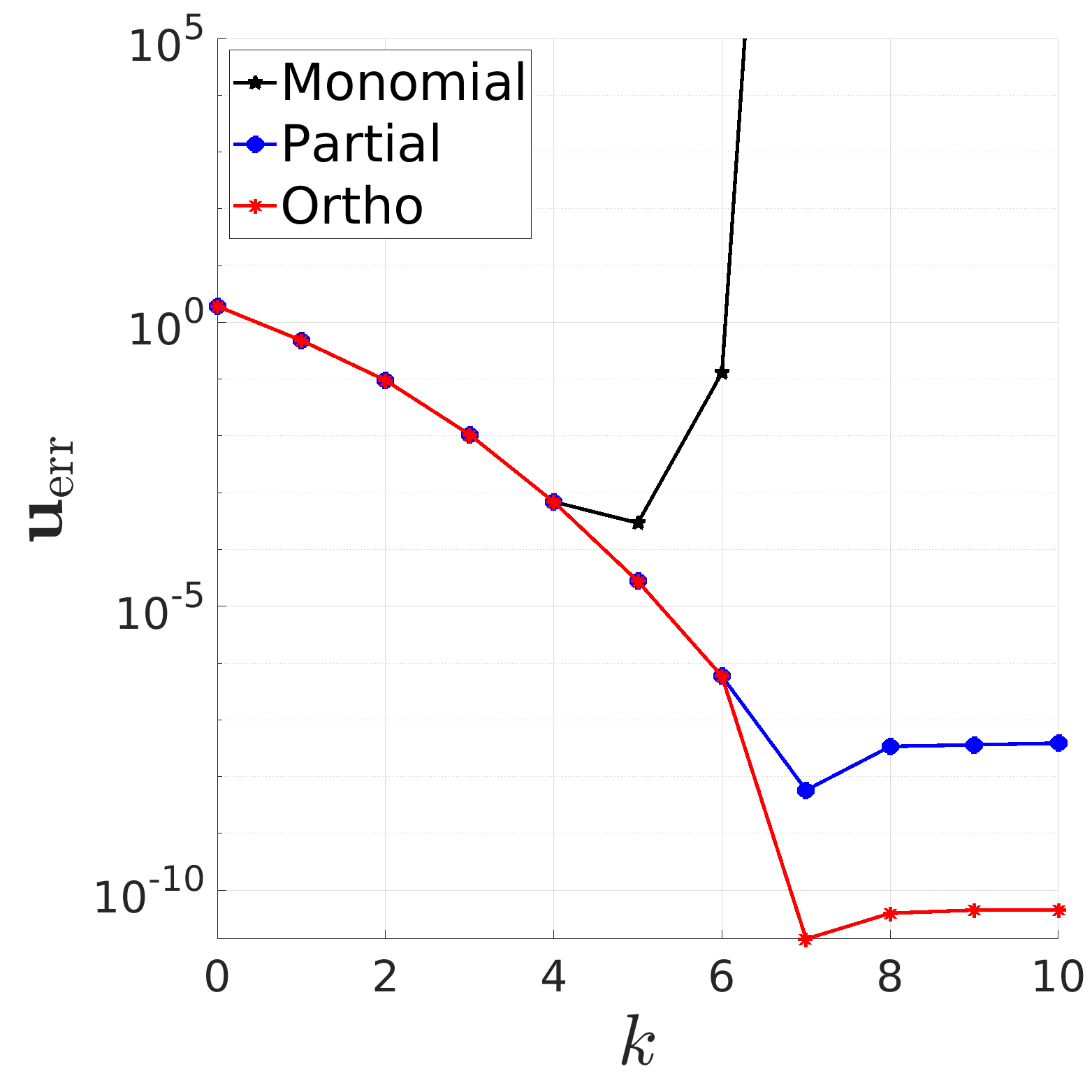}}
	\subfigure[\label{SuperConvergence_test9_DFN_abs_12_R1}]
	{\includegraphics[width=.3\textwidth, height = .17\textheight]{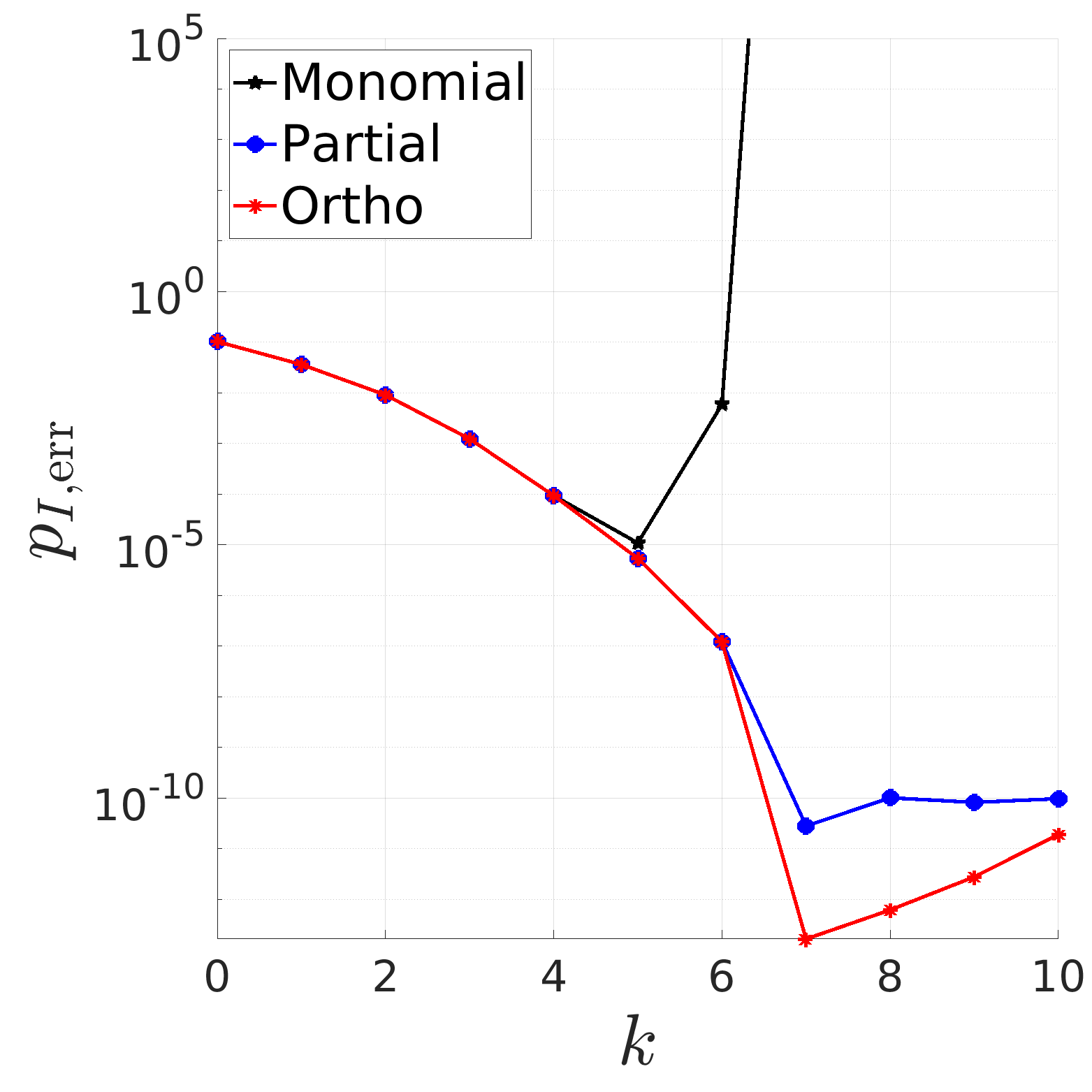}}
	\subfigure[\label{ErrorL2pressure_test9_DFN_abs_12_R2}]
	{\includegraphics[width=.3\textwidth, height = .17\textheight]{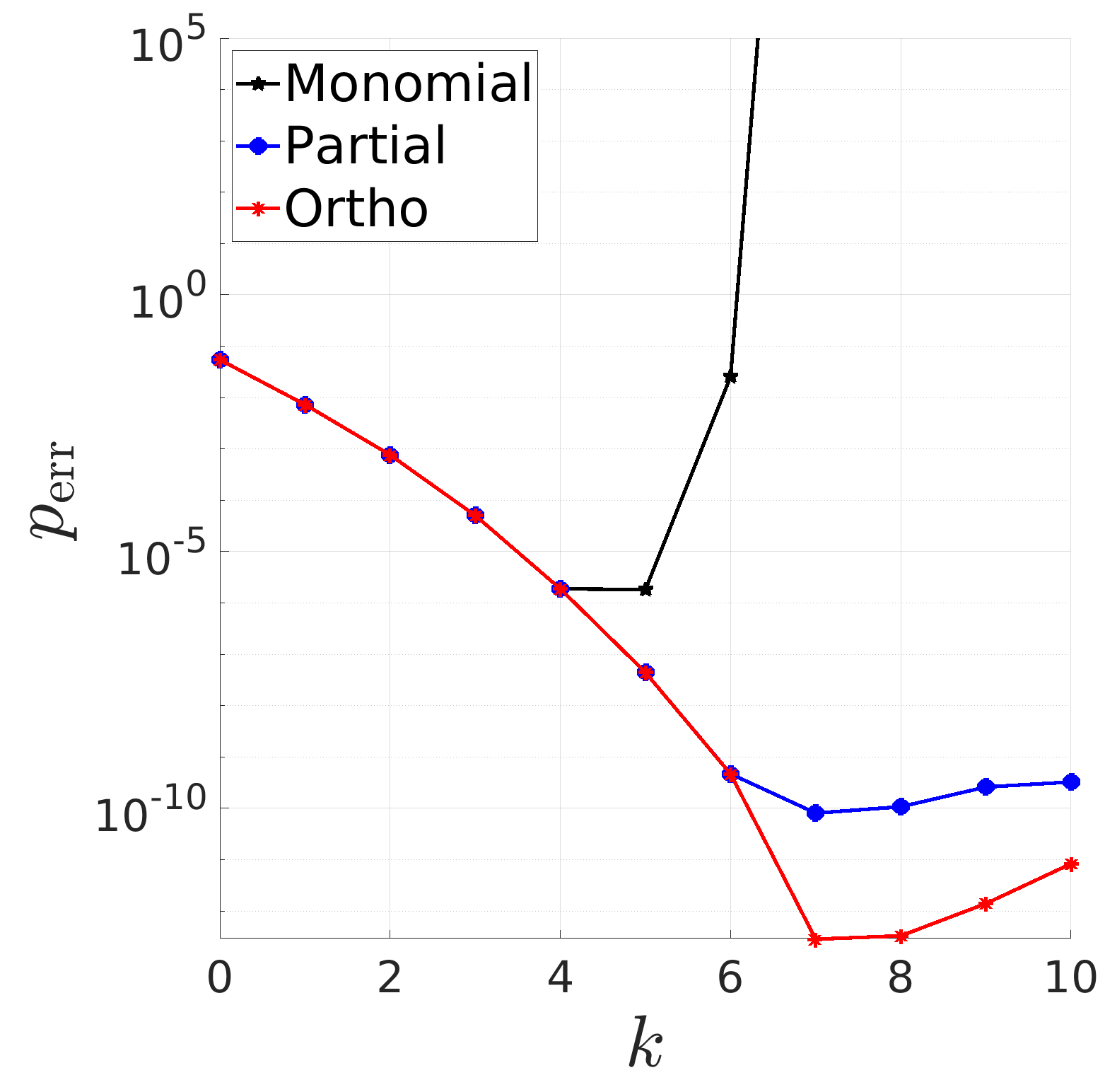}}
	\subfigure[\label{ErrorL2Velocity_test9_DFN_abs_12_R2}]
	{\includegraphics[width=.3\textwidth, height = .17\textheight]{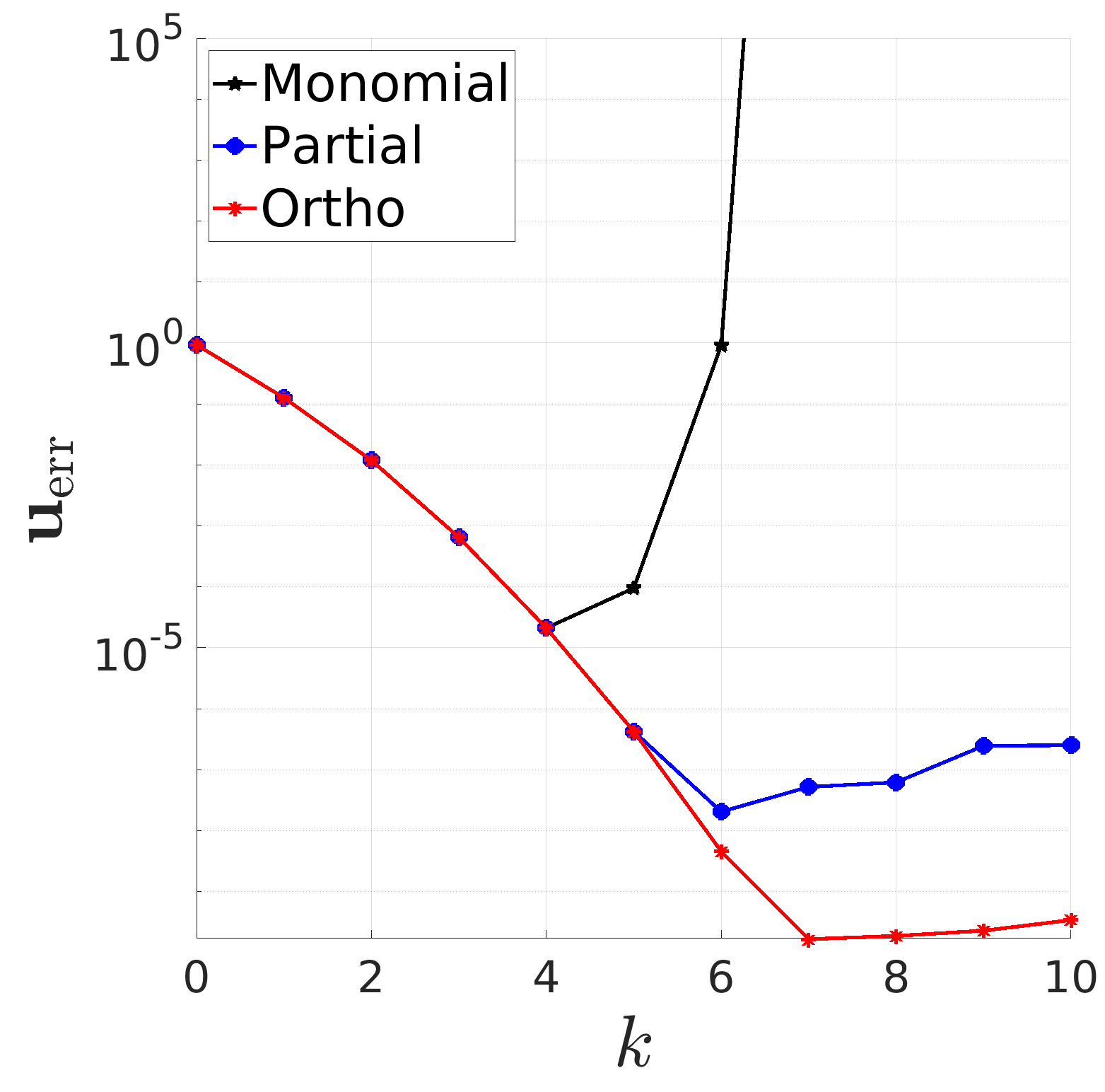}}
	\subfigure[\label{SuperConvergence_test9_DFN_abs_12_R2}]
	{\includegraphics[width=.3\textwidth, height = .17\textheight]{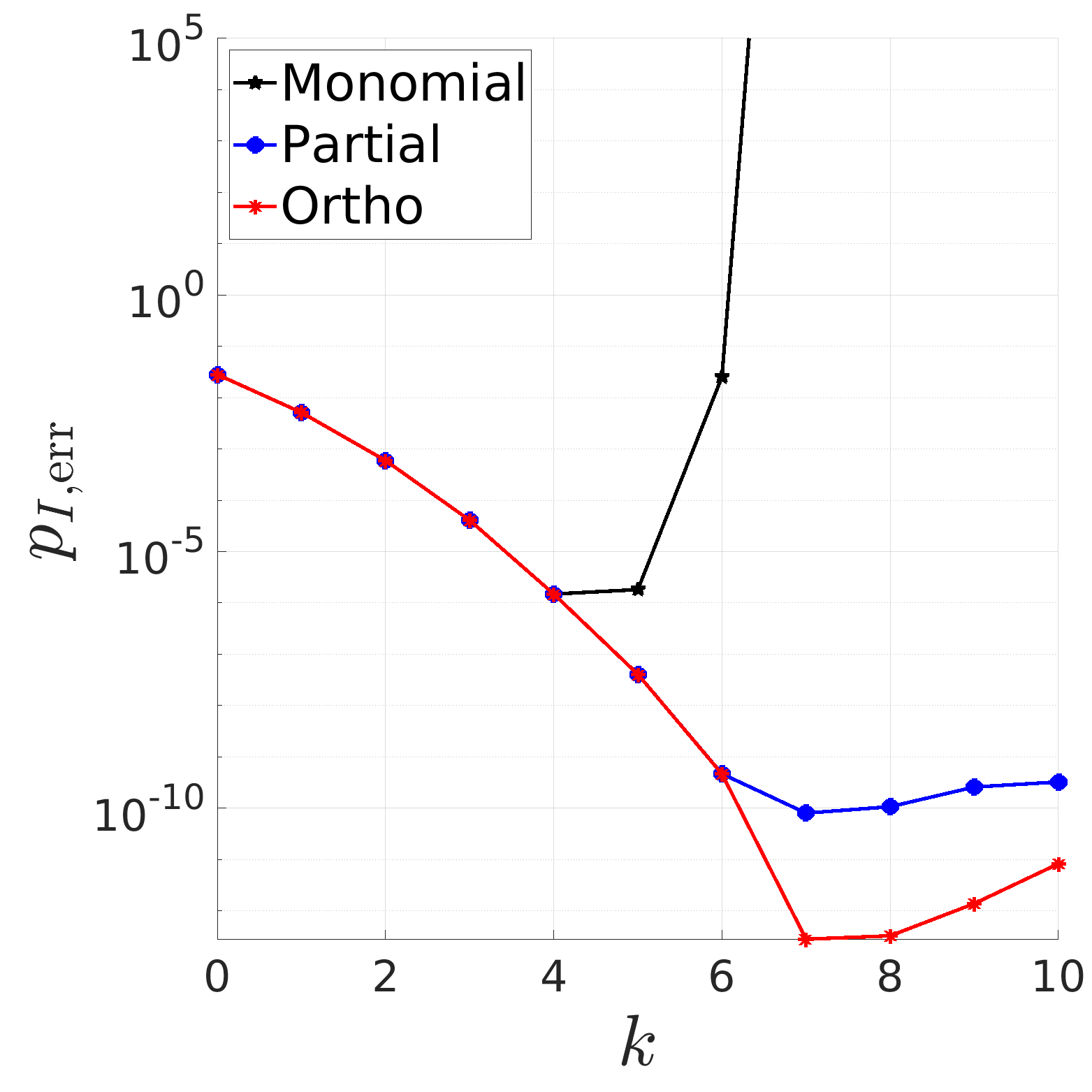}}
    \subfigure[\label{ErrorL2pressure_test9_DFN_abs_12_R3}]
	{\includegraphics[width=.3\textwidth, height = .17\textheight]{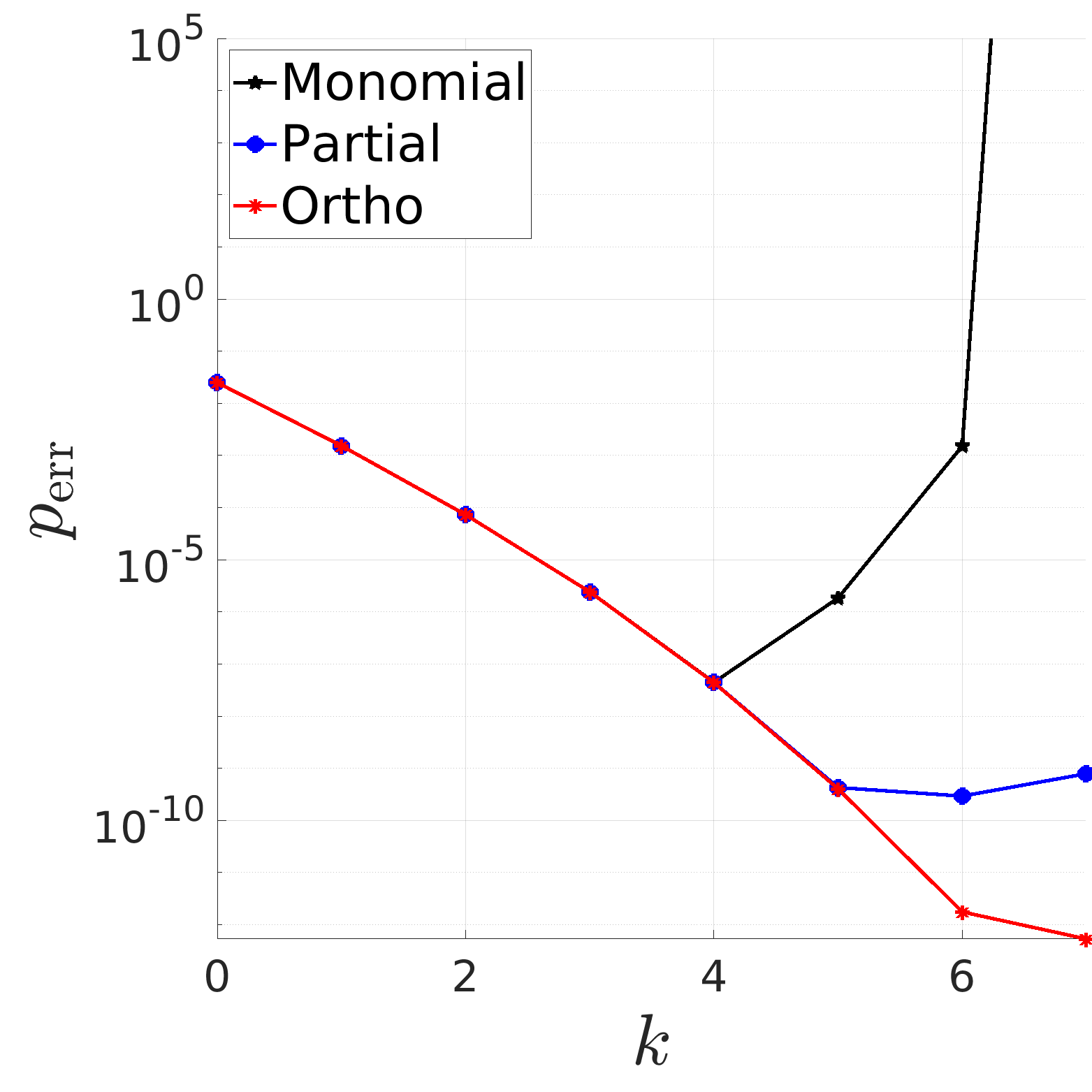}}
	\subfigure[\label{ErrorL2Velocity_test9_DFN_abs_12_R3}]
	{\includegraphics[width=.3\textwidth, height = .17\textheight]{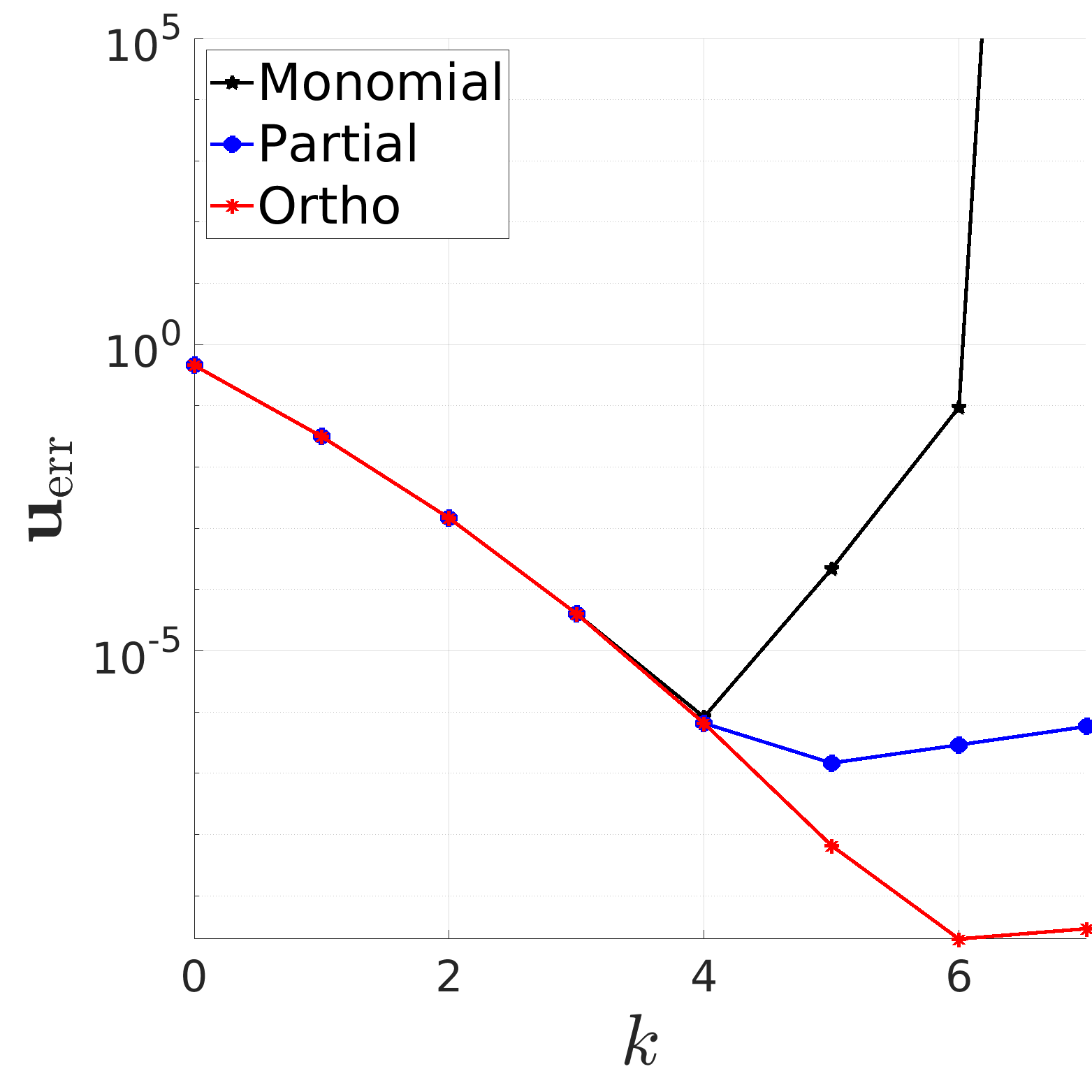}}
	\subfigure[\label{SuperConvergence_test9_DFN_abs_12_R3}]
	{\includegraphics[width=.3\textwidth, height = .17\textheight]{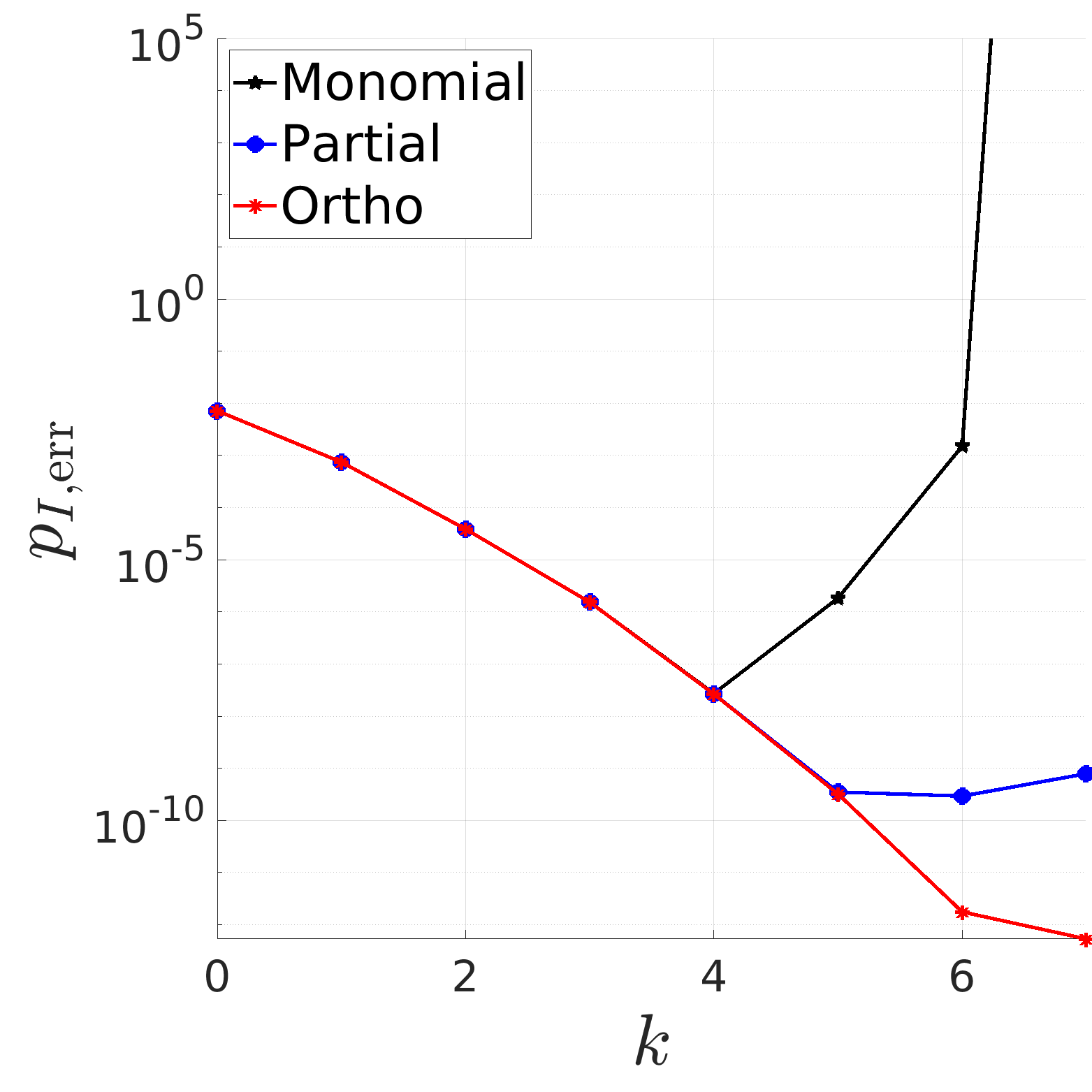}}
    \caption{Test3: Behaviour of errors \eqref{eq:L2pressure}, \eqref{eq:L2velocity} and \eqref{eq:superconvergence} at varying $k$ on conforming meshes. Each row represents a different refinement, from the coarsest mesh on top to the finest mesh at the bottom.}
    \label{fig:errors_thirdExp}
\end{figure}

Figures \ref{fig:cond_first_conf} and \ref{fig:cond_last_conf} report, as previously, the maximum condition number across mesh elements of the computed matrices, on the coarsest and finest considered meshes, showing the same behaviour as in the previous tests. 
Figure \ref{fig:errors_thirdExp} shows error convergence curves against polynomial accuracy $k$, for the four considered mesh refinement levels. Pictures on each row correspond to the same mesh. It can be seen that the curves obtained with the three approaches are almost indistinguishable for $k$ up to $4-5$. For higher values errors given by the monomial approach start growing rapidly, whereas errors given by partial and full orthonormal approaches still decrease up to stagnation due to finite precision arithmetic. 

\section{Conclusions}

In this paper, we presented a possible solution to cure the ill-conditioning of system matrix in the mixed formulation of the Virtual Element Method.

Since in the mixed formulation we need to introduce a discrete local space for both the pressure and the velocity variable, we have first introduced an orthonormal scalar-polynomial basis in the pressure space and then we have also orthonormalized the vector-polynomial basis used in the definition of the degrees of freedom related to the velocity variable.

Numerical experiments suggest that the introduction of orthonormal polynomial basis in both spaces allows to improve stability of mixed Virtual Elements for high order applications on distorted elements. 

It is worth to mention that the methods here suggested to build orthonormal polynomial bases improve the conditioning of the Vandermonde matrix defined with respect to the quadrature formulas in the interior of each element $E$. However, in general, this does not guarantee an improvement in the conditioning of the Vandermonde matrix defined with respect to quadrature formulas on the boundary of the elements. Nonetheless, this appears to be sufficient to recover optimal convergence trends in the considered cases.  

\section*{Acknowledgments}
\noindent
The author S.B. kindly acknowledges partial financial support provided by PRIN project “Advanced polyhedral discretisations of heterogeneous PDEs for multiphysics problems” (No. 20204LN5N5\_003) and by PNRR M4C2 project of CN00000013 National Centre for HPC, Big Data and Quantum Computing (HPC) (CUP: E13C22000990001). 
The author S.S. kindly acknowledges partial financial support provided by INdAM-GNCS through project “Sviluppo ed analisi di Metodi agli Elementi Virtuali per processi accoppiati su geometrie complesse” and that this  publication is part of the project NODES which has received funding from the MUR-M4C2 1.5 of PNRR with grant agreement no. ECS00000036. 
The author G.T. kindly acknowledges financial support provided by the MIUR programme ``Programma Operativo Nazionale Ricerca e Innovazione 2014 - 2020'' 
~ (CUP: E11B21006490005). Computational resources are partially supported by SmartData@polito.